\documentclass[12pt]{amsart}
\usepackage{amsbsy}
\usepackage{graphicx,epsfig,subfigure,psfrag}
\begin{document}

\newtheorem{theorem}{Theorem}
\newtheorem{proposition}{Proposition}
\newtheorem{lemma}{Lemma}
\newtheorem{corollary}{Corollary}
\newtheorem{definition}{Definition}
\newtheorem{remark}{Remark}
\newcommand{\beq}{\begin{equation}}
\newcommand{\eeq}{\end{equation}}
\numberwithin{equation}{section}
\numberwithin{theorem}{section}
\numberwithin{proposition}{section}
\numberwithin{lemma}{section}
\numberwithin{corollary}{section}
\numberwithin{definition}{section}
\numberwithin{remark}{section}
\newcommand{\re}{{\mathbb R}}
\newcommand{\n}{\nabla}
\newcommand{\ren}{{\mathbb R}^N}
\newcommand{\iy}{\infty}
\newcommand{\pa}{\partial}
\newcommand{\fp}{\noindent}
\newcommand{\ms}{\medskip\vskip-.1cm}
\newcommand{\mpb}{\medskip}
\newcommand{\BB}{{\bf B}}
\newcommand{\Am}{{\bf A}_{2m}}
\renewcommand{\a}{\alpha}
\renewcommand{\b}{\beta}
\newcommand{\g}{\gamma}
\newcommand{\G}{\Gamma}
\renewcommand{\d}{\delta}
\newcommand{\D}{\Delta}
\newcommand{\e}{\varepsilon}
\newcommand{\var}{\varphi}
\renewcommand{\l}{\lambda}
\renewcommand{\o}{\omega}
\renewcommand{\O}{\Omega}
\newcommand{\s}{\sigma}
\renewcommand{\t}{\tau}
\renewcommand{\th}{\theta}
\newcommand{\z}{\zeta}
\newcommand{\wx}{\widetilde x}
\newcommand{\wt}{\widetilde t}
\newcommand{\noi}{\noindent}
\newcommand{\inA}{\quad \mbox{in} \quad \ren \times \re_+}
\newcommand{\inB}{\quad \mbox{in} \quad}
\newcommand{\inC}{\quad \mbox{in} \quad \re \times \re_+}
\newcommand{\inD}{\quad \mbox{in} \quad \re}
\newcommand{\forA}{\quad \mbox{for} \quad}
\newcommand{\whereA}{,\quad \mbox{where} \quad}
\newcommand{\asA}{\quad \mbox{as} \quad}
\newcommand{\andA}{\quad \mbox{and} \quad}
\newcommand{\ssk}{\smallskip}
\newcommand{\LongA}{\quad \Longrightarrow \quad}
\def\com#1{\fbox{\parbox{6in}{\texttt{#1}}}}

\title 
 {\bf Variational approach to complicated  similarity solutions
of higher-order
 nonlinear
 PDEs. I}

\author{
V.A.~Galaktionov, E.~Mitidieri,  and S.I.~Pohozaev}

\address{Department of Mathematical Sciences, University of Bath,
 Bath BA2 7AY, UK}
\email{vag@maths.bath.ac.uk}

\address{Dipartimento di Scienze Matematiche,
Universit\`a degli Studi  di Trieste, Piazzale Europa 1, 34100
Trieste, ITALY} \email{mitidier@units.it}

\address{Steklov Mathematical Institute,
 Gubkina St. 8, 119991 Moscow, RUSSIA}
\email{pokhozhaev@mi.ras.ru}



\date{\today}




\begin{abstract}

   The Cauchy problem for  $(x,t) \in \ren
  \times \re_+$ for three higher-order degenerate  quasilinear PDEs, as
  basic models,
 $$
  \begin{matrix}
 u_t = (-1)^{m+1}\D^m(|u|^n u)+|u|^n u,\qquad\,\quad\ssk\\
 u_{tt} =(-1)^{m+1} \D^m(|u|^n u)+|u|^n u,\qquad\quad\ssk\\
 u_{t}= (-1)^{m+1}[\D^m(|u|^n u)]_{x_1}+(|u|^n
 u)_{x_1},
 \end{matrix}
 $$
  where $n>0$ is a fixed exponent and  $\D^m$ is the $(m \ge 2)$th iteration of the
  Laplacian,
  is studied. This
  diverse  class of degenerate PDEs embraces equations of different three
  types:
  parabolic, hyperbolic, and nonlinear dispersion.
  Such degenerate evolution equations  from various
 applications in mechanics and physics
  admitting
  compactly supported and blow-up solutions
  attracted
 attention of the mathematicians since
 the  1970-80s.

In the present paper, some general local, global, and blow-up
features of such PDEs on the basis of
 construction of their {\em blow-up} similarity and
{\em travelling wave} solutions are revealed.
 Blow-up, i.e., nonexistence of global in time solutions
 is proved by various methods.
 In particular, for
$m=2$ and $m=3$, such similarity patterns lead to the following
 semilinear fourth- and sixth-order elliptic PDEs with non-coercive
 operators and
  non-Lipschitz nonlinearities:
  \beq
  \label{a1}
   \mbox{$
  - \D^2 F +F - |F|^{-\frac{n}{n+1}} F=0 \andA \D^3 F + F- |F|^{-\frac n{n+1}}
  F=0
   \inA,
    $}
    \eeq
  which were not addressed  before in the
 mathematical literature. The goal is, using a variety of analytic variational,
 qualitative, and, often, numerical methods,
 to justify that  equations (\ref{a1}) admit an infinite at least countable set  of
countable families of compactly supported solutions that are
oscillatory near finite interfaces. In a whole, this solution set
 exhibit typical features of being of a chaotic structure.

 The present paper is an earlier extended version of
\cite{GMPSob}. In particular, 
here we pay  more attention to some category/fibering aspects of
critical values and points, as well as to sixth and higher-order
equations with $m \ge 3$, while in \cite{GMPSob} the case $m=2$
was under maximally detailed scrutiny.


\end{abstract}

\maketitle








\setcounter{equation}{0}
\section{Introduction:
 higher-order models and blow-up/compacton
solutions}
 \label{Sect1}
  \setcounter{equation}{0}




 \subsection{Three types of nonlinear PDEs under consideration}

 We describe   common local and
global properties of weak compactly solutions of classes of
nonlinear partial differential equations (PDEs) of parabolic,
hyperbolic, and nonlinear dispersion type:

 \ssk

 {\bf (I)} $2m$th-order {\em quasilinear parabolic equations with regional
blow-up},

 \ssk

{\bf (II)} $2m$th-order {\em quasilinear hyperbolic equations with
regional blow-up}, and

\ssk

{\bf (III)} ($2m$+1)th-order {\em nonlinear dispersion equations
 with compactons}.

\ssk

 Therefore,  we plan to study  some  common local and
global properties of weak compactly solutions of three classes of
quasilinear partial differential equations (PDEs) of parabolic,
hyperbolic, and nonlinear dispersion type, which, in general, look
like having nothing in common.
 Studying and better understanding of nonlinear degenerate PDEs of
higher-order including a new class of less developed {\em
nonlinear dispersion equations} (NDEs) from compacton theory are
striking features of modern general PDE theory  at the beginning
of  the twenty first century. It is worth
 noting and
realizing that several  key theoretical demands of modern
mathematics  are already associated and connected  with some
common local and global features of nonlinear evolution PDEs of
different types and orders, including higher-order parabolic,
hyperbolic, nonlinear dispersion, and others as typical
representatives.

Regardless the great progress of PDE theory achieved in the
twentieth century for many key classes of nonlinear equations
\cite{BrBra}, the transition process to higher-order degenerate
PDEs with more and more complicated non-monotone, non-potential,
and non-symmetric nonlinear operators will require different new
methods and mathematical approaches. Moreover, it seems that, for
some types of such nonlinear higher-order problems, the entirely
rigorous ``exhaustive" goal of developing a complete description
of solutions, their properties, functional settings of problems of
interest, etc., cannot be achieved in principle, in view of an
extremal variety of singular, bifurcation, and branching phenomena
that are contained in such multi-dimensional evolution.
 In
many cases, the main results should be extracted by a combination
of methods of various analytic, qualitative, and numerical
origins.
 In
many cases, the main results should be extracted by a combination
of methods of various analytic, qualitative, and numerical
origins\footnote{This is not a novelty in modern mathematics,
where several fundamental rigorous results have been already
justified with the aid of  hard, refined, and reliable  numerical
experiments; {\em q.v.} e.g., Tucker's proof of existence of a
robust strange attractor for the 3D Lorenz system \cite{Tuck02}.}.


In the present paper, we deal with  complicated pattern sets,
where,  for the elliptic problems in $\ren$ and even for the
corresponding   one-dimensional ODE reductions,  using the
proposed {\em analytic-numerical} approaches is necessary and
unavoidable.
 As a first illustration of such features, let us
 mention that, according to our current experience, for such classes of second-order $C^1$ variational
 problems,
 $$
 \fbox{$
 \begin{matrix}
  \mbox{distinguishing the classic Lusternik--Schnirel'man
countable sequence of$\quad\,\,\,\,\,\,$}\\
 \mbox{critical values and points is not possible without refined
numerical methods},
 \end{matrix}
 $}
  $$
  in view of  huge complicated
multiplicity of other admitted
 solutions.
 It is essential, that the arising problems
 do not admit, as customary for other classes of elliptic
 equations, any homotopy classification of solutions (say, on the
 hodograph plane), since all the compactly supported solutions are  infinitely
 oscillatory  that makes the homotopy rotational
 parameter infinite and hence the method non-applicable.

\smallskip

Let us now introduce these three classes, {\bf (I)--(III)}, of
PDEs and corresponding nonlinear phenomena to be studied by some
unified approaches.

\subsection{(I) Combustion models: regional blow-up, global stability, main goals,
and first discussion}

We begin with the following quasilinear degenerate $2m$th-order
parabolic equation of reaction-diffusion (combustion) type:
 \beq
    \label{S1}
 u_t = (-1)^{m+1}\D^{m}(|u|^n u)+|u|^n u \inA,
 \eeq
 where $n>0$ is a fixed exponent, $m \ge 2$ is integer, and $\D$ denoted the Laplace operator in
 $\re$.

\smallskip

\noi\underline{\em Globally asymptotically stable exact blow-up
solutions of S-regime}.
 In the simplest case $m=1$ and $N=1$, (\ref{S1}) is nowadays the canonical {\em quasilinear heat equation}
    \beq
    \label{RD.2}
    \mbox{$
    u_t=(u^{n+1})_{xx} + u^{n+1} \inC\quad (u \ge 0),
  $}
   \eeq
   which occurs in combustion theory.
The reaction-diffusion equation (\ref{RD.2}), playing a key role
in blow-up PDE theory, was under scrutiny since the middle of the
1970s. In 1976, Kurdyumov, with his former PhD students, Mikhailov
and Zmitrenko ({\em q.v.}  \cite{SZKM2} and \cite[Ch.~4]{SGKM} for
history) discovered the phenomenon of {\em heat} and {\em
combustion localization} by studying  the blow-up separate
variables {\em Zmitrenko--Kurdyumov solution} of
 (\ref{RD.2}):
 \beq
 \label{RD.31}
 u_{\rm S}(x,t)=(T-t)^{-\frac 1n} f(x) \quad \mbox{in}
 \quad \re \times (0,T),
  \eeq
  where $T>0$ is the blow-up time, and  
  $f$ satisfies the ODE
   \beq
   \label{RD.3}
   \mbox{$
 \frac 1n \, f= (f^{n+1})'' +  f^{n+1} \quad \mbox{for} \,\,\,\, x
 \in \re.
  $}
  \eeq
It turned out that (\ref{RD.3}) possesses the explicit compactly
supported solution
 \beq
 \label{RD.4}
 \mbox{$
   f(x) = \left\{ \begin{matrix}  \bigl[\frac {2(n+1)}{n(n+2)} \,  \cos^2
   \bigl(\frac  {n x}{2(n+1)}\bigr) \bigr]^{\frac 1n},
   \quad \,\mbox{if} \,\,\,\,\,\,\, |x| \le  \frac {n+1}n \, \pi, 
   \cr
   \qquad \, \qquad 0, \qquad \qquad \,\,\,\,\,\,\,
\quad \,
 \mbox{if} \,\,\,\,\,\,\,
    |x| >
   \frac {n+1}n \, \pi.
   \end{matrix} \right. $}
   \eeq
    This explicit integration of the ODE (\ref{RD.3}) was amazing and
rather surprising in the middle of the 1970s and led then to the
foundation of blow-up and heat localization theory. In dimension
$N>1$, the blow-up solution (\ref{RD.31}) does indeed exist
\cite[p.~183]{SGKM} but not in an explicit form (so that, it
seems, \ref{RD.4}) is the only available elegant form).

\smallskip

\noi\underline{\em Blow-up S-regime for higher-order parabolic
PDEs}. Evidently, the $2m$th-order counterpart (\ref{S1}) admits
the regional blow-up solution of the same form (\ref{RD.31}), but
the profile $f=f(y)$ then solves a more complicated ODE
 \beq
 \label{mm.561}
  \mbox{$
 (-1)^{m+1} \D^m (|f|^n f) + |f|^n f= \frac 1n \, f\quad \mbox{in} \,\,\, \ren.
 $}
  \eeq
After natural change, this gives the following equation with a
non-Lipschitz nonlinearity:
 $$ 
 \mbox{$
 F=|f|^n f \quad \Longrightarrow \quad
 \mbox{$(-1)^{m+1} \D^m F+F- \frac 1n \,
 \big| F \big|^{-\frac n{n+1}} F=0 \quad \mbox{in} \,\,\, \ren.
  $}
 $}
  $$ 
  Finally, we scale out the multiplier $\frac 1n$ in the nonlinear term,
\beq
 \label{S2NN}
 \mbox{$
 F \mapsto n^{-\frac{n+1}n} F \quad \Longrightarrow \quad
 \fbox{$(-1)^{m+1} \D^m F+F-
 \big| F \big|^{-\frac n{n+1}} F=0 \quad \mbox{in} \,\,\, \ren.
  $}
 $}
  \eeq
 In the one-dimensional case $N=1$, we obtain a simpler ODE,
\beq
 \label{S2}
 \mbox{$
 F \mapsto n^{-\frac{n+1}n} F \quad \Longrightarrow \quad
 \fbox{$(-1)^{m+1} F^{(2m)}+F-
 \big| F \big|^{-\frac n{n+1}} F=0 \quad \mbox{in} \,\,\, \re.
  $}
 $}
  \eeq
Thus, according to (\ref{RD.31}), the elliptic problems
(\ref{S2NN}) and the ODE (\ref{S2}) for $N=1$ are responsible for
the possible ``geometrical shapes" of regional blow-up described
by the higher-order combustion model (\ref{S1}).

\smallskip

\noi\underline{\em Plan and main goals of the paper related to
parabolic PDEs}.
Unlike the second-order case (\ref{RD.4}),
{\em explicit} compactly supported solutions $F(x)$ of (\ref{S2})
for any $m \ge 2$ are not available. Moreover, it turns out that
such profiles $F(x)$ have rather complicated local and global
structure. We are not aware of any rigorous or even formal
qualitative results concerning existence, multiplicity, or global
structure of solutions of ODEs such as (\ref{S2}).
 Our main goals are four-fold:

 \ssk

\noi{\bf (ii)} \underline{\sc Problem ``Blow-up"}: proving
finite-time blow-up in the parabolic (and hyperbolic) PDEs under
consideration (Section \ref{SBlow});

 \noi{\bf (ii)} \underline{\sc Problem ``Multiplicity"}:
 existence and
  multiplicity
  for elliptic PDEs (\ref{S2NN}) and the
ODEs (\ref{S2})
(Section \ref{SectVar});

\noi{\bf (iii)} \underline{\sc Problem ``Oscillations"}: the
generic structure of  oscillatory
 solutions of (\ref{S2}) near interfaces (Section \ref{Sect2});
 and


\noi{\bf (iv)} \underline{\sc Problem ``Numerics"}: numerical
study of various families of $F(x)$ (Sections \ref{Sect4},
\ref{Sectm4}).

\ssk

The research will be continued in the second half of this paper
\cite{GMPSobII}, where we intend to refine our results on the
multiplicity of solutions (especially, for $m>2$), pose the
problem on a ``Sturm index" of solutions (a homotopy
classification of some sub-families of solutions), and introduce
and study related analytic models with similar families of
solutions.


\ssk

Thus, in particular, we show that ODEs (\ref{S2}), as well as the
PDE (\ref{S2NN}), for any $m \ge 2$ admit infinitely many
countable families of compactly supported solutions in $\re$, and
the whole solution set exhibits certain {\em chaotic} properties.
Our analysis will be based on a combination of analytic
(variational and others), numerical, and various formal
techniques. Explaining existence,
 multiplicity, and asymptotics for the nonlinear problems
 involved, we leave several open mathematical problems.
 Some of these for higher-order equations are  extremely
 difficult.


 \subsection{(II) Regional blow-up in quasilinear hyperbolic equations}

  Secondly, consider the  $2m$th-order hyperbolic  counterpart of
  (\ref{S1}),
 \beq
 \label{S3}
 u_{tt} =(-1)^{m+1} \D^{m}(|u|^n u)+|u|^n u \inA.
  \eeq
We begin the discussion of its blow-up solutions in 1D, i.e., for
  \beq
    \label{RD.2H}
    u_{tt}=(u^{n+1})_{xx} +  u^{n+1}\quad \mbox{in} \quad \re \times \re_+ \quad (u \ge 0).
     \eeq
Here the blow-up solutions and the ODE take the form
 \beq
 \label{RD.31H}
 \mbox{$
 u_{\rm S}(x,t)=(T-t)^{-\frac 2n} \tilde f(x) \quad \Longrightarrow \quad
 \frac 2n \bigl(\frac 2n+1 \bigr) \, \tilde f= (\tilde f^{n+1})'' +  \tilde f^{n+1}.
 $}
  \eeq
Using extra scaling,
  \beq
 \label{S5}
 \mbox{$
 \tilde f(x)= \bigl[\frac {2(n+2)}n\bigr]^{\frac 1n} f(x)
   $}
  \eeq
yields the same ODE (\ref{RD.3}) and hence the exact localized
solution
 (\ref{RD.4}).

For the $N$-dimensional PDE (\ref{S3}),
 looking for  the same solution (\ref{RD.31H}),
 after scaling, leads to
  the elliptic equation (\ref{S2NN}).

 \subsection{(III) Nonlinear dispersion equations and compactons}

In a general setting, these rather unusual PDEs take the form
 \beq
 \label{NDEN}
 u_t= (-1)^{m+1}[\D^m(|u|^n u)]_{x_1} + (|u|^n u)_{x_1}
 \inA,
  \eeq
  where the right-hand side is just the derivation $D_{x_1}$ of
  that in the parabolic counterpart (\ref{S1}).
Then the elliptic problem (\ref{S2NN}) occurs when studying {\em
travelling wave} (TW) solutions of (\ref{NDEN}). As usual, we
explain this first in a simpler 1D case.

Let $N=1$ and $m=1$ in (\ref{NDEN}) that yields  the third-order 
   {\em Rosenau--Hyman} (RH)
{\em equation}
  \beq
  \label{Comp.4}
  \mbox{$
  u_t =  (u^2)_{xxx} + (u^2)_x,
  $}
  \eeq
  which  is known to model the effect of {\em nonlinear dispersion} in the pattern
  formation in liquid drops \cite{RosH93}. It is
   the $K(2,2)$ equation from the general $K(m,n)$ family of
   {\em nonlinear dispersion equations} (NDEs)
  \index{equation!$K(m,n)$}
  \beq
  \label{Comp.5}
   u_t =  (u^n)_{xxx} +  (u^m)_x \quad (u \ge 0),
   \eeq
   that also
 describe   phenomena of compact pattern
 formation, \cite{RosCom94, Ros96}. In addition, such PDEs  appear in curve motion and shortening
flows \cite{Ros00}.
 Similar to the previous models, the $K(m,n)$ equation (\ref{Comp.5}) with $n>1$
  is degenerated at $u=0$, and
 therefore may exhibit finite speed of propagation and admit solutions with finite
 interfaces. A permanent source of NDEs is integrable equation
 theory, e.g. look at the integrable
  fifth-order {\em
Kawamoto equation} \cite{Kaw85} (see \cite[Ch.~4]{GSVR} for other
models), which is of  NDE's type:
  \beq
  \label{Kaw1}
 u_t = u^5 u_{xxxxx}+ 5 \,u^4 u_x u_{xxxx}+ 10 \,u^5 u_{xx}
 u_{xxx}.
 \eeq
 Questions of local existence, uniqueness,
  regularity, shock and rarefaction wave formation, finite propagation and interfaces, including
  degenerate
 higher-order models to be studied are treated in  \cite{3NDEI, 3NDEII}; see also comments in
 \cite[Ch.~4.2]{GSVR}.

  We will study some  particular continuous
 solutions of the NDEs that give insight on several generic properties
 of such nonlinear PDEs.  
The crucial advantage of the RH equation
 (\ref{Comp.4}) is that it possesses
  {\em explicit}  moving compactly supported
 soliton-type solutions, called {\em compactons}
 \cite{RosH93}, which are {\em travelling wave} (TW) solutions:


\smallskip

\noi\underline{\em Compactons: manifolds of TWs and  blow-up
S-regime solutions coincide}.
 Let us show that such compactons are directly related to the
  blow-up patterns presented above.
 Actually, explicit TW compactons  exist  for the nonlinear dispersion KdV-type equations with
  arbitrary power nonlinearities (formulae will be given shortly)
   \beq
   \label{RD.1}
\mbox{$
   u_t= (u^{n+1})_{xxx} +  (u^{n+1})_x
 \inC.
  $}
  \eeq
   This is the $K(1+n,1+n)$ model, \cite{RosH93}.

   Thus, compactons as solutions of the equation (\ref{RD.1}) have the usual TW structure
 \beq
 \label{RD.5}
 u_{\rm c}(x,t)= f(y), \quad y= x- \l t,
  \eeq
so that, on substitution, $f$ satisfies the ODE
  \beq
 \label{RD.6SS}
 -\l f'=  (f^{n+1})''' +  (f^{n+1})',
  \eeq
 and, on integration once,  
 \beq
 \label{RD.6}
  -\l f=  (f^{n+1})'' +  f^{n+1}
 + D,
  \eeq
  where $D \in \re$ is a constant of integration. Setting
  $D=0$, which means the physical condition of
   zero flux at the interfaces,
   leads to  the blow-up ODE (\ref{RD.3}), so that
  the
  compacton equation (\ref{RD.6}) {\em coincide} with the blow-up one (\ref{RD.3})
  if 
  \beq
   \label{RD.7}
   \mbox{$
  \mbox{$ -\l =  \frac 1n $} \quad \bigl(\mbox{or \,$\mbox{$-\l= \frac 2n\bigl(\frac 2n+1\bigr)$}$\,  to
   match (\ref{RD.31H})}\bigr).
    $}
    \eeq
This yields  the compacton solution (\ref{RD.5}) with the same
compactly supported profile (\ref{RD.4}) with the translation $x
\mapsto y=x-\l t$.

Therefore,  in 1D, the 
blow-up solutions (\ref{RD.31}), (\ref{RD.31H}) of the parabolic
and hyperbolic PDEs and the compacton solution (\ref{RD.5})
 of the nonlinear dispersion equation (\ref{RD.1})
 are
essentially of a similar {\em mathematical} $($both the ODE and
PDE$)$ {\em nature}, and, possibly, more than that. This reflects
a certain {\em universality principle} of compact structure
formation in nonlinear evolution PDEs. Stability features of the
TW compacton (\ref{RD.5}) in the PDE setting (\ref{RD.1}) are
unknown, as well as for the higher-order counterparts to be posed
next.

In the $N$-dimensional geometry, i.e., for the PDE (\ref{NDEN}),
looking for a TW moving in the $x_1$-direction only,
 \beq
 \label{nn1}
  \mbox{$
 u_{\rm c}(x,t)=f(x_1-\l t,x_2,...,x_N) \quad \big(\l= - \frac
 1n\big)
  $}
 \eeq
 we obtain on integration in $x_1$ the elliptic problem (\ref{S2NN}).

 \ssk

Thus, we have introduced the necessary three classes, {\bf (I),
(II), (III)}, of nonlinear higher-order PDEs in $\ren \times
\re_+$, which, being representatives of very different  three
equations types, nevertheless will be shown to exhibit quite
similar evolution features (if necessary, up to replacing blow-up
by travelling wave moving), and the coinciding complicated
countable sets of evolution patterns. These common features reveal
an exiting concept of a certain unified principle of singularity
formation phenomena in general nonlinear PDE theory, that we seem
just begin to touch and study in the twenty first century. Several
classic mathematical concepts and techniques successfully
developed in the twentieth century including, of course,  Sobolev
 legacy continue to be  key, but also new ideas from different
 ranges of various rigorous and qualitative natures are
 desperately needed for tackling such fundamental difficulties and
 open problem arising.

\section{{\bf Problem ``Blow-up":} general blow-up analysis of parabolic and hyperbolic
PDEs}
 \label{SBlow}

\subsection{On global existence and blow-up in higher-order parabolic equations}

We begin with the parabolic model (\ref{S1}). Bearing in mind a
compactly supported nature of the solutions under consideration,
we consider  (\ref{S1}) in a bounded domain $\O \subset \ren$ with
a smooth boundary $\partial \O$, with  Dirichlet boundary
conditions
 \beq
 \label{1.2}
 u=Du=...=D^{m-1}u=0 \quad \mbox{on}
 \quad \partial \O \times \re_+,
  \eeq
 and a given sufficiently smooth and bounded initial data
 \beq
 \label{1.3}
 u(x,0)=u_0(x) \quad \mbox{in} \quad \O.
 \eeq
We will show that the phenomenon of blow-up depends essentially on
the size of the domain. Beforehand, let us observe that the
diffusion operator on the right-hand side in (\ref{S1}) is a
monotone operator in $H^{-m}(\O)$, so that the unique local
solvability of the problem in suitable Sobolev spaces is covered
by classic theory of monotone operators; see Lions' book
\cite[Ch.~2]{JLi}. We next show that, under certain conditions,
some of these solutions are global in time but some ones cannot be
globally extended and blow-up in finite time.

 For convenience, we
use the natural substitution
 \beq
 \label{1.0}
 v=|u|^n u \LongA v_0(x)=|u_0(x)|^n u_0(x),
  \eeq
that leads to the following parabolic equation with a standard
linear operator on the right-hand side:
 \beq
 \label{1.1}
  \mbox{$
   (\psi(v))_t= (-1)^{m+1} \D^m v + v \inA, \quad \mbox{with}
   \quad \psi(v)=|v|^{-\frac n{n+1}} v,
    $}
    \eeq
where $v$ satisfies the same Dirichlet boundary condition
(\ref{1.2}).

Multiplying (\ref{1.1}) by $v$ in $L^2(\O)$ and integrating by
parts by using (\ref{1.2}) yields
 \beq
 \label{1.4}
  \mbox{$
  \frac {n+1}
{n+2}\, \frac{\mathrm d}{{\mathrm d}t} \,
  \int\limits_{\O} |v|^{\frac{n+2}{n+1}} {\mathrm d}x=
-\int\limits_{\O} |\tilde D^m v|^2 {\mathrm d}x + \int\limits_{\O}
v^{2} {\mathrm d}x \equiv E(v),
 $}
 \eeq
where we use the notation $\tilde D^m= \D^{\frac m2}$ for even
  $m$ and $\tilde D^m= \n \D^{\frac {m-1}2}$ for odd $m$.
By Sobolev's embedding theorem, $H^m(\O) \subset L^2(\O)$
compactly, and moreover, the following sharp estimates holds:
 \beq
 \label{1.5}
  \mbox{$
\int\limits_{\O} v^{2} {\mathrm d}x \le \frac{1}{\l_1}\,
\int\limits_{\O} |\tilde D^m v|^2 {\mathrm d}x \quad \mbox{in}
\quad H^m_0(\O),
 $}
  \eeq
  where $\l_1=\l_1(\O)>0$ is the first simple eigenvalue of the
  poly-harmonic operator $(-\D)^m$ with the Dirichlet boundary
  conditions (\ref{1.2}):
   \beq
   \label{DD1}
   (-\D)^m e_1 = \l_1 e_1 \quad \mbox{in}
   \quad \O, \quad e_1 \in H^{2m}_0(\O).
    \eeq
    For $m=1$,
 since $(-\D) >0$ is strictly negative in the metric of $L^2(\O)$,
 we have, by   Jentzsch's classic theorem (1912) on the positivity
 of the first eigenfunction for linear integral operators with positive
 kernels,
 that
  \beq
  \label{eeee1}
   e_1(x) > 0 \quad \mbox{in} \quad \O.
    \eeq
For $m \ge 2$, (\ref{eeee1}) remains valid, e.g., for the unit
ball $\O=B_1$. Indeed, in the case of $\O=B_1$, the Green function
of the  poly-harmonic operator
 $(-\D)^m$
with Dirichlet boundary conditions is positive; see first results
by Boggio (1901-05) \cite{Boggio1, Boggio2} (see also Elias
\cite{Elias} for later related general results). Again,  by
Jentzsch's theorem,
 (\ref{eeee1}) holds.


  It follows from (\ref{1.4}) and (\ref{1.5}) that
 $$
  \mbox{$
  \frac {n+1}{n+2}\, \frac{\mathrm d}{{\mathrm d}t} \,
  \int\limits_{\O} |v|^{\frac{n+2}{n+1}} {\mathrm d}x \le
  \big(\frac 1{\l_1}-1\big)
\int\limits_{\O} |\tilde D^m v|^2 {\mathrm d}x.
 $}
 $$


 \noi\underline{\em Global existence for $\l_1>1$}.
Thus, we obtain the following inequality:
 \beq
 \label{1.6}
  \mbox{$
  \frac {n+1}{n+2}\, \frac{\mathrm d}{{\mathrm d}t} \,
  \int\limits_{\O} |v|^{\frac{n+2}{n+1}} {\mathrm d}x +
  \big(1-\frac 1{\l_1}\big)
\int\limits_{\O} |\tilde D^m v|^2 {\mathrm d}x \le 0.
 $}
  \eeq
 Consequently, for
  \beq
  \label{1.7}
   \l_1(\O)>1,
    \eeq
 (\ref{1.5}) yields good {\em a priori} estimates of solutions in
 $\O \times (0,T)$ for arbitrarily large $T>0$. Then, by the
 standard Galerkin method \cite[Ch.1~]{JLi}, we get global
 existence of solutions of the initial-boundary value problem
 (IBVP)
 (\ref{1.1}), (\ref{1.2}), (\ref{1.3}). This means no finite-time blow-up for
 the IBVP provided (\ref{1.7}) holds, meaning that the size of
 domain being sufficiently small.

\ssk

 \noi\underline{\em Global existence for $\l_1=1$}.
  Note that for $\l_1=1$, (\ref{1.6}) also yields an {\em a priori}  uniform
  bound, which is weaker, so the proof of global existence becomes
  more tricky and requires extra scaling to complete (this is not directly related
  to the present discussion, so we omit details). In this case, we have the conservation
 law
  \beq
  \label{cc1}
   \mbox{$
    \int\limits_\O \psi(v(t)) e_1\, {\mathrm d}x=c_0= \int\limits_\O \psi(v_0) e_1\, {\mathrm d}x
 \quad \mbox{for all} \quad t>0,
  $}
  \eeq
so that by the gradient system property (see below), the global
bounded orbit must stabilize to a unique stationary solution, which is characterized as
follows (recall that $\l_1$ is always a simple eigenvalue, so the eigenspace is 1D):
 \beq
 \label{cc2}
  \mbox{$
  v(x,t) \to C_0 e_1(x) \asA t \to + \iy \whereA \int\limits_\O \psi(C_0 e_1) e_1\, {\mathrm
  d}x=c_0.
   $}
   \eeq

\ssk

 \noi\underline{\em Blow-up for $\l_1 < 1$}.
 Let us now show that for the opposite inequality
 \beq
  \label{1.9N}
   \l_1(\O)<1,
    \eeq
the solutions blow-up in finite time.

\ssk

\noi\underline{\em Blow-up of nonnegative solutions for $m=1$}. We
begin with the simpler case $m=1$, where, by the Maximum
Principle, we can restrict to the class of nonnegative solutions
 \beq
 \label{v11}
 v=v(x,t) \ge 0, \quad \mbox{i.e., assuming that $u_0(x) \ge 0$}.
 \eeq
 In this case,
  we can directly study  the
evolution of the first Fourier coefficient of the function
$\psi(v(\cdot,t))$. To this end, we multiply (\ref{1.1}) by the
positive eigenfunction $e_1$ in $L^2(\O)$ to obtain that
 \beq
 \label{1.9}
  \mbox{$
\frac{\mathrm d}{{\mathrm d}t} \,
  \int\limits_{\O} \psi(v) e_1 {\mathrm d}x =
  (1-\l_1)
\int\limits_{\O} v e_1 {\mathrm d}x.
 $}
  \eeq
In view of (\ref{v11}), we apply H\"older's inequality in the
right-hand side of (\ref{1.9}) to derive the following ordinary
differential inequality for the Fourier coefficient:
 \beq
 \label{2.0}
  \mbox{$
\frac{{\mathrm  d}J}{{\mathrm d}t}  \ge
  (1-\l_1)c_2 J^{n+1} \whereA
 J(t)= \int\limits_\O v^{\frac 1{n+1}}(x,t) e_1(x)\, {\mathrm d}x,
  \,\,\, c_2=\big( \int\limits_\O e_1\,{\mathrm d}x\big)^{-n}.
 $}
  \eeq
Hence, for any nontrivial nonnegative initial data
 $$
  \mbox{$
 u_0(x) \not \equiv 0 \LongA J_0= \int_\O v_0 e_1\, {\mathrm d}x
 >0,
 $}
  $$
we have finite-time blow-up of the solution with  the following
lower estimate on the Fourier coefficient:
 \beq
 \label{JJ1}
  \mbox{$
 J(t) \ge A (T-t)^{-\frac 1n} \whereA A=\big( \frac
 1{nc_2(1-\l_1)}\big)^{\frac 1n}, \quad T= \frac {J_0^{-n}}{n
 c_2(1-\l_1)}.
  $}
  \eeq

\ssk

\noi\underline{\em On unbounded orbits and blow-up for $m \ge 2$}.
 It is curious that we do not know a similar simple proof of blow-up for
 the higher-order equations with $m \ge 2$.  The main technical
 difficulty is that the set of nonnegative solutions (\ref{v11})
 is not invariant of the parabolic flow, so we have to deal with
 solutions $v(x,t)$ of changing sign. Then,  (\ref{2.0}) cannot be
 derived from (\ref{1.9}) by the H\"older inequality.
 Nevertheless, we easily obtain the following result as a first
 step to blow-up of the orbits:

  \begin{proposition}
   \label{Pr.Un}
   Let $m \ge 2$, $(\ref{1.9N})$ hold, and
 \beq
 \label{EE1}
 E(v_0)>0.
  \eeq
Then the solution of the
   IBVP  $(\ref{1.1})$, $(\ref{1.2})$, $(\ref{1.3})$ is not
   uniformly bounded for $t>0$.
    \end{proposition}

     \noi {\em Proof.} We use the obvious fact that (\ref{1.1})
 is a gradient system in $H^m_0(\O)$. Indeed, multiplying
 (\ref{1.1}) by $v_t$ yields, on sufficiently smooth local solutions,
 \beq
 \label{SS1}
  \mbox{$
   \frac 12 \frac{\mathrm d}{{\mathrm d}t} E(v(t)) = \frac 1{n+1}
   \int\limits_\O |v|^{-\frac n{n+1}} (v_t)^2 \, {\mathrm d}x \ge
   0.
   $}
    \eeq
Therefore, under the hypothesis (\ref{EE1}) we have from
(\ref{1.4}) that
 \beq
 \label{SS2}
  \mbox{$
  E(v(t)) \ge E(v_0)>0 \LongA
  \frac {n+1}
{n+2}\, \frac{\mathrm d}{{\mathrm d}t} \,
  \int\limits_{\O} |v|^{\frac{n+2}{n+1}} {\mathrm d}x=
 E(v) \ge E(v_0)>0, \quad \mbox{i.e.,}
 $}
 \eeq
  \beq
  \label{SS3}
   \mbox{$
 \int\limits_\O |v(t)|^{\frac{n+2}{n+1}} {\mathrm d}x \ge
  \frac{n+2}{n+1}\,  E(v_0) \, t \to + \iy \asA t \to + \iy. \qed
   $}
 \eeq

  \ssk

Concerning  the  hypothesis (\ref{EE1}),
 recall that by classic dynamical system theory \cite{Hale}, the
 $\o$-limit set of bounded orbits of gradient systems
 consists of equilibria only, i.e.,
  \beq
   \label{nnn1}
    \o(v_0) \subseteq  {\mathcal S}=\big\{V \in H^m_0(\O): \,\,\, -(-\D)^m V +
    V=0\big\}.
     \eeq
   Therefore, stabilization to a nontrivial equilibrium is possible if
    $$
    \l_l=1 \quad \mbox{for some} \,\,\, l \ge 2.
    $$
    Otherwise, we have that
   \beq
   \label{RR10}
  {\mathcal S}=\{0\} \quad (\l_l \not =1 \,\,\, \mbox{for any}\,\,\, l \ge 1).
  \eeq
   Then,
   formally, by
 the gradient structure of (\ref{1.1}),
  one should take into account solutions that
  decay to 0
  as $t \to +\iy$. One can check that (at least formally, a
  necessary functional framework could take some time), the trivial solution 0
 has the empty stable manifold, so that, under the assumption (\ref{RR10}), the result in Proposition
 \ref{Pr.Un} is naturally expected to be true for any nontrivial
 solution.






Thus, we have that in the case (\ref{EE1}), i.e., for sufficiently
large domain $\O$, solutions become arbitrarily large in any
suitable metric, including $H^m_0(\O)$ or the uniform one
$C_0(\O)$. Then, it is a technical matter to show that
 such large solutions  must next blow-up in finite time. In fact, often, this is not
 that straightforward, and omitting this blow-up analysis, we would like
 to attract the attention of the interested Reader.

\ssk

\noi\underline{\em Blow-up for $m \ge 2$ in
 a similar modified model}. On
the other hand, the previous proof of blow-up is easily adapted
for the following slightly modified equation (\ref{1.1}):
 \beq
 \label{1.1M}
 (\psi(v))_t=(-1)^{m+1} \D^m v+ |v|,
  \eeq
  where the source term is replaced by $|v|$. Actually, for
  ``positively dominant" solutions (i.e., for those of a non-zero integral
   $\int u(x,t)\, {\mathrm d}x$), this is not a big change, and
  most of our self-similar patterns perfectly exist for
  (\ref{1.1M}) and the oscillatory properties of solutions near
  interfaces remain practically untouched (since the source term
  plays no role there).

  Take $\O=B_1$, so that (\ref{eeee1}) holds. Then, instead of
  (\ref{1.9}), we will get a similar inequality,
 \beq
 \label{1.9M}
  \mbox{$
\frac{\mathrm d}{{\mathrm d}t} \,
  \int\limits_{\O} \psi(v) e_1 {\mathrm d}x =
\int\limits_{\O} |v| e_1 {\mathrm d}x - \l_1 \int\limits_{\O} v
e_1 {\mathrm d}x \ge (1-\l_1)\int\limits_{\O} |v| e_1 {\mathrm d}x
>0,
 $}
  \eeq
 where $J(t)$ is defined without the positivity sign restriction,
  \beq
  \label{PP1}
   \mbox{$
J(t)= \int_\O (|v|^{-\frac n{n+1}}v)(x,t) e_1(x)\, {\mathrm d}x.
 $}
  \eeq
It follows from (\ref{1.9M}) that, for $\l_1<1$,
 \beq
 \label{PP2}
 J(0)>0 \quad \Longrightarrow \quad J(t)>0 \,\,\,\mbox{for}
 \,\,\,t>0.
  \eeq
Therefore, by the H\"older inequality,
  \beq
  \label{MMM1}
  \mbox{$
   \int |v|e_1 \, {\mathrm d}x \ge c_2 \big( \int |v|^{\frac 1{n+1}}
   e_1 \, {\mathrm d}x\big)^{n+1} \ge c_2 \big( \int |v|^{-\frac n{n+1}}v
   e_1 \, {\mathrm d}x\big)^{n+1} \equiv c_2 J^{n+1}.
   $}
   \eeq
This allows us to
get the inequality (\ref{2.0}) for the function
(\ref{PP1}). Hence, the blow-up estimate (\ref{JJ1}) holds.

\subsection{Blow-up data for higher-order parabolic and hyperbolic PDEs}

We have seen above that, in general, blow-up occurs for some
initial data, since in many cases, small data can lead to globally
existing sufficiently small solutions (of course, if 0 has a
nontrivial stable manifold).

Below, we introduce classes of such ``blow-up data", i.e., initial
functions generating finite-time blow-up of solutions. Actually,
studying such crucial data will eventually require to perform a
detailed study of the corresponding elliptic systems with
non-Lipschitz nonlinearities.

\ssk

\noi\underline{\em Parabolic equations}. To this end, again
beginning with the transformed parabolic equation (\ref{1.1}), we
consider the separate variable solutions
 \beq
 \label{2.1PP}
  \mbox{$
  v(x,t)
= (T-t)^{ -\frac{n+1}n} F(x).
   $}
    \eeq
 Then $F(x)$ solves the elliptic equation (\ref{S2NN}) in $\O$,
 i.e.,
  \beq
  \label{2.2PP}
   \left\{
    \begin{matrix}
   (-1)^{m+1} \D^m F + F- \frac 1n \,|F|^{-\frac n{n+1}} F=0 \inB \O, \ssk\ssk\\
    F=DF=...=D^{m-1}F=0 \quad \mbox{on} \quad \partial \O.
    \qquad\quad
 \end{matrix}
     \right.
 \eeq
 Let $F(x) \not \equiv 0$ be a solution of the problem
 (\ref{2.2PP}), which is a key object in the present paper.  Hence,
 it follows from (\ref{2.1PP}) that initial data
  \beq
  \label{2.4PP}
   v_0(x) = C F(x),
    \eeq
    where $C \not =0$ is an arbitrary constant to be scaled out,
    generate blow-up of the solution of (\ref{1.1}) according to
    (\ref{2.1PP}).

\ssk

\noi\underline{\em Hyperbolic equations}. Similarly, for the
hyperbolic counterpart of (\ref{1.1}),
 \beq
 \label{HH1}
  (\psi(v))_{tt}= (-1)^{m+1} \D^m v+v,
   \eeq
 we take initial data in the form
  \beq
  \label{HH2}
   v(x,0)= c F(x) \andA v_t(x,0)= c_1 F(x),
    \eeq
  with some constants $c$ and $c_1$ such that $cc_1>0$.
  Then the solution blows up in finite time.
  In particular, choosing
   $$
    \mbox{$
   c>0 \andA c_1 = \frac {2(n+1)}n \, B^\frac 1 \b c^{1-\frac 1
   \b},
    $}
    $$
    with $\b=- \frac{2(n+1)}n$ and $B= \big[
    \frac{2(n+2)}{n^2}\big]^{\frac {n+1}n}$, we have the blow-up
    solution of (\ref{HH1}) in the separable form
     $$
      \mbox{$
     v(x,t)=(T-t)^\b B F(x) \whereA T=\big( \frac c B\big)^ {\frac
     1 \b}.
      $}
      $$


\subsection{Blow-up rescaled equation as a gradient system:
towards the generic blow-up behaviour for parabolic PDEs}

Let us briefly discuss another important issue associated with the
scaling (\ref{2.1PP}). Consider a general solution $v(x,t)$ of the
IBVP for (\ref{1.1}), which blows up first time at $t=T$.
Introducing the rescaled variables
 \beq
 \label{in1}
  v(x,t)=(T-t)^{-\frac {n+1}{n}} w(x,\t), \quad \t= - \ln(T-t) \to +\iy \asA t \to T^-,
   \eeq
one can see that $w(x,\t)$ then solves the following rescaled
equation:
 \beq
 \label{in2}
  \mbox{$
  (\psi(w))_\t =  (-1)^{m+1} \D^m w + w- \frac 1n \,|w|^{-\frac n{n+1}} w,
   $}
   \eeq
   where on the right-hand side we observe the same operator with
   a non-Lipschitz nonlinearity as in (\ref{S2NN}) or (\ref{2.2PP}).
By an analogous argument, (\ref{in2}) is a gradient system and
admits a Lyapunov functions that is strictly monotone on
non-equilibrium orbits:
 \beq
 \label{in3}
  \mbox{$
   \frac{\mathrm d}{{\mathrm d}\t} \big( -\frac 12 \int |\tilde
   D w|^2 + \frac 12 \, \int w^2 - \frac{n+1}{n(n+2)}\, \int
   |w|^{\frac{n+2}{n+1}} \big)= \frac 1{n+1} \, \int |w|^{- \frac
   n{n+1}} |w_t|^2 > 0.
    $}
     \eeq
Therefore, the corresponding to (\ref{nnn1}) conclusion holds,
i.e., all bounded orbits can approach stationary solutions only:
 \beq
 \label{in4}
  \mbox{$
  \o(w_0) \subseteq {\mathcal S}= \big\{ F \in H^m_0(\O): \,\,\, (-1)^{m+1} \D^m F + F- \frac 1n \,
  |F|^{-\frac n{n+1}}
  F=0\big\}.
   $}
   \eeq
Moreover, since under natural smoothness parabolic properties,
$\o(w_0)$ is connected and invariant \cite{Hale}, the omega-limit
set reduces to a single equilibrium provided that ${\mathcal S}$
is disjoint, i.e., consists of isolated points. Here, the
structure of the stationary rescaled set ${\mathcal S}$ becomes
key for understanding blow-up behaviour of general solutions of
the higher-order parabolic flow (\ref{S1}).

\ssk

Thus,  the above analysis shows again  that the ``stationary"
elliptic problems (\ref{S2NN}) and (\ref{2.2PP}) are crucial for
revealing various local and global evolution properties of all
three classes of PDEs involved. We begin this study with an
application of classic variational techniques.

\section{\underline{\bf Problem ``Existence"}: variational approach and
 countable families of solutions
by Lusternik--Schnirel'man  category and fibering theory}
 \label{SectVar}

\subsection{Variational setting and compactly supported solutions}

Thus, we need to study, in a general multi-dimensional geometry,
existence and   multiplicity of compactly supported solutions of
the elliptic problem in (\ref{S2NN}).


 Since all the operators in (\ref{S2NN}) are potential, the problem
admits the variational setting in $L^2$, so the solutions can be
obtained as critical points of the following $C^1$ functional:
 \beq
 \label{V1}
  \mbox{$
 {E}(F)= - \frac 12  \int |\tilde D^m F|^2 + \frac 12 \int
 F^2 -\frac{1}{\b} \, \int |F|^{\b}, \quad \mbox{with} \,\,\, \b=\frac {n+2}{n+1}\in (1,2),
  $}
  \eeq
  where, as above,  $\tilde D^m= \D^{\frac m2}$ for even
  $m$ and $\tilde D^m= \n \D^{\frac {m-1}2}$ for odd $m$.
  In general, we have to look for critical points in $W^2_m(\ren)
  \cap L^2(\ren) \cap L^{\b}(\ren)$.
  Bearing in mind 
 compactly supported solutions, we choose a sufficiently
  large radius $R>0$ of the ball $B_R$ and consider the variational problem for (\ref{V1})
  in $W_{m,0}^2(B_R)$, where we assume Dirichlet boundary
  conditions on $S_R= \partial B_R$. Then
   both spaces $L^2(B_R)$ and $L^{p+1}(B_R)$ are
  compactly embedded into $W_{m,0}^2(B_R)$ in the subcritical Sobolev
  range
   \beq
   \label{Sob1}
    \mbox{$
   1<p<p_S= \frac {N+2m}{N-2m} \quad (\b < p_S).
    $}
    \eeq
    In general,
we have to use the following preliminary result:



\begin{proposition}
 \label{Pr.CS}
  Let $F$ be a continuous weak solution of the
  equation $(\ref{S2NN})$ such that
   \beq
   \label{ffy1}
   F(y) \to 0 \asA |y| \to \infty.
    \eeq
Then $F$ is compactly supported in $\ren$.
 \end{proposition}

 Notice that continuity of $F$ is guaranteed for $N<2m$ directly
 by Sobolev embedding $H^m(\ren) \subset C(\ren)$, and, in the
 whole range (\ref{Sob1}), by local elliptic regularity theory; see
 necessary embeddings of functional spaces in Maz'ja \cite[Ch.~1]{Maz}.

\smallskip

\noi {\em Proof.} We consider the corresponding parabolic equation
with the same elliptic operator,
 \beq
 \label{Par1}
  \mbox{$
 w_t = (-1)^{m+1} \D ^m w +w   -
 \big| w \big|^{- \frac n{n+1}} w \quad \mbox{in} \,\,\, \ren
 \times\re_+,
  $}
  \eeq
  with initial data $F(y)$. Setting $w = {\mathrm e}^t \hat w$
yields the equation
 $$
  \mbox{$
 \hat w_t = (-1)^{m+1} \D ^m \hat w -{\mathrm e}^{-\frac n{n+1} \, t}\,
 | \hat w |^{p-1}\hat w,
 \quad \mbox{where} \quad p = \frac 1{n+1} \in (0,1),
  $}
 $$
where the operator is monotone in $L^2(\ren)$. Therefore, the
Cauchy problem (CP) with initial data $F$ has a unique weak
solution, \cite[Ch.~2]{JLi}.
 Thus, (\ref{Par1}) has the unique solution
 $w(y,t) \equiv F(y)$,
  which then
  must be  compactly supported for arbitrarily small $t>0$.
  Indeed, such nonstationary  instant compactification phenomena for quasilinear
  absorption-diffusion equations
 with singular absorption $-|u|^{p-1}u$, with $p < 1$,
   have been known
  since  the 1970s and are also called the {\em instantaneous
 shrinking} of the support
     of
 solutions. These phenomena
  have been  proved for
  quasilinear higher-order parabolic equations with non-Lipschitz
  absorption terms; see \cite{Shi2}.
   $\qed$

\smallskip

 Thus, to revealing  compactly supported
patterns $F(y)$, we can pose the problem in bounded balls that are
sufficiently large.
 Indeed, one can see from (\ref{V1}) that, in small domains, nontrivial solutions
are impossible.



 \subsection{L--S theory and direct
 application of fibering method}

 The functional (\ref{V1}) is $C^1$ and is
uniformly differentiable and weakly continuous, so we can apply
classic
 Lusternik--Schnirel'man (L--S) theory of calculus of variations
\cite[\S~57]{KrasZ} in the form of the fibering method \cite{Poh0,
PohFM},  as a convenient generalization of previous versions
\cite{Clark, Rabin} of variational approaches.




 Namely, following  L--S theory and the fibering
approach \cite{PohFM}, the number of critical points of the
functional (\ref{V1}) depends on the {\em category} (or {\em
genus}) of the functional subset on which fibering is taking
place. Critical points of ${E(F)}$ are obtained by
 {\em spherical
fibering}
 \beq
 \label{f1}
 F= r(v) v \quad (r \ge 0),
  \eeq
  where $r(v)$ is a scalar functional, and $v$ belongs to a subset
  in  $W_{m,0}^2(B_R)$ given as follows:
   \beq
   \label{f2}
    \mbox{$
    {\mathcal H}_0=\bigl\{v \in W_{m,0}^2(B_R): \,\,\,H_0(v)
     \equiv  -  \int |\tilde D^m v|^2 +  \int
 v^2 =1\bigr\}.
    $}
    \eeq
The new functional
 \beq
 \label{f3}
  \mbox{$
H(r,v)= \frac 12 \, r^2 - \frac 1{\b}\, r^{\b} \int |v|^{\b}
 $}
  \eeq
 has the absolute minimum point, where
 \beq
 \label{f31}
  \mbox{$
 H'_r \equiv r-  r^{\b-1} \int |v|^{\b} =0
  \,\,\Longrightarrow \,\,
   r_0(v)=\bigl(\int |v|^{\b}\bigr)^{\frac 1{2-\b}}.
   $}
   \eeq 
   We then obtain the following functional:
 \beq
 \label{f400}
 \mbox{$
  \tilde H(v)= H(r_0(v),v)=- \frac {2-\b}{2\b} \,  r_0^2(v)
  \equiv - \frac {2-\b}{2\b} \, \big( \int |v|^\b\big)^{\frac
  2{2-\b}}.
   $}
    \eeq
 The critical points of the functional (\ref{f400}) on the
 set (\ref{f2}) coincide with those for
 \beq
  \label{f4}
   \mbox{$
    \tilde H(v)= \int |v|^\b,
    $}
    \eeq
 so
  we arrive at even, non-negative, convex,
 and uniformly differentiable functional, to which
L--S  theory applies, \cite[\S~57]{KrasZ}; see also
\cite[p.~353]{Deim}.
   Following \cite{PohFM}, searching
  for critical points of $\tilde H$ on the set ${\mathcal H}_0$
  one needs to estimate the category $\rho$
  of the set ${\mathcal H}_0$.
The details on this notation and basic results can be  found in
Berger \cite[p.~378]{Berger}. Notice that the Morse index $q$ of
the quadratic form $Q$ in Theorem 6.7.9 therein is precisely the
dimension of the space where the corresponding form is negatively
definite. This includes all the multiplicities of eigenfunctions
involved in the corresponding subspace.
 See also genus and cogenus definitions and applications to variational problems in
 \cite{BB81} and \cite{BL88}.

It follows that, by this variational construction, $F$ is an
eigenfunction satisfying
 $$
 \mbox{$
 (-1)^{m+1} \D^m F + F - \mu \big|F \big|^{-\frac n{n+1}}F=0,
 $}
  $$
where $\mu >0$ is Lagrange's multiplier. Then scaling $F \mapsto
\mu^{(n+1)/n} F$ yields the original equation in (\ref{S2NN}).


For further discussion of geometric shapes of patterns, it is
convenient to recall that utilizing Berger's version
\cite[p.~368]{Berger} of this minimax analysis of L--S category
theory \cite[p.~387]{KrasZ},  the critical values $\{c_k\}$ and
the corresponding critical points $\{v_k\}$ are given by
 \beq
 \label{ck1}
  \mbox{$
 c_k = \inf_{{\mathcal F} \in {\mathcal M}_k} \,\, \sup_{v \in {\mathcal
 F}} \,\, \tilde H(v),
  $}
  \eeq
where  ${\mathcal F} \subset {\mathcal H}_0$ are  closed sets,
 and
 ${\mathcal M}_k$ denotes the set of all subsets of the form
  $
  B S^{k-1}
\subset {\mathcal H}_0,
 $
 where $S^{k-1}$ is a suitable sufficiently
smooth $(k-1)$-dimensional manifold (say, sphere) in ${\mathcal
H}_0$ and $B$ is an odd continuous map.
 Then each member of ${\mathcal M}_k$ is of  genus at least $k$
 (available in ${\mathcal H}_0$).
   It is also important to remind that the
definition of genus \cite[p.~385]{KrasZ} assumes  that
$\rho({\mathcal F})=1$, if no {\em component} of ${\mathcal F}
\cup {\mathcal F}^*$, where
 $
 {\mathcal F}^*=\{v: \,\, v^*=-v \in {\mathcal F}\},
 $
 is the {\em reflection} of ${\mathcal F}$ relative to 0,
 contains a pair of antipodal points $v$ and $v^*=-v$.
 Furthermore, $\rho({\mathcal F})=n$ if each compact subset of
${\mathcal F}$ can be covered by, minimum, $n$ sets of genus one.
According to (\ref{ck1}),
 $$
 c_1 \le c_2 \le ... \le c_{l_0},
 $$
 where $l_0=l_0(R)$ is the category of ${\mathcal H}_0$ (see an estimate
 below) such that
  \beq
  \label{l01}
  l_0(R) \to + \infty \quad \mbox{as} \quad R \to \infty.
  \eeq
  Roughly speaking,
since the dimension of the sets ${\mathcal F}$ involved in the
construction of ${\mathcal M}_k$ increases with $k$, this
guarantees that the critical points delivering critical values
(\ref{ck1}) are all different.
  It follows from (\ref{f2}) that the category
$l_0=\rho({\mathcal H}_0)$  of the set ${\mathcal H}_0$ is equal
to the number (with multiplicities) of the eigenvalues $\l_k < 1$,
  \beq
  \label{pp1nn}
 \mbox{$
  l_0= \rho({\mathcal H}_0)= \sharp \{ \l_k < 1\}
 $}
  \eeq
 of the linear poly-harmonic operator $(-1)^{m} \D^m > 0$,
 \beq
 \label{f55}
 (-1)^{m} \D^m \psi_k= \l_k \psi_k, \quad \psi_k \in W^2_{m,0}(B_R);
  \eeq
  see \cite[p.~368]{Berger}.
  Since the dependence of the spectrum on $R$ is, obviously,
   \beq
   \label{f56}
   \l_k(R)= R^{-2m} \l_k(1), \quad k=0,1,2,... \, ,
    \eeq
we have that the category $\rho({\mathcal H}_0)$ can be
arbitrarily large for $R \gg 1$, and (\ref{l01}) holds. We fix
this in the following:

\begin{proposition}
 \label{Pr.MM}
The elliptic problem $(\ref{S2NN})$ has at least a countable set
of different solutions denoted by $\{F_l, \, l \ge 0\}$, each one
$F_l$ obtained as a critical point of the functional $(\ref{V1})$
in $W^2_{m,0}(B_R)$ with a sufficiently large $R=R(l)>0$.
 \end{proposition}

Indeed, in view of Proposition \ref{Pr.CS},  we  choose $R \gg 1$
such that ${\rm supp}\, F_l \subset B_R$.


\subsection{On a model supplying with explicit description of the
L--S sequence}

As we will see shortly, detecting the L--S sequence of critical
values for the original functional (\ref{V1}) is a hard problem,
where numerical estimates of the functional will be key.

However, there exist similar models, for which this can be done
much easier. We now perform a slight modification in (\ref{V1})
and consider the functional
 \beq
 \label{V1N}
 \mbox{$
 {E}_1(F)= - \frac 12  \int |\tilde D^m F|^2 + \frac 12 \int
 F^2 -\frac{1}{\b} \, \big( \int F^{2}\big)^{\frac \b 2} \quad \big(\b=\frac {n+2}{n+1}\in (1,2)
  \big).
  $}
  \eeq
  This corresponds to the following non-local elliptic problem:
   \beq
   \label{ee1}
   \mbox{$
   -(-\D)^m F + F - F \big( \int F^{2}\big)^{\frac \b 2-1}=0 \quad
   (\mbox{in} \,\,\, B_R, \,\, \mbox{etc.})
    $}
    \eeq
 Denoting by $\{\l_k\}$ the spectrum in (\ref{f55}) and by
 $\{\psi_k\}$ the corresponding eigenfunction set, we can solve
 the problem (\ref{ee1}) explicitly: looking for solutions
  \beq
  \label{ee2}
   \mbox{$
   F= \sum_{(k \ge 1)} c_k \psi_k \,\,
    \Longrightarrow \,\,c_k \big[ -\l_k+1 - \big( \sum_{(j \ge 1)}
    c_j^2 \big)^{\frac \b2-1} \big]=0, \,\, k=1,2,... \, .
     $}
     \eeq
 The algebraic system in (\ref{ee2}) is easy and yields precisely
the number (\ref{pp1nn})
  of various nontrivial basic solutions $F_l$ having the form
 \beq
 \label{pp2}
  F_l(y)= c_l \psi_l(y), \quad \mbox{where} \quad |c_l|^{\b-2}= -\l_l+1
  >0, \,\,\, l=1,2,...,\,l_0.
   \eeq

\subsection{Preliminary analysis of geometric shapes of patterns}

The forthcoming  discussions and conclusions should be understood
in conjunction with the results obtained in Section \ref{Sect4} by
numerical and other analytic and formal methods.
 In particular, we use here the concepts of the index and Sturm
 classification of various basic and other patterns.

Thus, we now  discuss   key questions of the spatial structure of
patterns constructed by the L--S method. Namely, we would like to
know how the genus $k$ of subsets involved in the minimax
procedure (\ref{ck1}) can be attributed to the ``geometry" of the
critical point $v_k(y)$ obtained in such a manner. In this
discussion, we assume to explain how to  merge the L--S genus
variational aspects with the actual practical structure of
``essential zeros and extrema" of basic patterns
 $\{F_l\}$. Recall that in the second-order case $m=1$, $N=1$
this is easy: by Sturm's Theorem, the genus $l$, which can be
formally ``attributed" to the function $F_l$, is equal to the
number of zeros (sign changes) $l-1$ or the number $l$ of isolated
local extrema points. Though, even for $m=1$, this is not that
univalent: there are other structures that do not obey the
Sturmian order (think about the solution via gluing
$\{F_0,F_0,...,F_0\}$ without sign changes); see more comments
below.

 For $m \ge 2$, this question is
 more difficult, and seems does not admit a clear rigorous
treatment. Nevertheless, we will try to clarify some its aspects.

Given a solution $F$ of (\ref{S2NN}) (a critical point of
(\ref{V1})), let us calculate the corresponding critical value
$c_F$ of (\ref{f4}) on the set (\ref{f2}), by taking
  \beq
 \label{dd1}
 \begin{matrix}
  \mbox{$
  v= C F \in {\mathcal H}_0 \quad \Longrightarrow \quad
   C = \frac 1{\sqrt{ - \int |\tilde D^m F|^2 + \int F^2}},
   $}
    \ssk\ssk\\
\mbox{so that}\quad
  \mbox{$
 c_F \equiv \tilde H(v) = \frac { \int |F|^\b}{( - \int |\tilde D^m F|^2 + \int
  F^2)^{ \b/ 2}} \quad \bigl(\b = \frac{n+2}{n+1} \bigr).
   $}
   \end{matrix}
   \eeq
This formula is important in what follows.

\ssk

\noi\underline{\em Genus one}. As usual in many variational
elliptic problems, the first pattern $F_0$ (typically, called a
{\em ground state}) is always of the simplest geometric shape, is
radially symmetric, and is a localized profile such as those in
Figure \ref{G1}. Indeed, this simple shape with a single dominant
maximum is associated with the variational formulation for $F_0$:
 \beq
 \label{FF81}
  \mbox{$
  F_0=r(v_0)v_0, \quad \mbox{with $v_0:$} \quad {\rm inf} \, \tilde H(v) \equiv {\rm inf} \,
   \int
  |v|^\b, \,\,\, v \in {\mathcal H}_0.
   $}
   \eeq
   This is (\ref{ck1}) with the simplest choice of closed
   sets as points, ${\mathcal F}=\{v\}$.

   Let us illustrate why a localized pattern like $F_0$ delivers the
   minimum to $\tilde H$ in (\ref{FF81}).
Take e.g. a two-hump structure,
 $$
 \hat v(y)=C\big[v_0(y) +
v_0(y+a)\big], \quad C \in \re,
 $$
  with sufficiently large $|a|\ge {\rm diam} \, {\rm supp} \, F_0$,
 so that supports of these two functions do not overlap.
 Then, evidently,
 $
  \hat v \in {\mathcal H}_0$ implies that
  $C = \frac
  1{\sqrt 2},
   $
   and, since $\b \in (1,2)$,
    $$
 \tilde H(\hat v) = 2^{\frac{2-\b}2} \tilde H(v_0) >  \tilde H(v_0).
  $$
By a similar reason, $F_0(y)$ and $v_0(y)$ cannot have ``strong
nonlinear oscillations" (see next sections for related concepts
developed in this direction), i.e., the positive part $(F_0)_+$
must be dominant, so that the negative part $(F_0)_-$ cannot be
considered as a separate dominant 1-hump structure. Otherwise,
deleting it will diminish $\tilde H(v)$ as above. In other words,
essentially non-monotone patterns such as in Figure \ref{G4} or
\ref{G6} cannot correspond to the 
 variational
problem (\ref{FF81}), i.e., the genus of the functional sets
involved is $\rho=1$.

Radial symmetry of $v_0$ is often standard in elliptic theory,
though is not straightforward at all in view of the lack of the
Maximum Principle and moving plane/spheres tools based on
Aleksandrov's Reflection Principle.
 We
just note that small non-radial deformations of this structure,
$v_0 \mapsto \hat v_0$ will more essentially affect (increase) the
first differential term $\int | \tilde D^m \hat v_0|^2$ rather
than the second one in the formula for $C$ in (\ref{dd1}).
Therefore, a standard scaling to keep this function in ${\mathcal
H}_0$ would mean taking $C \hat v_0$ with a constant
 $C>1$. Hence,
 $$
 \tilde H (C \hat v_0) = C^\b  \tilde H (\hat v_0) \approx C^\b  \tilde H ( v_0)
 >  \tilde H ( v_0),
  $$
  so infinitesimal non-radial perturbations do not provide us with
  critical points of (\ref{FF81}).

For $N=1$, this shows that $c_1$ cannot be attained at another
``positively dominant" pattern $F_{+4}$, with a  shape shown in
Figure \ref{G8}(a). See Table 1 below, where
 for $n=1$,
 $$
 c_{F_{+4}}= 1.9488...> c_2= c_{F_{1}} =1.8855...>c_1=  c_{F_{0}}=1.6203... \, .
  $$


\ssk

\noi\underline{\em Genus two}. Let now again for simplicity $N=1$,
and let $F_0$ obtained above for the genus $\rho=1$ be a simple
compactly supported pattern  as in Figure \ref{G1}.
 By $v_0(y)$ we denote the corresponding critical point given by
 (\ref{FF81}). We now take  the function corresponding to the difference
  (\ref{F0-}),
   \beq
   \label{v0-}
    \mbox{$
\hat v_2(y)=  \frac 1{\sqrt 2} \, \big[-v_0(y-y_0) + v(y+y_0)\big]
\in {\mathcal H}_0 \quad \bigl({\rm supp} \, v_0=[-y_0,y_0]\bigr),
 $}
 \eeq
which approximates the basic profile $F_1$ given in Figure
\ref{G4}.
One can see that
 \beq
 \label{HH12}
 \tilde H(\hat v_2)=2^{\frac{2-\b}2} \tilde H(v_0)= 2^{\frac{2-\b}2} c_1,
 \eeq
 so that, by (\ref{ck1}) with $k=2$,
  \beq
  \label{c21}
 c_1 < c_2 \le  2^{\frac{2-\b}2} c_1.
  \eeq


On the other hand, the sum as in (\ref{F0+}) (cf. Figure
\ref{G6}),
 \beq
 \label{ll21}
  \mbox{$
 \tilde
v_2(y)=  \frac 1{\sqrt 2} \, \big[v_0(y-y_0) + v(y+y_0)\big] \in
{\mathcal H}_0
 $}
 \eeq
 also delivers  the same value (\ref{HH12}) to the functional $\tilde
 H$.

\begin{figure}
 \centering
 \psfrag{F}{$F$}
 \psfrag{tF0}{$\sim F_0$}
 \psfrag{tF4}{$\sim F_{+4}$}
 \psfrag{tF222}{$\sim F_{+2,2,+2}$}
 \psfrag{tF1}{$\sim F_{1}$}
  \psfrag{v(x,t-)}{$v(x,T^-)$}
  \psfrag{final-time profile}{final-time profile}
   \psfrag{tapp1}{$t \approx 1^-$}
\psfrag{y}{$y$}
 \psfrag{0}{$0$}
 \psfrag{y}{$y$}
\includegraphics[scale=0.5]{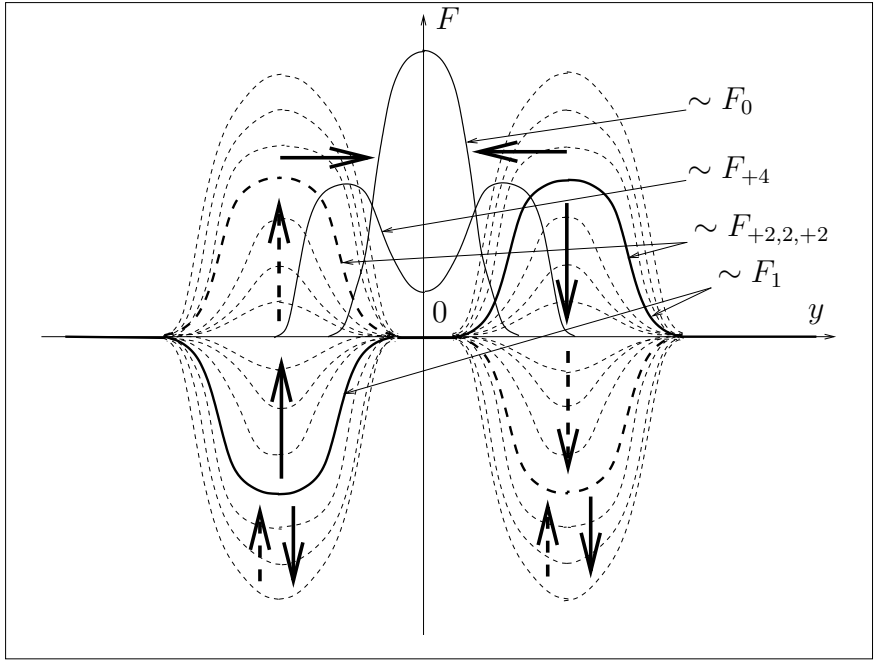}     
\caption{\small Profiles $F_1 \equiv F_{-2,1,+2} $, $F_{+2,2,+2}$,
$F_{+4}$, etc., can be embedded into a 1D subset in ${\mathcal
H}_0$ of genus two.}
 \label{Fgen2N}
\end{figure}

 It is easy to see that these patterns $F_1$ and $F_{+2,2,+2}$ as
 well as $F_{+4}$ and many others with two dominant extrema can be
  embedded into a 1D subset of genus two on ${\mathcal H}_0$. We show such a schematic
  picture if Figure \ref{Fgen2N}. Arrows there indicate the
  directions of deformations of patterns  on ${\mathcal H}_0$
  that can lead to any other profile from such a family.

  embedded into a 1D subset of genus two on ${\mathcal H}_0$.

\begin{figure}
 \centering
\includegraphics[scale=0.75]{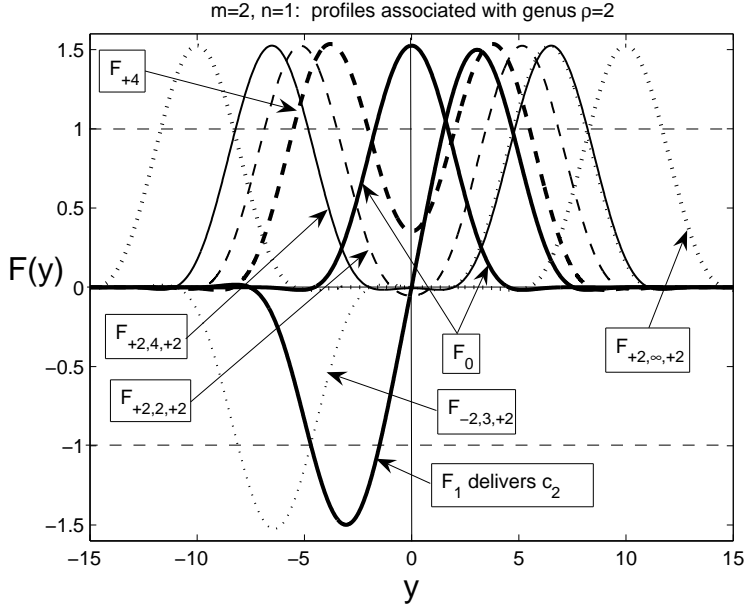}
 \vskip -.4cm
\caption{\rm\small Seven  patterns $F(y)$ indicated in Table 1;
$m=2$ and $n=1$.}
 \label{Fgen22}
\end{figure}

  It seems that, with such a huge, at least, countable variety of similar patterns,
we firstly  distinguish  the profile that delivers the critical
value $c_2$ given by (\ref{ck1}) by comparing the values
(\ref{dd1}) for each  pattern. The results are
presented in Table 1 
  for $n=1$, for which the critical values
(\ref{dd1}) are
 \beq
 \label{dd11}
  \mbox{$
 c_F= \tilde H(C F) = \frac { \int |F|^{3/2}}{( - \int |\tilde D^m F|^2 + \int
  F^2)^{3/4}}\quad \bigl(\b = \frac{3}{2} \bigr).
   $}
   \eeq
The corresponding profiles are shown in Figure \ref{Fgen22}.
Calculations have been performed with the enhanced values
Tols\,$=\e=10^{-4}$. Comparing the critical values in Table 1, we
thus arrive at the following conclusion based on this
analytical-numerical evidence:
 \beq
 \label{c211}
 \mbox{for genus $k =2$,
 the S--L critical
value $c_2=1.855...$ is delivered by $F_1$}.
 \eeq
Notice that the critical values $c_F$ for $F_1$ and $F_{+2,2,+2}$
are close  by just two percent.

\begin{table}[h]
\caption{Critical values  of $\tilde H(v)$; genus two}
\begin{tabular}{@{}lll}
 $F$ & $c_F$
 \\\hline
 $F_0$ & $1.6203...=c_1$\\ $F_1$ & $1.8855...=c_2$
\\ $F_{+2,2,+2}$ & $1.9255...$\\
 $F_{-2,3,+2}$ & $1.9268...$\\
$F_{+2,4,+2}$ & $1.9269...$
\\
$F_{+2,\infty,+2}$ & $1.9269...$\\
 $F_{+4}$ & $1.9488...$\\
\end{tabular}
\end{table}


 Thus, by  Table 1 the second critical value $c_2$ is achieved
 at the 1-dipole solution $F_1(y)$ having  the  transversal
zero at $y=0$, i.e., without any part of the oscillatory tail for
$y \approx 0$. Therefore, the neighbouring profile $F_{-2,3,+2}$
(see the dotted line in Figure \ref{Fgen22}) which has a small
remnant of the oscillatory tail  (see details in Section
\ref{Sect2}) with just 3 extra zeros, delivers another, worse
value
 $$
 c_F=1.9268... \quad \mbox{for} \quad F=F_{-2,3,+2}.
  $$
In addition, the lines from second to fifth in Table 1 clearly
show how $\tilde H$ increases with  the number of zeros in between
the $\pm F_0$-structures involved.

\ssk

\noi{\bf Remark: even for $m=1$ profiles are not variationally
recognizable.}
 Recall that for $m=1$, i.e., for the
ODE (\ref{4.4}), the $F_0(y)$ profile is not oscillatory at the
interface, so that the future rule (\ref{comp1}) fails. This does
not explain the difference between $F_1(y)$ and, say,
$F_{+2,0,+2}$, which, obviously deliver the same critical S--L
values by (\ref{ck1}). This is the case where we should
conventionally attribute the S--L critical point to $F_1$. Of
course, for $m=1$, existence of profiles $F_l(y)$ with precisely
$l$ zeros (sign changes) and $l+1$ extrema follows from Sturm's
Theorem.



Checking the accuracy of numerics and using (\ref{HH12}), we take
the critical values in the first and the fifth lines in Table 1 to
get for the profile $F_{+2, \infty,+2}$, consisting of two
independent $F_0$'s, to within $10^{-4}$,
 $$
c_F= 2^{\frac{2-\b}2} \tilde H(v_0)=2^{\frac 14} c_1= 1.1892...
\times 1.6203...=1.9269... \, .
 $$

\ssk

\begin{figure}
 \centering
\includegraphics[scale=0.75]{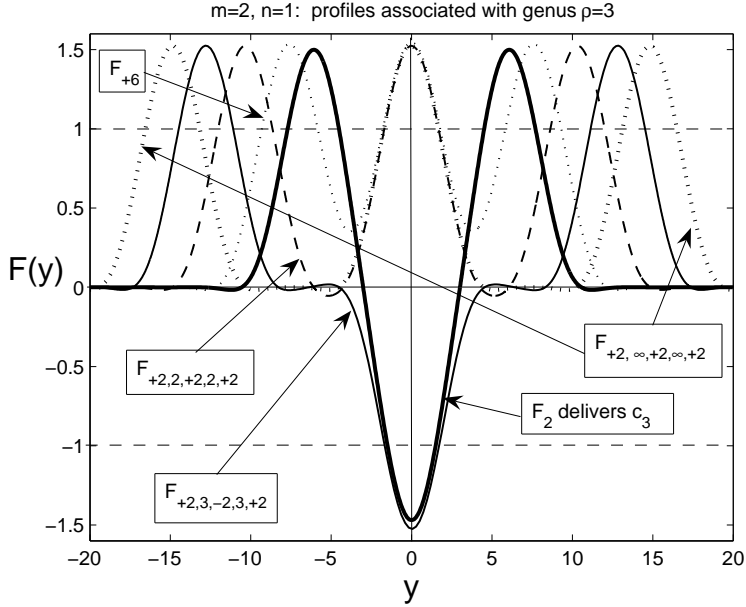}
 \vskip -.4cm
\caption{\rm\small Five patterns $F(y)$ indicated in Table 2;
$m=2$ and $n=1$.}
 \label{Fg3NN}
\end{figure}

\noi\underline{\em Genus three}. Similarly, for $k=3$  (genus
$\rho =3$), there are also several patterns that can pretend to
deliver the L--S critical value $c_3$. These are shown in Figure
\ref{Fg3NN}.

The corresponding critical values (\ref{dd11}) for $n=1$ are shown
in Table 2, which allows to conclude as follows:
 \beq
 \label{c2113}
 \mbox{$k =3$:
 the S--L critical
value $c_3=2.0710...$ is again delivered by the basic $F_2$}.
 \eeq
All critical values in Table 2 are very close to each other.
 Again
checking the accuracy of numerics and
   taking
the critical values $c_1$ in Table 1 and $c_F$ for $F_{+2,
\infty,+2, \infty,+2}$ in Table 2, consisting of three independent
$F_0$'s, yields, to within $10^{-4}$,
 $$
c_F= 3^{\frac{2-\b}2} \tilde H(v_0)= 3^{\frac 14} c_1= 1.31607...
\times 1.6203...=2.1324... \, .
 $$

\begin{table}[h]
\caption{Critical values  of $\tilde H(v)$; genus three}
\begin{tabular}{@{}ll}
$F$ & $c_F$
 \\\hline
 $F_2$ & $2.0710...=c_3$
\\ $F_{+2,2,+2,2,+2}$ & $2.1305...$\\
  $F_{+2,3,-2,3,+2}$ & $2.1322...$\\
 $F_{+2,\infty,+2, \infty,+2}$ & $2.1324...$
\\
 $F_{+6}$ & $2.1647...$\\
\end{tabular}
\end{table}


Note that the S--L category-genus construction (\ref{ck1}) itself
guarantees that all solutions $\{v_k\}$ as critical points will be
(geometrically) distinct; see \cite[p.~381]{KrasZ}. Here we stress
upon two important conclusions:


 \noi (I) First,  what is key for us, is that closed subsets
in ${\mathcal H}_0$ of functions of the sum type  in (\ref{ll21})
{\em  do not deliver
 S--L critical values in}
(\ref{ck1});


 \noi (II) On the other hand, patterns of the
$\{F_0,F_0\}$-interaction, i.e., those associated with the sum
structure as in (\ref{ll21}) {\em do exist}; see Figures \ref{G6}
for $m=2$;
 and


 \noi (III) Hence, these patterns (different from the basic ones $\{F_l\}$)
  as well as many others {\em are not
obtainable by a direct S--L approach.} Therefore, we will need
another version of S--L and fibering theory, with different type
of decomposition of functional spaces to be introduced below.

\ssk

\noi\underline{\em Genus $k$}. Similarly taking a proper sum of
shifted and reflected functions $\pm v_0(y \pm l y_0)$, we obtain
from (\ref{ck1}) that
 \beq
 \label{ck11}
c_{k-1} < c_k \le  k^{\frac{2-\b}2} c_1.
 \eeq


\ssk

\noi\underline{\em Conclusions: conjecture and an open problem}.
As a conclusion, we mention that, regardless such close critical
values $c_F$, the above numerics confirm that there is a
geometric-algebraic way to distinguish the S--L patterns
delivering (\ref{ck1}). It can be
 seen from (\ref{dd11}) that, destroying  internal oscillatory
``tail", or even any two-three zeros
 between two $F_0$-like patterns in the
complicated pattern $F(y)$,
 \beq
 \label{hh1}
  \mbox{$
\mbox{decreases two main terms
 $-\int |\tilde D^m F|^2$ and
  $\int|F|^{\frac 32}$ in $c_F$ in  (\ref{dd1}).
   }
   $}
   \eeq
Recall that
 precisely these terms in the ODE
  $$
  F^{(4)} = -|F|^{-\frac 12}F +...
 \quad (\mbox{see (\ref{4.3}) for $n=1$})
 $$
   are
 responsible for formation of the tail as shown in Section
 \ref{Sect2}, while the $F$-term, giving
 $\int F^2$,  is negligible in the tail.
  Decreasing both terms, i.e., eliminating the tail
 in between  $F_0$'s, will decrease the value $c_F$, since in (\ref{dd1})  the numerator
 gets less and the denominator larger.
 Therefore, composing a complicated
 pattern $F_l(y)$ from several elementary profiles looking like
 $F_0(y)$, by using $(k-1)$-dimensional manifolds of genus $k$,
  we  follow

\ssk

 \noi{\bf Formal Rule of
 Composition} (FRC) {\bf of patterns:} {\em performing  maximization of $\tilde H(v)$ of
 any $(k-1)$-dimensional manifold ${\mathcal F} \in {\mathcal
 M}_k$,}
  \beq
  \label{comp1}
\mbox{\em the S--L point $F_{k-1}(y)$ is obtained by minimizing
all internal tails and zeros,}
 \eeq
{\em  i.e., making the minimal number of internal transversal
zeros between single structures.
   }


Regardless such a simple variational-oscillatory meaning
(\ref{hh1}) of this FRC, we do not know how to give to such a rule
a rigorous sounding.


Concerning the actual critical S--L points, we
 end up with
  the
following conjecture,
 which well-corresponds to the
FRC (\ref{comp1}):

\ssk

\noi {\bf Conjecture \ref{SectVar}.4.} {\em For $N=1$ and any $m
\ge 2$, the critical S--L value $(\ref{ck1})$, $k \ge 1$, is
delivered by the basic pattern $F_{k-1}$, obtained by minimization
on the corresponding $(k-1)$-dimensional manifold ${\mathcal F}
\in {\mathcal M}_k$, which is
   the interaction
 \beq
 \label{f011}
 F_{k-1}=(-1)^{k-1}\{+F_0,-F_0,+F_0,...,(-1)^{k-1}F_0\},
  \eeq
where each neighbouring pair $\{F_0,-F_0\}$ or  $\{-F_0,F_0\}$ has
a single transversal zero in between the structures.
   }


 \ssk

We also would like to formulate the following conclusion that is
again associated with the specific structure of the  L--S
construction (\ref{ck1}) over suitable  subsets ${\mathcal F}$ as
smooth $(k-1)$-dimensional manifolds of genus $k$:

\ssk

\noi{\bf Open problem \ref{SectVar}.5.} {\em For $N=1$ and $m \ge
2$, there are no any purely geometric-topology arguments
establishing
that the conclusion of Conjecture $\ref{SectVar}.4$ holds.
Naturally, the same remains true in $\ren$.
  }

\ssk

In other words, we claim that the metric ``tail"-analysis of the
functionals involved in the FRC (\ref{comp1}) cannot be dispensed
with by any geometric-type arguments. Actually, the geometric
analysis is nicely applied for $m=1$ and this is perfectly covered
by Sturm's Theorem on zeros for second-order ODEs. If such a Sturm
Theorem is nonexistent, this emphasizes the end of
geometry-topology (or purely homotopy if tails are oscillatory)
  nature of the variational problem under consideration.







   \ssk

 \noi\underline{\sc On patterns in $\ren$}.
In the elliptic setting in $\ren$, such a clear picture of basic
patterns $F_l(y)$ obtained via (\ref{ck1}) with $k=l+1$ is not
available. As usual in elliptic theory, nodal structure of
solutions in $\ren$ is very difficult to reveal. Nevertheless, we
strongly believe that the L--S minimax approach (\ref{ck1}) also
can be used to detect the geometric shape of patterns in the
$N$-dimensional geometry.

For instance, it is most plausible  that the 1-dipole profile
$F_1(y)$, which {\em is not radially symmetric} is essentially
composed from two radial $F_0(y)$-type profiles via an
  $\{-F_0,F_0\}$-interaction ({\em q.v.} Figure \ref{G7} in 1D).
Therefore, $F_1(y)$ has two dominant extrema, in a natural way,
similar to the second eigenfunction
  \beq
  \label{He1}
   \mbox{$
  \frac {{\mathrm d}}{{\mathrm d}y_1}\,{\mathrm e}^{-|y|^2/4} \sim
 y_1 {\mathrm e}^{-|y|^2/4} \quad \mbox{in} \quad \ren
  $}
  \eeq
 of the self-adjoint second-order Hermite operator
 $$
  \mbox{$
 H_2=\D + \frac 12\,  y \cdot \n + \frac N2 \, I.
  $}
  $$
 Such a comparison assumes a bifurcation phenomenon
 at $n=0$ from eigenfunctions of the linearized operator;
 see applications
 in \cite[Appendix]{GHCo} to second-order quasilinear porous medium operators.
 Compactification of the pattern
(\ref{He1}) and making it oscillatory at the interface surface
would lead to
 correct understanding what the
1-dipole profile looks like, at least, for small $n>0$. Similar
analogy is
 developed for all odd patterns $F_{2k+1}$. For
instance, $F_3(y)$ has a dominant ``topology" similar to the
fourth eigenfunction of $H_2$,
 $$
  \mbox{$
 \frac {{\mathrm d}^3}
{{\mathrm d}y_1^3}\,{\mathrm e}^{-|y|^2/4} \sim
 \frac {y_1}2 \big( \frac {y_1^2}2-3 \big) {\mathrm e}^{-|y|^2/4},
  $}
  $$
  which has precisely four extrema on the $y_1$-axis, etc.

Concerning even patterns $F_{2k}$, we believe that the above L--S
algorithm leads to simple radially-symmetric solutions of
(\ref{S2NN}), i.e., solutions of ODEs; see further comments below.

\ssk

\noi {\bf Remark: on radial and 1D geometry.} Of course, the
elliptic equation (\ref{S2NN}) also admits a countable family of
radially symmetric solutions
 $$
 \{F_l^{\rm rad}(|y|), \, l=0,2,4,...\}
  $$
   satisfying
the corresponding ODE. These are constructed in a similar manner
by L--S and fibering theory. In 1D, this gives the basic set
$\{F_l^{\rm 1D}, \, l=0,1,2,...\}$ that was described in the
previous sections. We expect that the first members of all three
families coincide,
 \beq
 \label{f57}
 F_0=F_0^{\rm rad}=F_0^{\rm 1D}.
 \eeq
Further correspondence of the L--S spectrum of patterns in
Proposition \ref{Pr.MM} and the 1D one $\{F_l^{\rm 1D}, \, l \ge
0\}$ is discussed in \cite[\S~4]{GMPSobII}.

\section{\underline{\bf Problem ``Oscillations"}: local oscillatory structure of solutions
close to interfaces and  periodic connections with singularities}
 \label{Sect2}

As we have seen,
the first principal feature of the ODEs (\ref{S2}) (and the elliptic counterparts)
 is that these
admit compactly supported solutions. Indeed, all interesting for
us patterns have finite interfaces.
 This has been proved
in Proposition \ref{Pr.CS} in a general elliptic setting.

Therefore, we are going to study typical local behaviour of the
solutions of (\ref{S2}) close to the singular points, i.e., to
finite interfaces. We will reveal extremely oscillatory structure
of such behaviour to be compared with global oscillatory behaviour
obtained above by variational techniques.

The phenomenon of oscillatory changing sign behaviour of solutions
of the Cauchy problem has been detected for various classes of
evolution PDEs; see a general view in \cite[Ch.~3-5]{GSVR} and
 various results for different PDEs in
\cite{Gl4, GBl6, Galp1}.
   For the present $2m$th-order equations, the
oscillatory behaviour exhibits special features to be revealed.
 We expect that the
presented oscillation analysis makes sense for more general
solutions of the parabolic equation (\ref{S1}) and explains their
generic behaviour close to the moving interfaces.

\subsection{Autonomous ODEs for oscillatory components}

Assume that the finite interface of $F(y)$ is situated at the
origin $y=0$, so that we can use the trivial  extension $F(y)
\equiv 0$ for $y<0$. We then are interested in describing the
behaviour of solutions as $y \to 0^+$, so we consider the ODE
(\ref{S2}) written in the form
 \beq
 \label{2.1}
  \mbox{$
  F^{(2m)}=  (-1)^{m+1} |F|^{-\a} F + (-1)^{m}F\,\,\, \mbox{for} \,\, y>0, \,
\,\, F(0)=0
  \quad\big(\a=\frac n{n+1}\in (0,1)\big).
   $}
  \eeq
   In view of the scaling structure of the first two terms, for
   convenience, we perform extra rescaling and
   introduce the {\em oscillatory component} $\var(s)$ of $F$ by
   \beq
   \label{2.2}
    \mbox{$
   F(y) = y^\g \var(s) \whereA s= \ln y \andA
    \g = \frac{2m}\a \equiv \frac {2m(n+1)}n.
    $}
    \eeq
Therefore, since $s \to -\infty$ (the new interface position) as
$y \to 0^-$, the monotone function $y^\g$ in (\ref{2.2}) plays the
role of the {\em envelope} to the oscillatory function $F(y)$.
Substituting (\ref{2.2}) into (\ref{2.1}) yields the following
equation for $\var$:
 \beq
 \label{2.3}
 P_{2m}(\var) = (-1)^{m+1} |\var|^{-\a} \var +(-1)^m {\mathrm e}^{2m s}
 \var.
  \eeq
Here $\{P_{k}, \, k \ge 0\}$ are linear differential operators
defined by the recursion
 \beq
  \label{LH.7}
  P_{k+1}(\var)= (P_k(\var))' + (\g-k)P_k(\var) \,\,\, \mbox{for}
  \,\,\, k=0,1,... \, , \,\,\,
  P_0(\var)=\var.
   \eeq
 Let us present the first five operators, which are sufficient for
 further use:
 \begin{align}
 & P_1(\var)= \var' + \g \var; 
 \notag\\
 & P_2(\var)=\var''+(2\g-1) \var' + \g(\g-1)\var;
  \notag\\
 &P_3(\var)= \varphi''' + 3(\g-1) \varphi'' + (3 \g^2 - 6 \g
+2) \varphi'
  + \g(\g-1)(\g-2)\varphi;
   \notag\cr
& P_4(\var)= \var^{(4)} + 2(2\g -3) \var''' + (6 \g ^2- 18 \g +11)
\var''  \smallskip
 \cr
& \qquad \quad \,\,+  2(2 \g ^3
   - 9 \g ^2
    +  11 \g  -3) \var'
 + \g (\g -1)(\g -2)(\g -3) \var;\notag
  \ssk \\
&  P_5(\var)= \varphi^{(5)}+5(\g -2)\varphi^{(4)} +5
(2\g^2-8\g+7)\varphi'''
 + 5 (\g-2)(2\g^2-8\g +5)\varphi'' \notag \ssk\\
& \qquad \quad\,\, + (5\g^4-40\g^3+105 \g^2 - 100\g+24)\varphi' +
\, \g(\g-1)(\g-2)(\g-3)(\g-4)\varphi;
  \,\,\,\mbox{etc.}
 \notag
 \end{align}

According to (\ref{2.2}), the interface at $y=0$ now corresponds
to $s=-\infty$, so that (\ref{2.3}) is an exponentially (as $s \to
-\infty$) perturbed autonomous ODE
 \beq
 \label{2.4}
 \mbox{$
 P_{2m}(\var) = (-1)^{m+1} |\var|^{-\a} \var \quad \mbox{in} \quad
 \re
 \quad \big( \a= \frac n{n+1}\big),
 $}
  \eeq
  which we will concentrate upon.
   By
 classic ODE theory
  \cite{CodL}, one can expect that for $s \ll -1$ typical
 (generic)  solutions of (\ref{2.3}) and (\ref{2.4}) asymptotically
differ by exponentially small factors. Of course, we must admit
that (\ref{2.4}) is a singular ODE with a non-Lipschitz term, so
the results on continuous dependence need extra justification in
general.

Thus, in two principal cases, the ODEs for the oscillatory
component $\var(s)$ are
\begin{align}
 m=2: \quad P_4(\var)= - |\var|^{-\a} \var, \label{m2} \\
 m=3: \quad P_6(\var)= + |\var|^{-\a} \var, \label{m3}
  \end{align}
 that exhibit rather different properties because comprise even and odd $m$'s.
For instance, (\ref{2.4}) for any odd $m \ge 1$
 (including (\ref{m3}))
has two constant equilibria, since
 \beq
 \label{2.5}
 \begin{matrix}
\g(\g-1)...(\g-(2m-1)) \var= |\var|^{-\a}
  \var \quad
  \Longrightarrow \ssk \\
\var(s)= \pm \var_0 \equiv \pm [ \g(\g-1)...(\g-(2m-1))]^{-\frac 1
\a}
 \quad \mbox{for all \, $n>0$}.
 \end{matrix}
  \eeq
For even $m$ including (\ref{m2}), such equilibria for (\ref{2.4})
do not exist at least for  $n \in (0,1]$. We will show how this
affects the oscillatory properties of solutions for odd and even
$m$'s.

\subsection{Periodic oscillatory components}

We now look for {\em periodic} solutions of (\ref{2.4}), which are
the simplest nontrivial bounded solutions that can be continued up
to the interface at  $s= - \infty$. Periodic solutions together
with their stable manifolds are  simple {\em connections} with the
interface, as a singular point of the ODE (\ref{S2}).

Note that (\ref{2.4}) does not admit variational setting, so we
cannot apply well developed  potential \cite[Ch.~8]{MitPoh} (see a
large amount of related existence-nonexistence results and further
 references therein),  or degree \cite{Kras, KrasZ} theory.
For $m=2$, the proof of existence of $\varphi_*$ can be done by
shooting; see \cite[\S~7.1]{Gl4} that can be extended to $m=3$ as
well. Nevertheless, uniqueness of a periodic orbit is still open,
so we conjecture the following result supported by various
numerical and analytical evidence (cf. \cite[\S~3.7]{GSVR}):

\ssk

 \noi {\bf Conjecture \ref{Sect2}.1.}  {\em
For any $m \ge 2$ and $\a \in (0,1]$, the ODE $(\ref{2.4})$ admits
a unique nontrivial periodic solution $\varphi_*(s)$ of changing
sign.}

















\subsection{Numerical construction of periodic orbits for $m=2$}

\begin{figure}
\centering \subfigure[$n=2$]{
\includegraphics[scale=0.52]{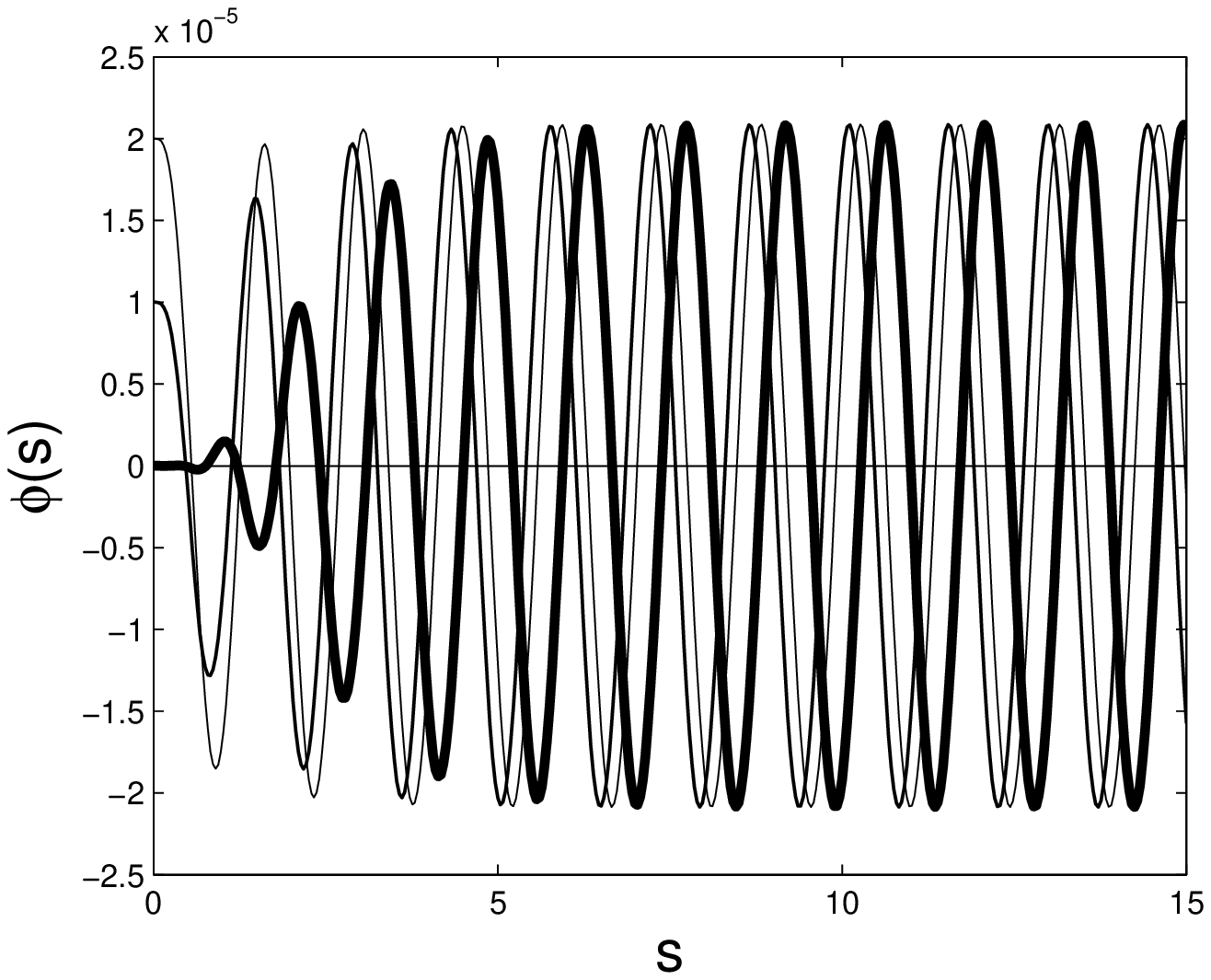}
} \subfigure[$n=4$]{
\includegraphics[scale=0.52]{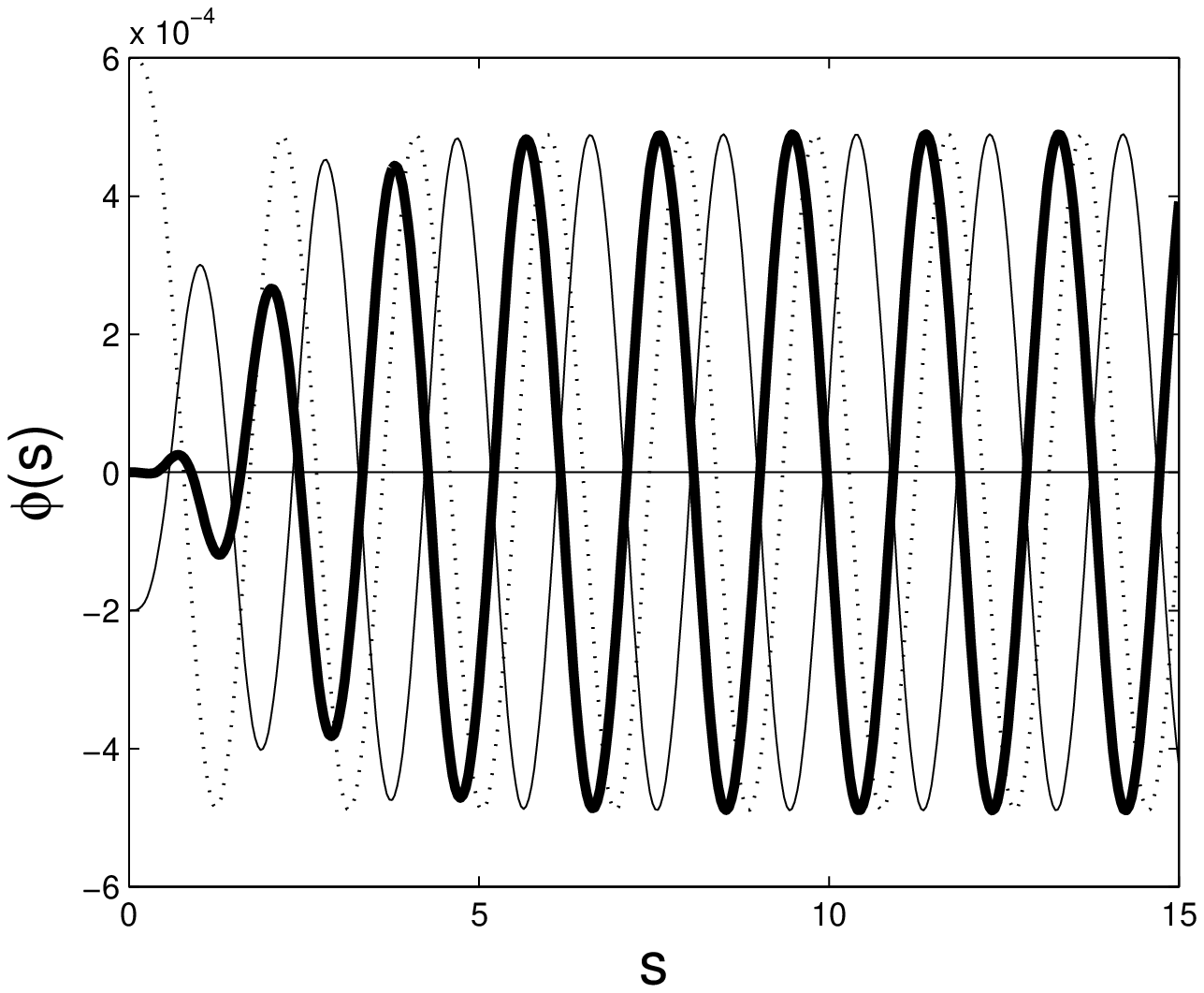}
}
 \vskip -.2cm
\caption{\rm\small Convergence to the stable periodic solution of
(\ref{m2})  for $n=2$ (a) and $n=4$ (b).}
 \label{F1}
\end{figure}

Numerical results clearly suggest that (\ref{m2}) possesses a
unique periodic solution $\var_*(s)$, which is stable in the
direction opposite to the interface, i.e., as $s \to +\infty$; see
Figure \ref{F1}. The proof of exponential stability and
hyperbolicity of $\var_*$ is straightforward by estimating the
eigenvalues of the linearized operator.
This  agrees with the obviously correct similar result for $n=0$,
namely, for the linear equation (\ref{2.1}) for $\a=0$
 \beq
 \label{2.9}
 F^{(4)}=-F \quad \mbox{as} \quad y \to -\infty.
  \eeq
Here the interface is infinite, so its position corresponds to
$y=-\infty$. Indeed, setting $F(y)={\mathrm e}^{ \mu y}$ gives the
characteristic equation and a unique exponentially decaying
pattern
 \beq
 \label{2.10}
 \mu^4=-1 \quad \Longrightarrow \quad
  \mbox{$
  F(y) \sim 
 {\mathrm e}^{ \frac y {\sqrt 2} } \,\,\bigl[A \cos\bigl(\frac y {\sqrt 2}
 \bigr)
  +B \sin\bigl(\frac y {\sqrt 2}\bigr)\bigr]
 \,\,\, \mbox{as} \,\,\, y \to -\iy
  . $}
  \eeq
Continuous dependence on $n \ge 0$ of typical solutions of
(\ref{2.4}) with ``transversal" zeros only will continue to be key
in our analysis, that actually means existence of a ``homotopic"
connection between the nonlinear  and the linear ($n=0$)
equations. The passage to the limit $n \to 0$ in similar
degenerate ODEs from thin film equations (TFEs) theory is
discussed in \cite[\S~7.6]{Gl4}.


 The oscillation
amplitude becomes very
 small for $n \approx 0$, so we perform
extra scaling.

\smallskip

\begin{figure}
\centering \subfigure[$n=0.2$]{
\includegraphics[scale=0.52]{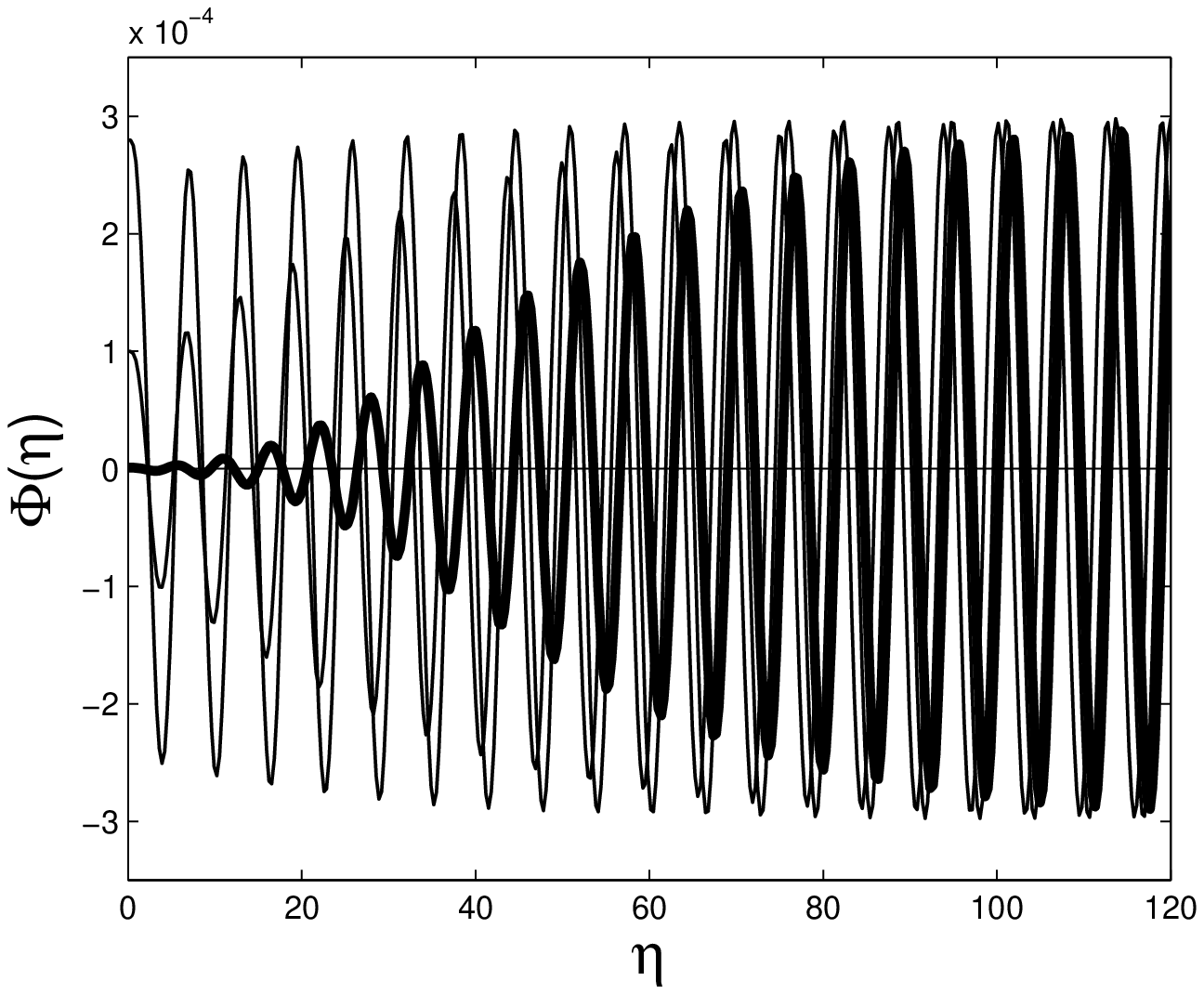}
} \subfigure[$n=0.08$]{
\includegraphics[scale=0.52]{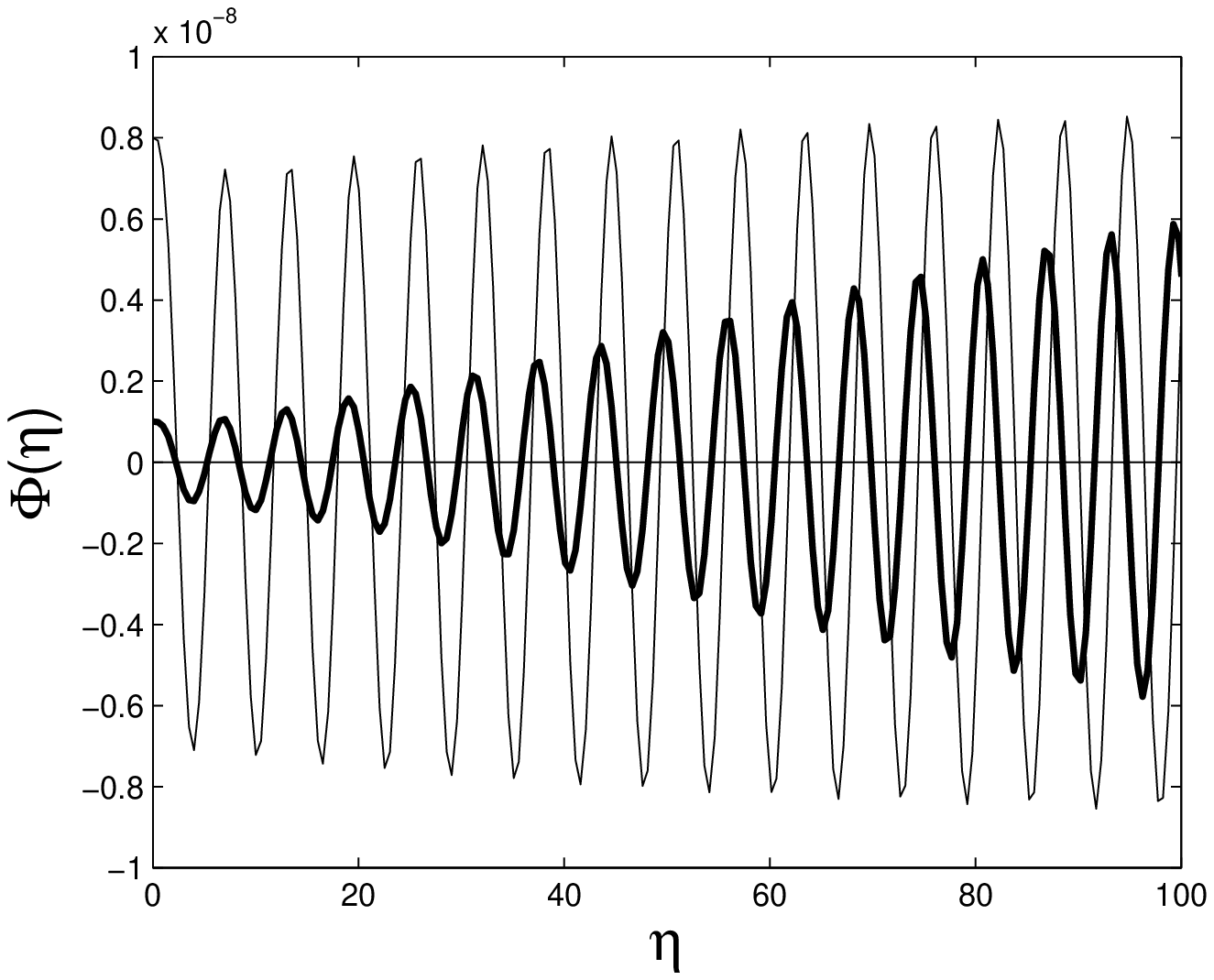}               
}
 \vskip -.2cm
\caption{\rm\small Stable periodic oscillations  in the ODE
(\ref{GGGF.1}) for $n=0.2$ (a) and $n=0.08$ (b).}
 \label{F2}
\end{figure}

\noi\underline{\em Limit ${n \to 0}$}. This scaling is
 \beq
 \label{FGF.222}
 \mbox{$
 \var(s)= \bigl(\frac n4 \bigr)^{\frac 4n} \Phi(\eta), \quad \mbox{where} \quad
 \eta= \frac {4 s}n,
  $}
  \eeq
  where $\Phi$ solves a simpler {\em limit} binomial ODE,
   \beq
   \label{GGGF.1}
   \mbox{$
   {\mathrm e}^{-\eta}({\mathrm e}^{\eta}\Phi)^{(4)} \equiv
\Phi^{(4)}+ 4 \Phi'''+6\Phi''+4 \Phi'+ \Phi=
-\big|\Phi\big|^{-\frac n{n+1}} \Phi.
 $}
 \eeq
 The stable oscillatory patterns of (\ref{GGGF.1})
   are shown in Figure
 \ref{F2}.
  For such small $n$ in Figure \ref{F2}(a) and (b), by
scaling (\ref{FGF.222}),  the periodic components $\varphi_*$ get
really small, e.g.,
 $$
 \begin{matrix}
\max |\var_*(s)| \sim 3 \cdot  10^{-4} \bigl(\frac n4
\bigr)^{\frac 4n} \sim 3 \cdot  10^{-30} \quad \mbox{for \,\,
$n=0.2$ \,in \,(a)}, \ssk \ssk\\
\mbox{and} \quad
  \max |\var_*(s)|  \sim 10^{-93}\quad \mbox{for
\,\, $n=0.08$ \,in
  \,(b)}. 
   \end{matrix}
   $$



\smallskip

\begin{figure}
 \centering
\includegraphics[scale=0.6]{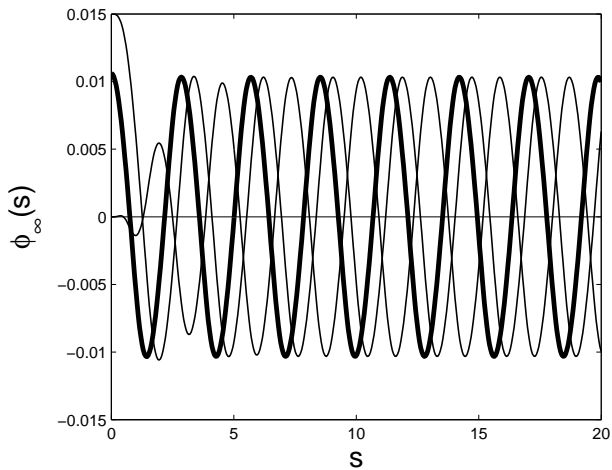}
 \vskip -.4cm
\caption{\rm\small Convergence to the stable periodic solution of
(\ref{GGGF.2}) ($n= +\infty$).} \vskip -.2cm
 \label{Finf1}
\end{figure}

 \noi\underline{\em Limit ${ n \to \infty}$}. Then 
$\a \to 1$, so the original ODE (\ref{m2})  approaches the
following equation with a discontinuous sign-nonlinearity:
 \beq
   \label{GGGF.2}
\var_\infty^{(4)}+ 10 \var_\infty'''+35 \var_\infty''+50
\var_\infty'+ 24 \var_\infty=- {\rm sign} \, \var_\infty.
 \eeq
 This also admits a stable periodic solution, as shown in Figure
 \ref{Finf1}.


\subsection{Numerical construction of periodic orbits for $m=3$}

 Consider now equation (\ref{m3}) that admits constant equilibria
(\ref{2.5}) existing for all $n>0$. It is easy to check that the
equilibria $\pm \var_0$ are asymptotically stable as $s \to
+\infty$. Then the necessary periodic orbit is  situated in
between of these stable equilibria, so it is unstable as $s \to
+\infty$.

\begin{figure}
 \centering
\includegraphics[scale=0.65]{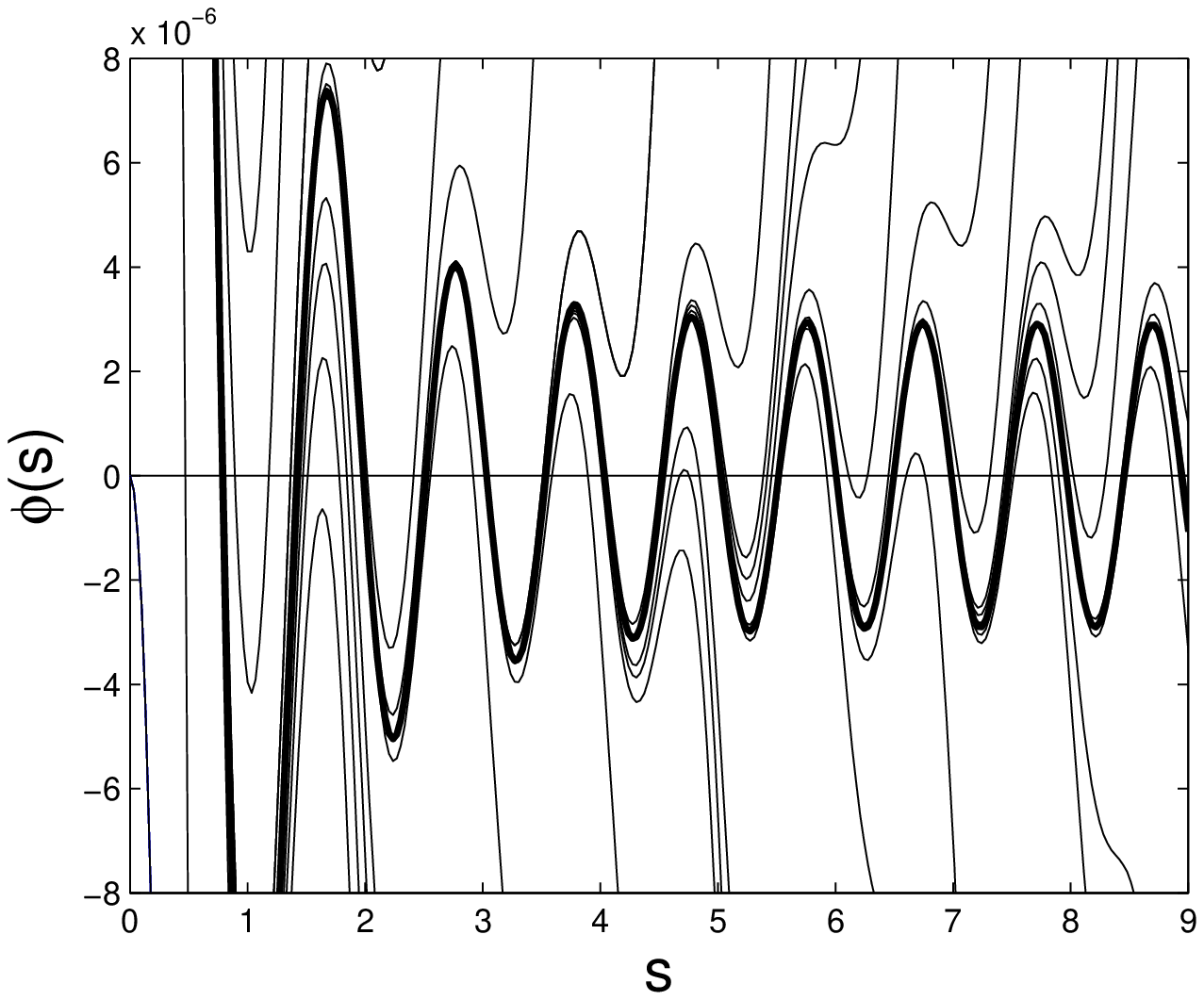}
 \vskip -.4cm
\caption{\rm\small Unstable periodic behavior of the ODE
(\ref{m3})  for $n= 15$. Cauchy data are given by
$\var(0)=10^{-4}$, $\var'(0)=\var'''(0)=...=\var^{(5)}(0)=0$, and
 $\var''(0)=-5.0680839826093907...\times 10^{-4}$ .}
   \vskip -.3cm
 \label{FUn15}
\end{figure}

 Such an unstable periodic solution of (\ref{m3}) is
 shows in Figure \ref{FUn15} for $n=15$, which is obtained by shooting
 from $s=0$ with prescribed Cauchy data.


\begin{figure}
 \centering
\includegraphics[scale=0.65]{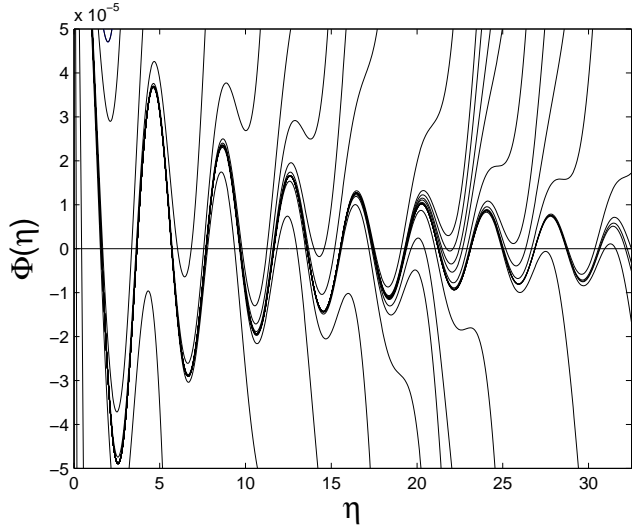}
 \vskip -.4cm
\caption{\rm\small Unstable periodic behavior of the ODE
(\ref{Sc.n.2}) for $n= \frac 12$. Cauchy data are given by
$\var(0)=10^{-4}$, $\var'(0)=\var'''(0)=...=\var^{(5)}(0)=0$, and
$\var''(0)=-9.456770333415...\times 10^{-4}$ .}
   \vskip -.3cm
 \label{FNN1}
\end{figure}

As for $m=2$,  in order to reveal periodic oscillations for
smaller $n$ (actually, there is a numerical difficulty already for
$n \le 4$), we apply the scaling
 \beq
\label{Sc.n.1} \mbox{$ \var(s)= \bigl(\frac n6 \bigr)^{\frac 6n}
\Phi(\eta), \quad \mbox{where} \quad\eta= \frac {6 s}n.
 $}
 \eeq
This gives in the limit a simplified ODE with the binomial linear
operator,
 \beq
\label{Sc.n.2}
 {\mathrm
e}^{-\eta}({\mathrm e}^\eta \Phi)^{(6)} \equiv \Phi^{(6)} + 6
\Phi^{(5)} + 15  \Phi^{(4)} +20 \Phi''' + 15
 \Phi'' + 6 \Phi' + \Phi
  = \big|\Phi \big|^{-\frac
n{n+1}}\Phi.
 \eeq
Figure \ref{FNN1} shows the trace of the periodic behaviour for
equation (\ref{Sc.n.2}) with $n=\frac 12$.
 According to scaling (\ref{Sc.n.1}),
the periodic oscillatory  component $\var_*(s)$ gets very small,
 $$
 \max |\var_*| \sim 1.1 \times 10^{-18}  \,\,\, \mbox{for} \,\,\,
n=0.5.
 $$
A more detailed study of the behaviour of the oscillatory
component as $n \to 0$ is available in \cite[\S~12]{GBl6}.

The passage to the limit $n \to +\infty$ leads to the equation
with discontinuous nonlinearity that is easily obtained from
(\ref{m3}). This  admits a periodic solution, which is rather
close to the periodic  orbit in Figure \ref{FUn15} obtained for
$n=15$.



We claim that the above two cases $m=2$ (even) and $m=3$ (odd)
exhaust all key types of periodic behaviours in ODEs like
(\ref{S2}). Namely, periodic orbits are stable for even $m$ and
are unstable for odd, with typical stable and unstable manifolds
as $s \to \pm \infty$. So we observe a dichotomy relative to all
orders $2m$ of the ODEs under consideration.

\section{\underline{\bf Problem ``Numerics"}: numerical construction and first classification of
basic types  of localized blow-up or compacton patterns for $m=2$}
   \label{Sect4}

We need  a careful   numerical description of various families of
solutions of the ODEs (\ref{S2}). In practical computations, we
have to use the regularized version of the equations,
 \beq
 \label{4.1}
  \mbox{$
(-1)^{m} F^{(2m)}=F- 
 \big(\e^2+ F^2 \big)^{-\frac n{2(n+1)}} F \quad \mbox{in} \,\,\,
 \re,
 $}
 \eeq
 which, for $\e>0$, contains smooth analytic nonlinearities.
 In numerical analysis, we typically take $\e = 10^{-4}$ or, at least,  
$10^{-3}$ which is sufficient to revealing global structures.


It is worth mentioning that detecting in Section \ref{Sect2} a
highly oscillatory structure of solutions close to interfaces
makes it impossible to use  well-developed {\em homotopy} theory
\cite{KKVV00, VV02} that was successfully applied to another class
of fourth-order ODEs with coercive operators; see also
Peletier--Troy \cite{PelTroy}.
   Roughly speaking, our non-smooth problem
cannot be used in  a homotopy classification, since the
oscillatory behaviour close to interfaces destroys such a standard
 homotopy parameter as the number of rotations on the hodograph
plane $\{F,F'\}$. Indeed, for any solution of (\ref{S2}),
 the rotation number about the origin is always infinite.
Then as $F \to 0$, i.e., as $y \to \pm \infty$, the linearized
equation is (\ref{2.9}), which  admits the  oscillatory behaviour
 (\ref{2.10}).


\subsection{Fourth-order equation: $m=2$}

We will describe main families of solutions.

\ssk

\noi \underline{\em First basic pattern and structure of zeros}.
For $m=2$, (\ref{S2}) reads
 \beq
 \label{4.3}
  \mbox{$
 F^{(4)}=F- 
 \big|F\big|^{-\frac n{n+1}} F \quad \mbox{in} \,\,\,
 \re.
 $}
 \eeq
 We are looking for compactly supported patterns $F$ (see
Proposition \ref{Pr.CS}) satisfying
 \beq
 \label{4.3NN}
  \begin{matrix}
 {\rm meas} \, {\rm supp} \,\, F > 2 R_*,
 \quad \mbox{where $R_*> \frac \pi 2$ is the first} \ssk\ssk \\
 \mbox{ positive root
 of the equation \,\,$\tanh R=- \tan R$.}
  \end{matrix}
  \eeq

\begin{figure}
 \centering
\includegraphics[scale=0.8]{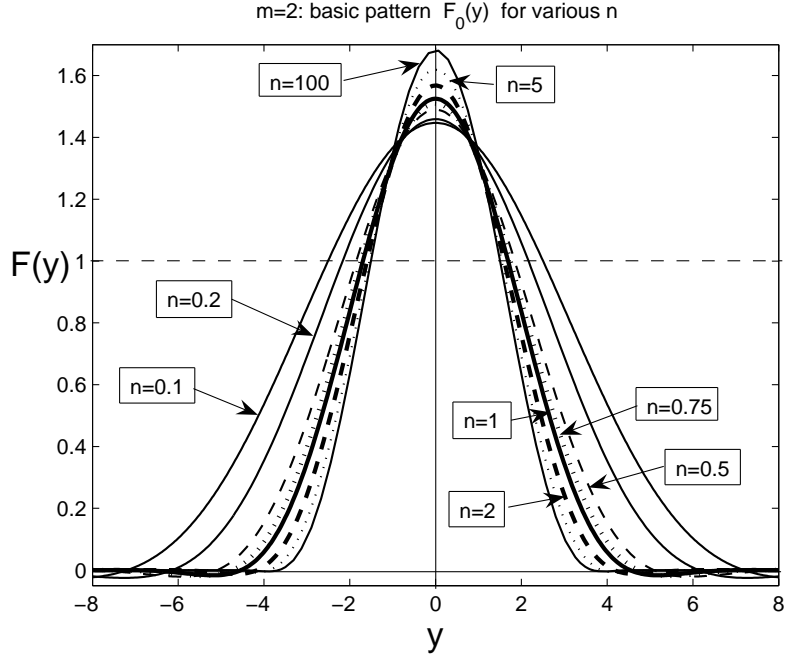} 
 \vskip -.4cm
\caption{\rm\small The first (stable) solution $F_0(y)$  of
(\ref{4.3}) for various $n$.}
   \vskip -.3cm
 \label{G1}
\end{figure}

 In Figure
\ref{G1}, we show the first basic 
pattern called the  $F_0(y)$ for various $n \in [\frac
1{10},100]$. Concerning the last profile $n=100$, note that
(\ref{4.3}) admits a natural passage to the limit $n \to +\infty$
that gives the ODE with a discontinuous nonlinearity,
 \beq
 \label{sign1}
 F^{(4)}= F- {\rm sign} \, F \equiv \left\{
 \begin{matrix}
 F-1 \quad\,\,\,\, \mbox{for} \quad F  > 0, \ssk\\
 F+1 \quad \mbox{for} \quad F <0.
  \end{matrix}
 \right.
  \eeq
A unique oscillatory solution of (\ref{sign1}) can be treated by
an algebraic approach; cf. \cite[\S~7.4]{Gl4}. For $n=1000$ and
$n=+\infty$, the profiles are close to that for $n=100$ in Figure
\ref{G1}.

\begin{figure}
\centering \subfigure[scale $10^{-3}$]{
\includegraphics[scale=0.52]{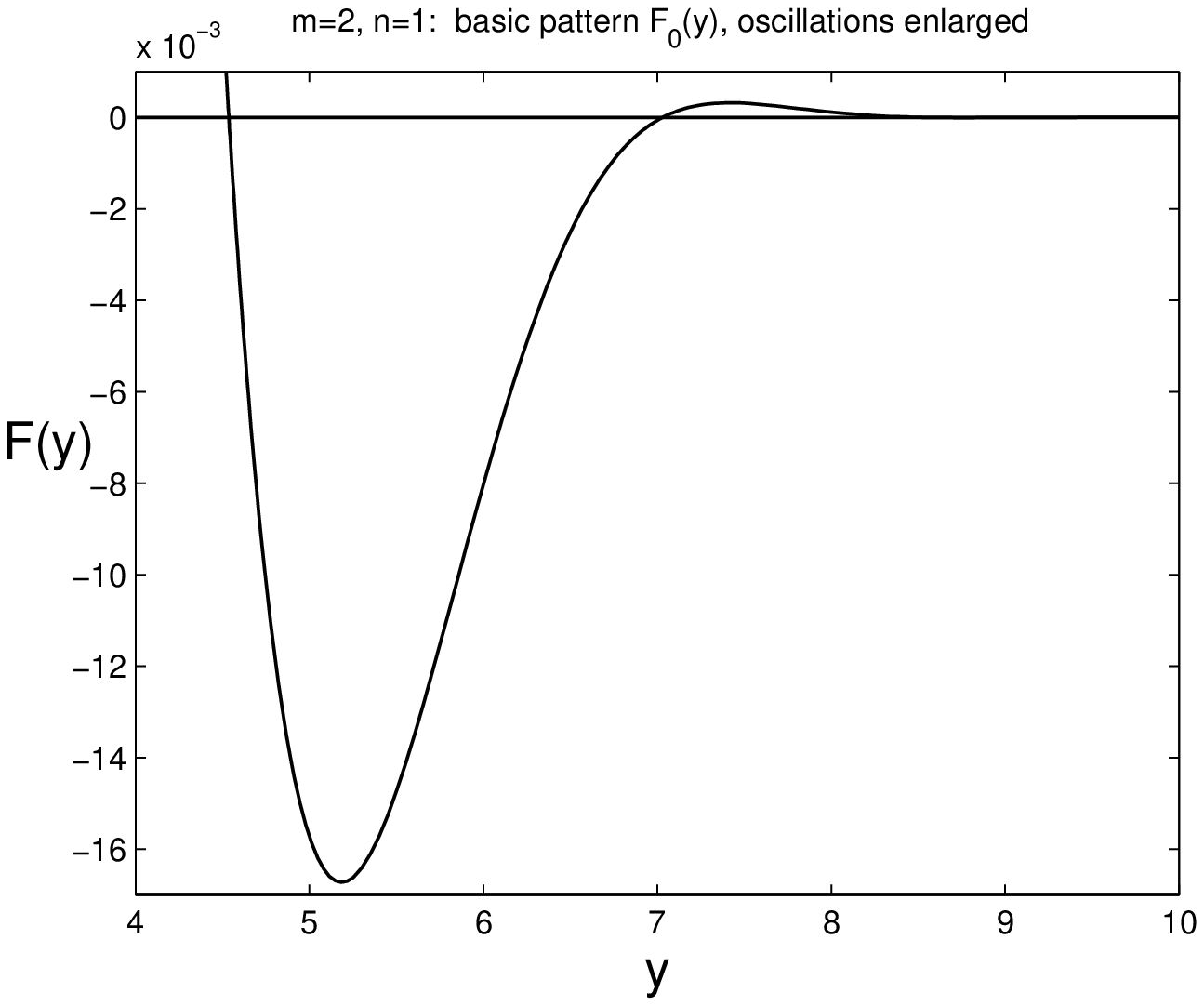}
} \subfigure[scale $10^{-4}$]{
\includegraphics[scale=0.52]{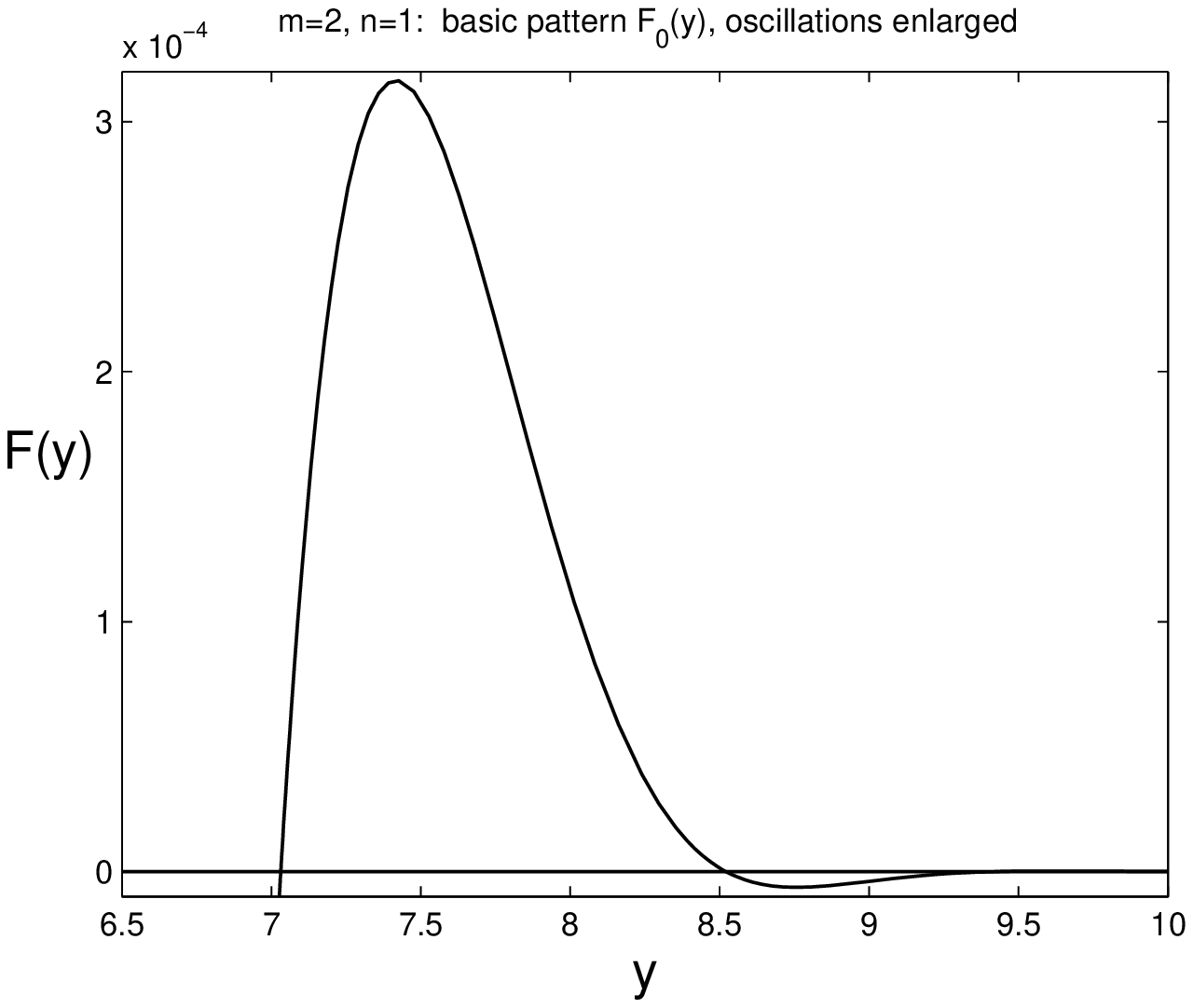}               
}
 \vskip -.2cm
\caption{\rm\small  Enlarged zero structure of the profile
$F_0(y)$ for $n=1$ in  the linear scale.}
 \label{G2}
\end{figure}

The profiles in Figure \ref{G1} are  constructed by MATLAB with
extra accuracy, where $\e$ in (\ref{4.1}) and both tolerances in
the {\tt bvp4c} solver have been enhanced and took the values
 $
  \e= 10^{-7}$ and ${\rm Tols}=10^{-7}.
 $
 This allows us also to check the refined local structure of
 multiple zeros at the interfaces. Figure \ref{G2} corresponding to $n=1$ explains how
the zero structure repeats itself  from one zero to another in the
usual linear scale.

\ssk

\begin{figure}
 \centering
\includegraphics[scale=0.8]{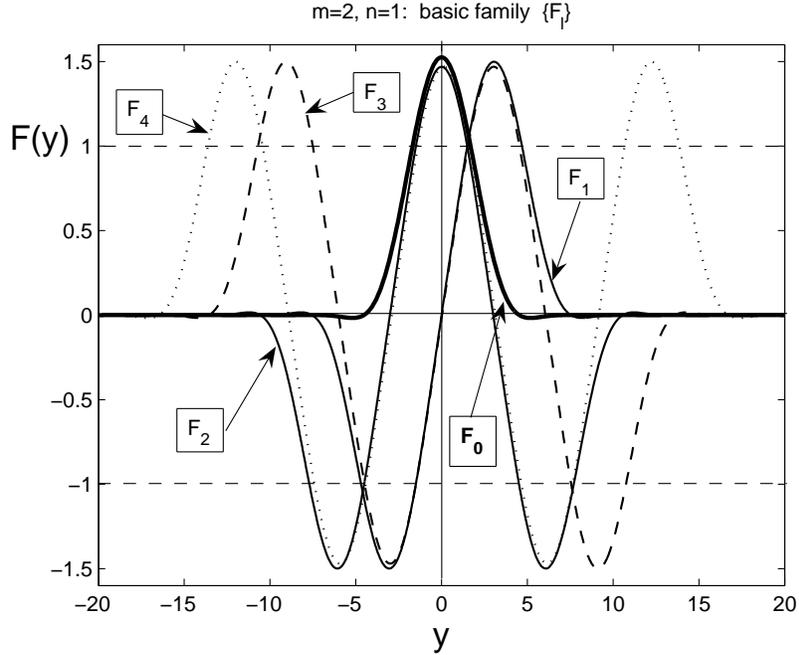}
 \vskip -.4cm
\caption{\rm\small The first five patterns of the basic family
$\{F_l\}$ of the ODE (\ref{4.3}) for $n=1$.}
   \vskip -.3cm
 \label{G4}
\end{figure}

\noi\underline{\em Basic countable family: approximate Sturm's
property}. In Figure \ref{G4}, we show the basic family denoted by
$\{F_l, \, l=0,1,2,...\}$ of solutions of (\ref{4.3}) for $n=1$.
Each profile $F_l(y)$ has precisely $l+1$ ``dominant" extrema and
$l$ ``transversal" zeros; see further discussion below and
\cite[\S~4]{GHUni} for other details. It is important that
 {\em all the internal zeros of $F_l(y)$ are clearly {
transversal}}
  (obviously,  excluding the oscillatory end points of the support).
  In other words, each profile $F_l$ is
approximately obtained by a simple ``interaction" (gluing
together) of $l+1$ copies of the first pattern $\pm F_0$ taking
with necessary signs; see further development below.

Actually, if we forget for a moment about the complicated
oscillatory structure of solutions near interfaces, where an
infinite number of extrema and zeros occur, the dominant geometry
of profiles in Figure \ref{G4} looks like it approximately obeys
Sturm's classic zero set property, which is true rigorously for
$m=1$ only, i.e., for the second-order ODE
 \beq
 \label{4.4}
  \mbox{$
 F''=-F + \big|F \big|^{-\frac n{n+1}}F \inB \re.
 $}
  \eeq
For (\ref{4.4}), the basic family $\{F_l\}$ is constructed by
direct gluing together the explicit patterns (\ref{RD.4}), i.e.,
$\pm F_0$. Therefore, each $F_l$ consists of precisely $l+1$
patterns (\ref{RD.4}) (with signs $\pm F_0$), so that Sturm's
property is clearly true.
 In Section \ref{SectVar}, we presented
some analytic evidence showing that precisely this basic family
$\{F_l\}$ is obtained
    by direct  application of L--S
    category theory.


 \subsection{Countable family of $\{F_0,F_0\}$-interactions}

 We now show that the actual nonlinear interaction of the two first
 patterns $+F_0(y)$ leads to a new family of profiles.

\begin{figure}
 \centering
\includegraphics[scale=0.7]{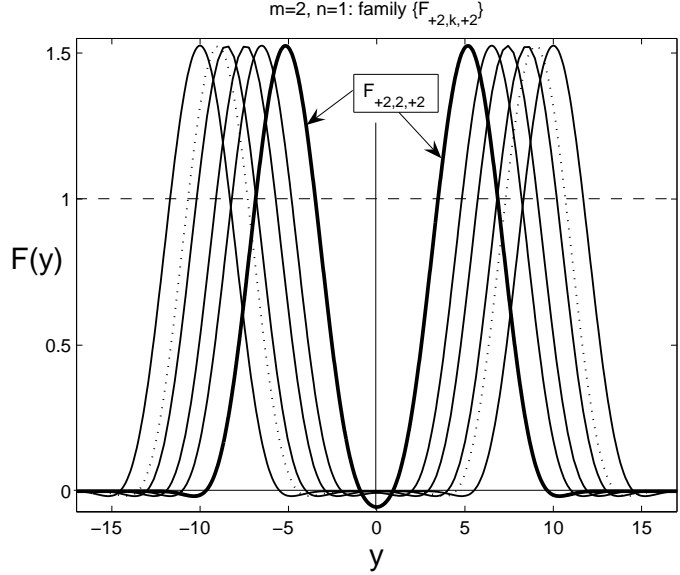}
 \vskip -.4cm
\caption{\rm\small First  patterns from the family
$\{F_{+2,k,+2}\}$ of the  $\{F_0,F_0\}$-interaction;  $n=1$.}
   \vskip -.3cm
 \label{G6}
\end{figure}

In Figure \ref{G6}, $n=1$, we show the first profiles from this
family denoted by $\{F_{+2,k,+2}\}$, where in each function
$F_{+2,k,+2}$ the multiindex $\s=\{+2,k,+2\}$ means, from left to
right, +2 intersections with the equilibrium +1, then next $k$
intersections with zero, and final +2 stands again for 2
intersections with +1. Later on, we will use such a multiindex
notation to classify other  patterns obtained.

\begin{figure}
\centering \subfigure[zeros: scale $10^{-2}$]{
\includegraphics[scale=0.52]{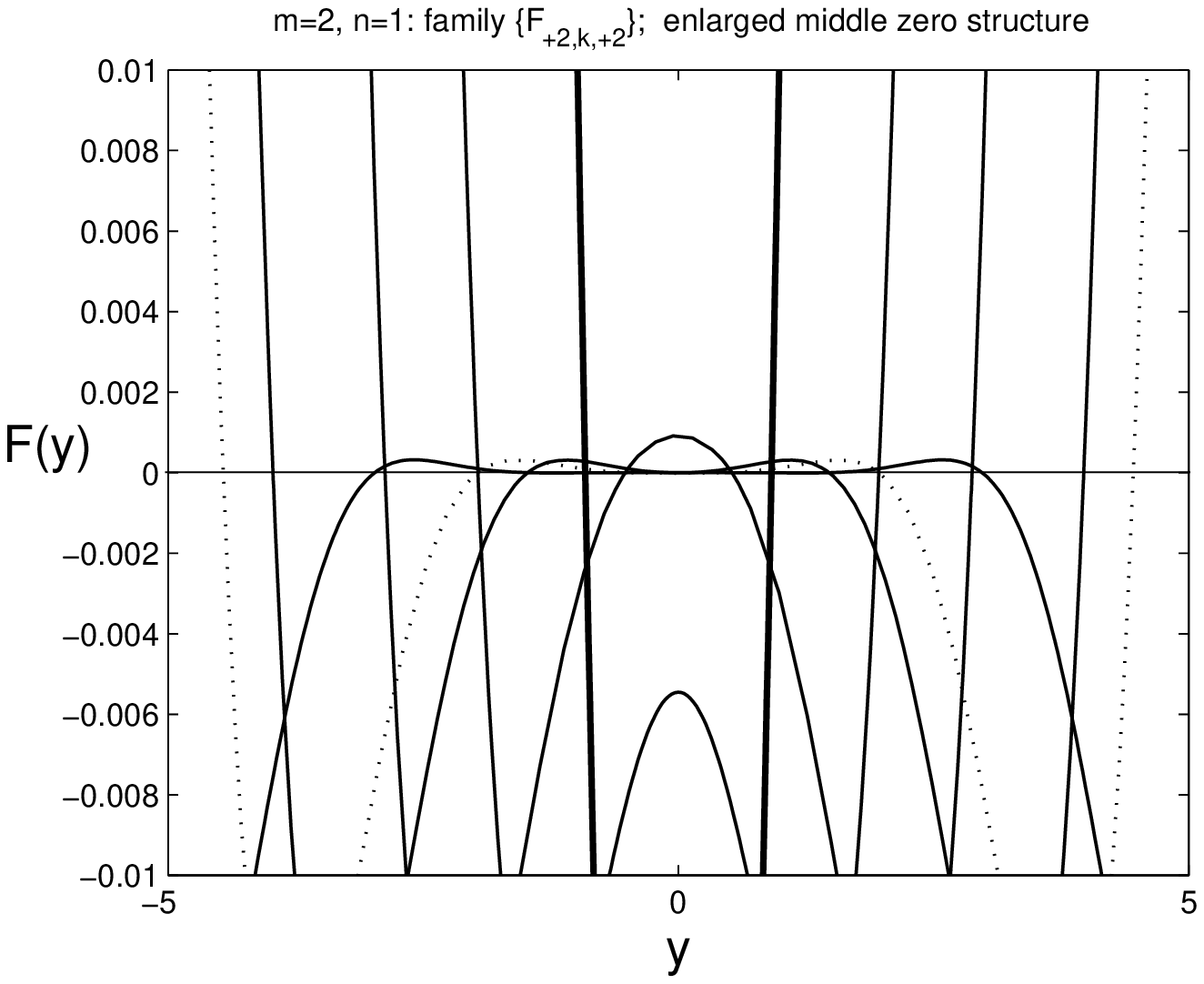}
} \subfigure[zeros: scale $10^{-4}$]{
\includegraphics[scale=0.52]{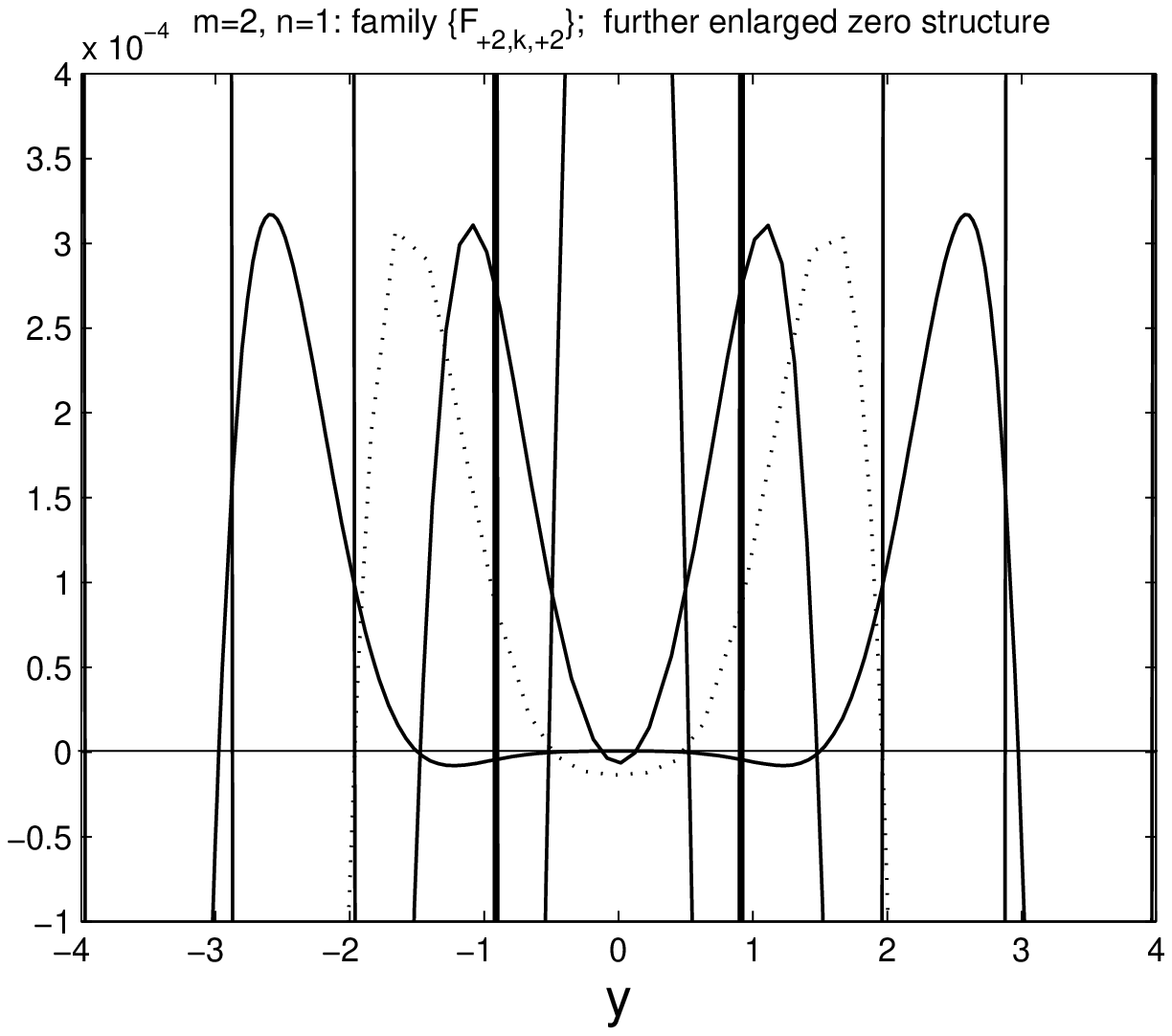}               
}
 \vskip -.2cm
\caption{\rm\small  Enlarged middle zero structure of the profiles
$F_{+2,k,+2}$ from Figure \ref{G6}.}
 \label{G61}
\end{figure}

\begin{figure}
 \centering
\includegraphics[scale=0.52]{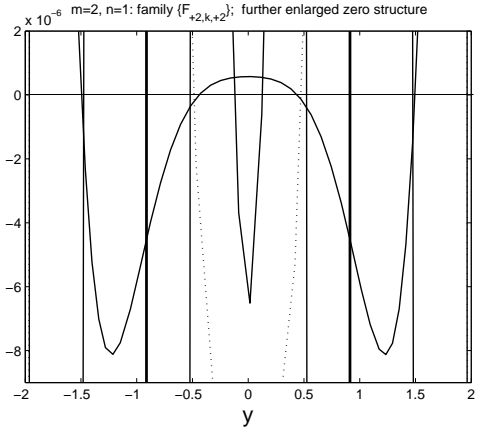}
 \vskip -.4cm
\caption{\rm\small Enlarged middle zero structure of the profiles
$F_{+2,6,+2}$ from Figure \ref{G6}.}
 \label{G62}
\end{figure}

In Figure \ref{G61}, we present the enlarged behaviour of zeros
explaining the structure of the interior layer of connection of
two profiles $\sim +F_0(y)$. In particular, (b) shows that there
exist  {\em two} profiles $F_{+2,6,+2}$, these are given by  the
dashed line and the previous one, both having two zeros on
$[-1,1]$. Therefore, the identification and the classification of
the profiles just by the successive number of intersections with
equilibria 0 and $\pm 1$ is not always acceptable (in view of a
non-homotopic nature of the problem), and some extra geometry of
curves near intersections should be taken into account.
 In fact, precisely this proves that a standard homotopy classification  of
 patterns is not consistent for such non-coercive and oscillatory equations.
 Anyway, whenever possible without  confusion,
we will continue use such a multiindex classification, though now
meaning that in general a profile $F_\s$ with a given multiindex
$\s$ may denote actually a {\em class} of profiles with the given
geometric characteristics.   Note that the last profile in Figure
\ref{G6} is indeed $F_{+2,6,+2}$, where the last two zeros are
seen in the scale $\sim 10^{-6}$ in  Figure \ref{G62}.
 Observe here a clear non-smoothness of two last profiles as a
 numerical discrete mesh phenomenon, which nevertheless does nor
 spoil at all this differential presentation.


In view of the oscillatory character of $F_0(y)$ at the
interfaces, we expect that the family $\{F_{+2,k,+2}\}$ is
countable, and such functions exist for any even $k=0,2,4,...\, $.
Then $k=+\infty$ corresponds to the non-interacting pair
  \beq
 \label{F0+}
 F_0(y+ y_0) + F_0(y-y_0), \quad \mbox{where} \,\,\,\,\,{\rm supp}\,
 F_0(y) = [-y_0,y_0].
  \eeq
Of course, there exist various triple $\{F_0,F_0,F_0\}$ and any
multiple interactions $\{F_0,...,F_0\}$ of $k$ single profiles,
with different distributions of zeros between any pair of
neighbours.

 \subsection{Countable family of  $\{-F_0,F_0\}$-interactions}

\begin{figure}
 \centering
\includegraphics[scale=0.65]{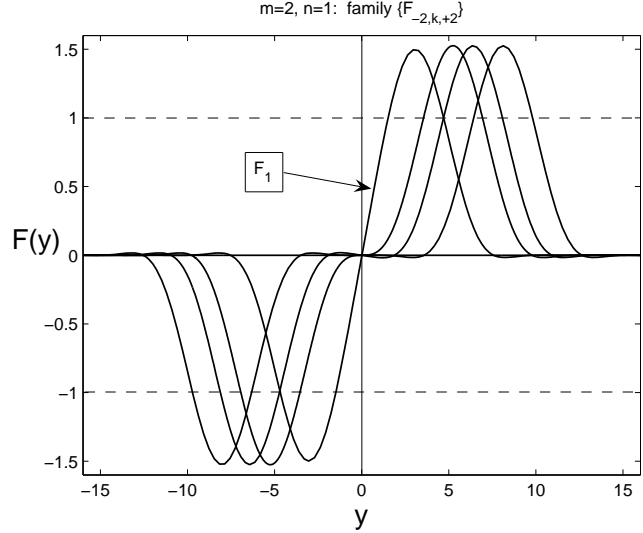}
 \vskip -.4cm
\caption{\rm\small First four patterns from the family
$\{F_{-2,k,+2}\}$ of the  $\{-F_0,F_0\}$-interactions;  $n=1$.}
   \vskip -.3cm
 \label{G7}
\end{figure}

\begin{figure}
\centering \subfigure[zeros: scale $10^{-2}$]{
\includegraphics[scale=0.52]{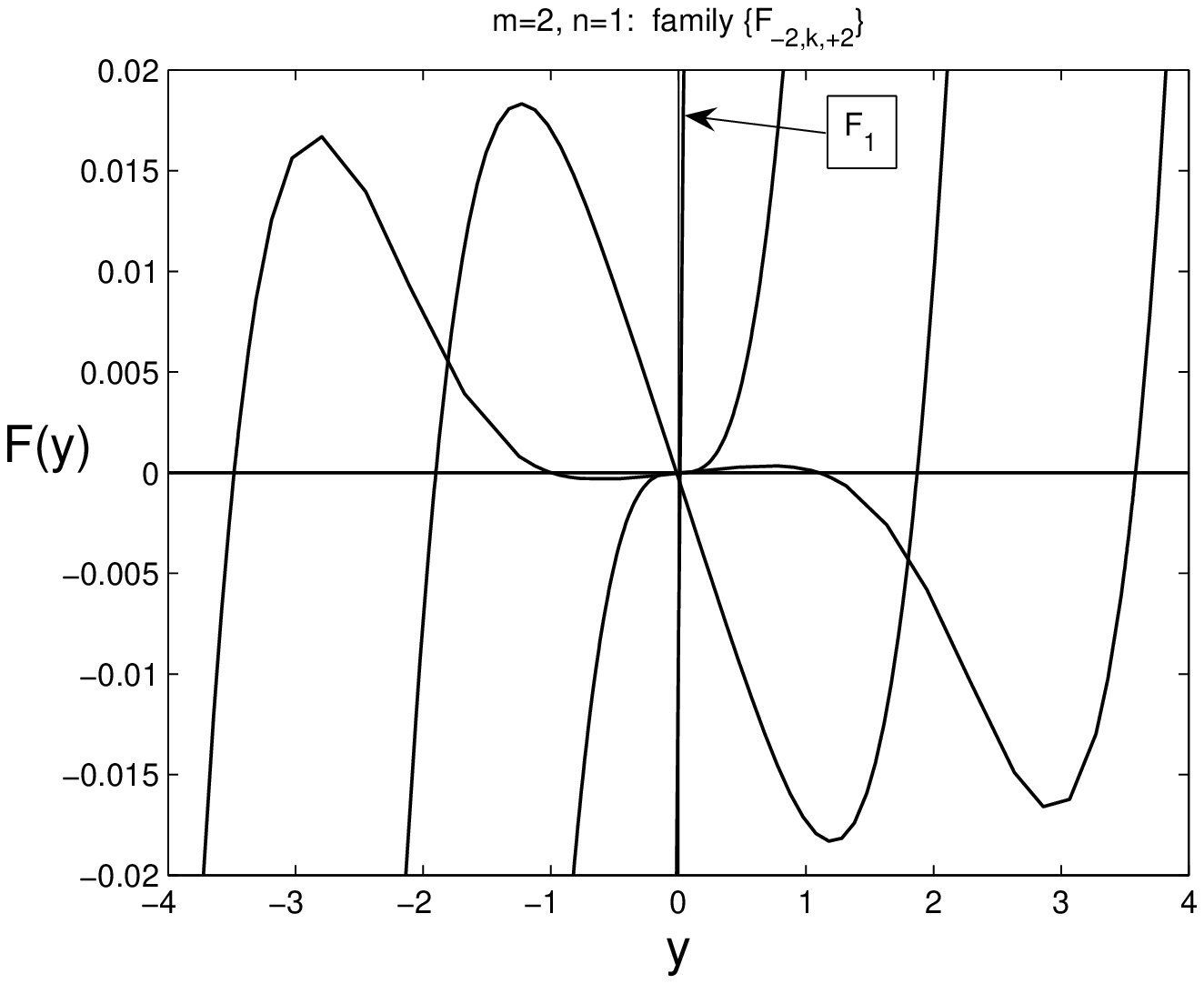}
} \subfigure[zeros: scale $10^{-4}$]{
\includegraphics[scale=0.52]{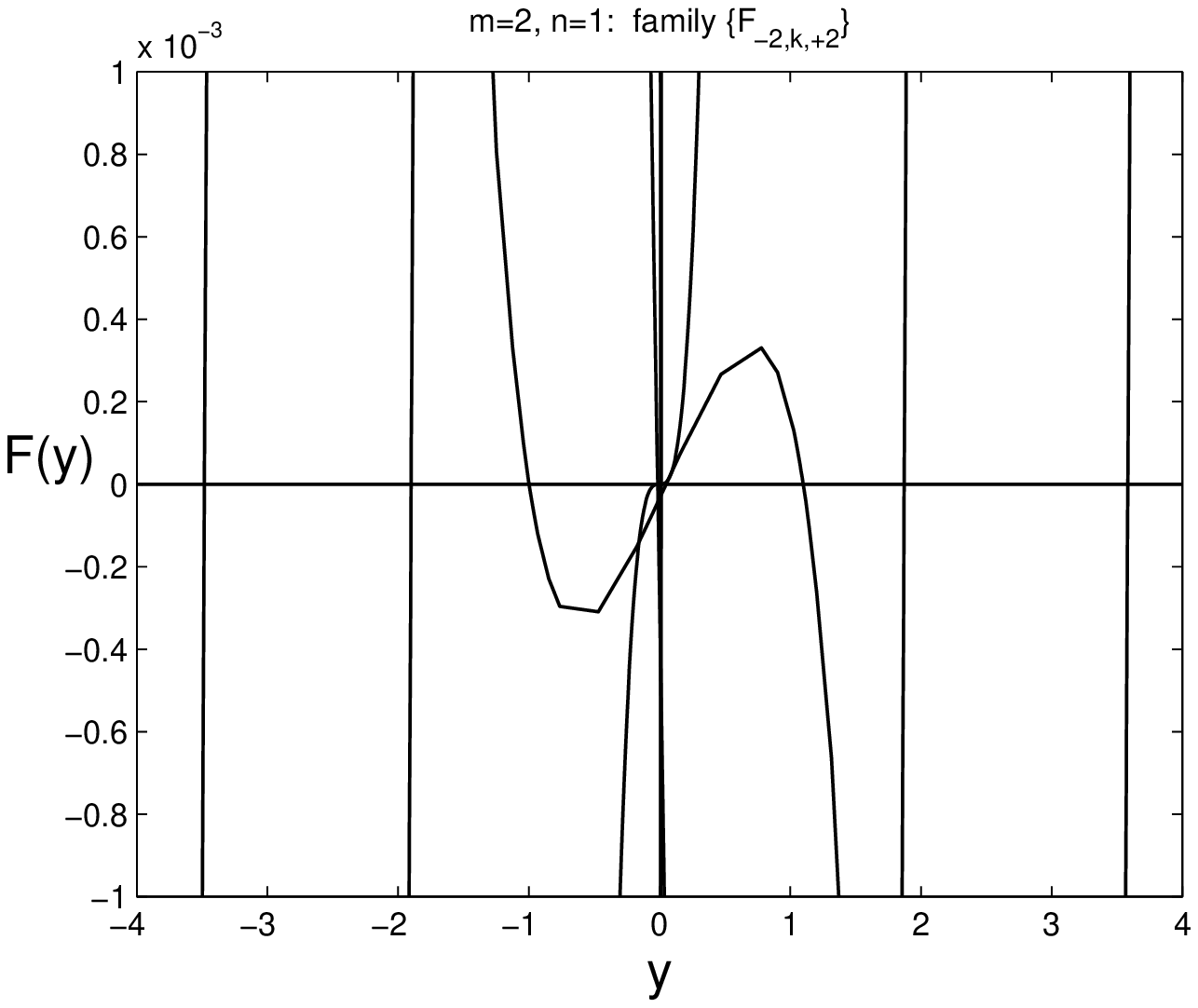}               
}
 \vskip -.2cm
\caption{\rm\small  Enlarged middle zero structure of the profiles
$F_{-2,k,+2}$ from Figure \ref{G7}.}
 \label{G71}
\end{figure}

 We now describe the interaction of $-F_0(y)$ with
  $F_0(y)$.
In Figure \ref{G7}, $n=1$, we show the first profiles from this
family denoted by $\{F_{-2,k,+2}\}$, where for  
 the
multiindex  $\s=\{-2,k,+2\}$,
  the first number $-2$ means 2
intersections with the equilibrium $-1$, etc. The zero structure
close to $y=0$ is presented in Figure \ref{G71}. It follows from
(b) that the first two profiles belong to the class
 $
 F_{-2,1,2},
 $
i.e., both have a single zero for $y \approx 0$.
  The last solution shown is $F_{-2,5,+2}$.
 Again, we expect that the family $\{F_{-2,k,+2}\}$ is countable,
and such functions exist for any odd $k=1,3,5,...$, and
$k=+\infty$ corresponds to the non-interacting pair
 \beq
 \label{F0-}
  \mbox{$
- F_0(y+ y_0) + F_0(y-y_0) \quad \bigl({\rm supp}\,
 F_0(y) = [-y_0,y_0]\bigr).
  $}
  \eeq
There exist families of an arbitrary number of interactions such
as $\{\pm F_0, \pm F_0,..., \pm F_0\}$ consisting of any $k \ge 2$
members.

\subsection{Periodic solutions in $\re$}

\begin{figure}
 \centering
\includegraphics[scale=0.65]{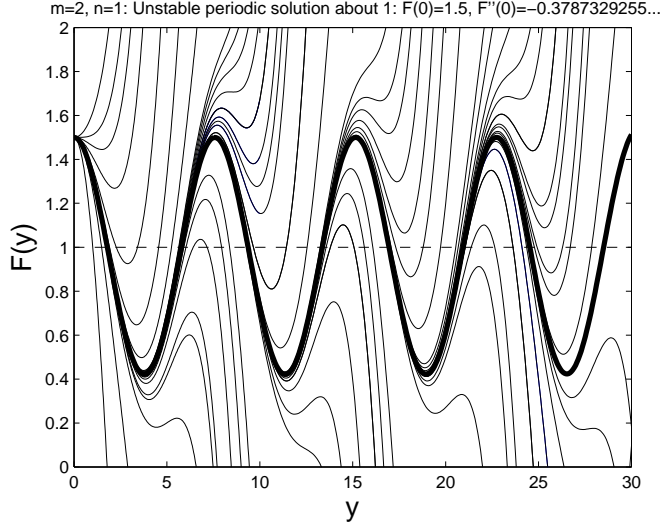}
 \vskip -.4cm
\caption{\rm\small An example of a periodic solution of the ODE
(\ref{4.3}) for $n=1$.}
   \vskip -.3cm
 \label{GP1}
\end{figure}

Before introducing new types of patterns, we need to describe
other non-compactly supported solutions in $\re$. As a variational
problem, equation (\ref{4.3}) admits an infinite number of
periodic solutions; see e.g. \cite[Ch.~8]{MitPoh}. In Figure
\ref{GP1} for $n=1$, we present a special  unstable periodic
solution  obtained by shooting from the origin with conditions
 $$
 F(0)=1.5, \quad F'(0)=F'''(0)=0, \quad F''(0)=-0.3787329255... \,
 .
 $$
 We will show next that precisely the periodic orbit $F_*(y)$ with
 \beq
 \label{**1}
 F_*(0) \approx 1.535...
 \eeq
  plays an important part in the
 construction of other families of compactly supported
patterns. Namely, all the variety of solutions of (\ref{4.3}) that
have oscillations about equilibria $\pm 1$ are  close to $\pm
F_*(y)$ there.

\subsection{Family $\{F_{+2k}\}$}

\begin{figure}
\centering \subfigure[ $F_{+2k}(y)$ ]{
\includegraphics[scale=0.52]{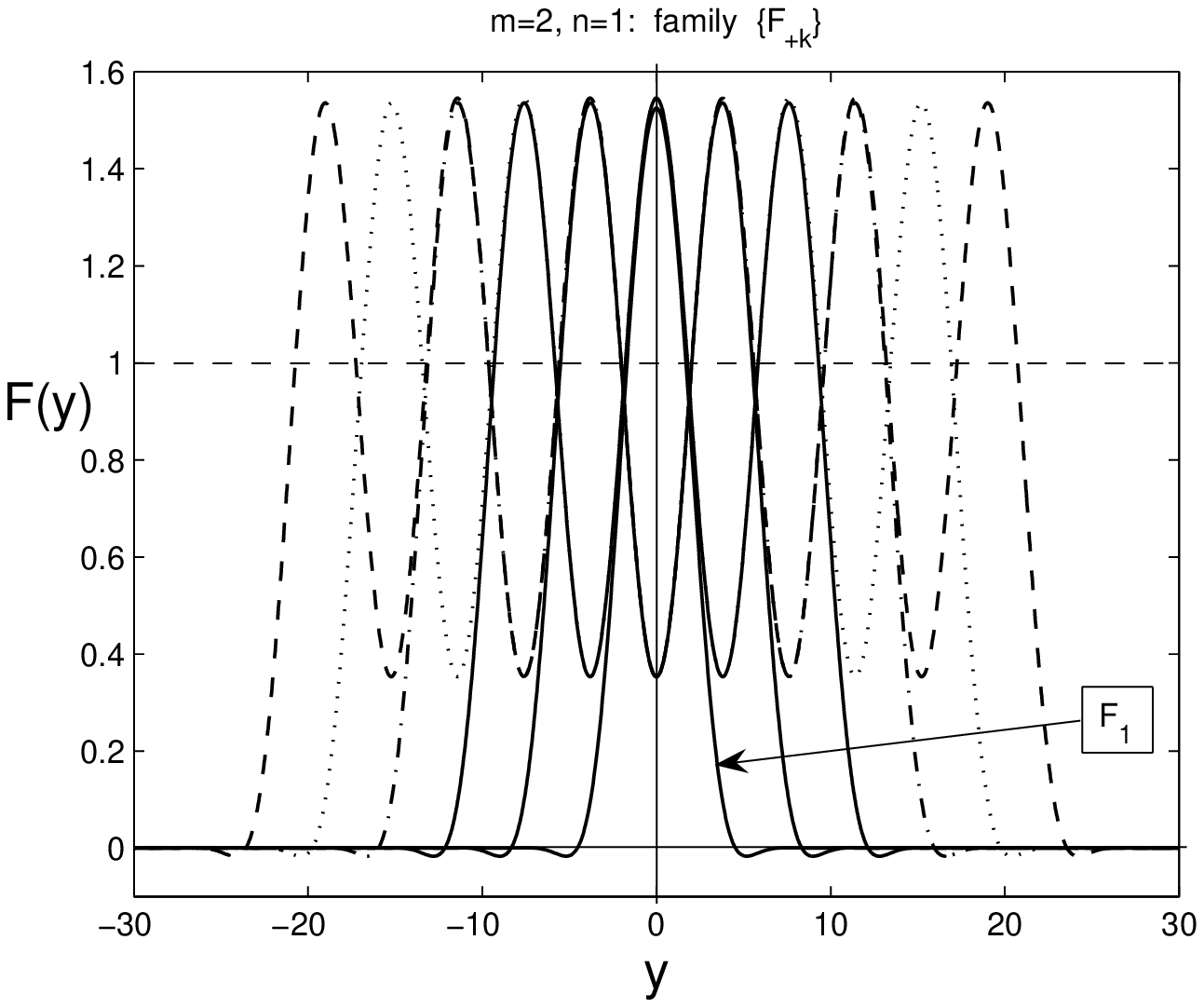}
} \subfigure[$F_{+k,l,-m,l,+k}$]{
\includegraphics[scale=0.52]{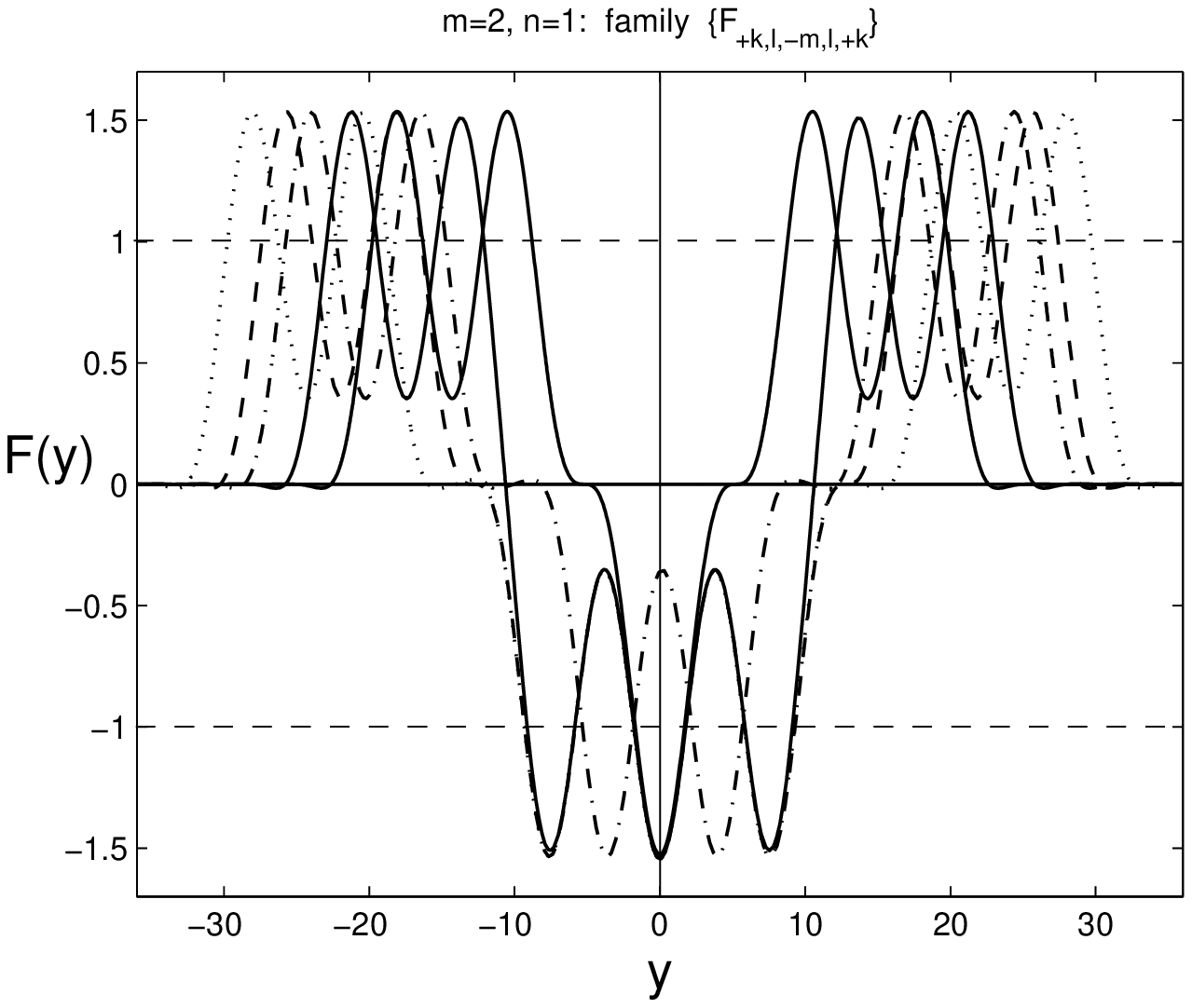}
}
 \vskip -.2cm
\caption{\rm\small Two families of solutions  of (\ref{4.3}) for
$n=1$; $F_{+2k}(y)$ (a) and $F_{+k,l,-m,l,+k}$ (b).}
 \label{G8}
\end{figure}

Such functions $F_{+2k}$ for $k \ge 1$ have $2k$ intersection with
the single equilibrium +1  only and have a clear ``almost"
periodic structure of oscillations about; see Figure \ref{G8}(a).
The number of intersections denoted by $+2k$ gives an extra Strum
index to such a pattern. In this notation,
 $
 F_{+2}=F_0.
  $

\subsection{More complicated patterns: towards  chaotic structures}

\begin{figure}
\centering \subfigure[$\s=\{+6,2,+2,2,+6\}$]{
\includegraphics[scale=0.52]{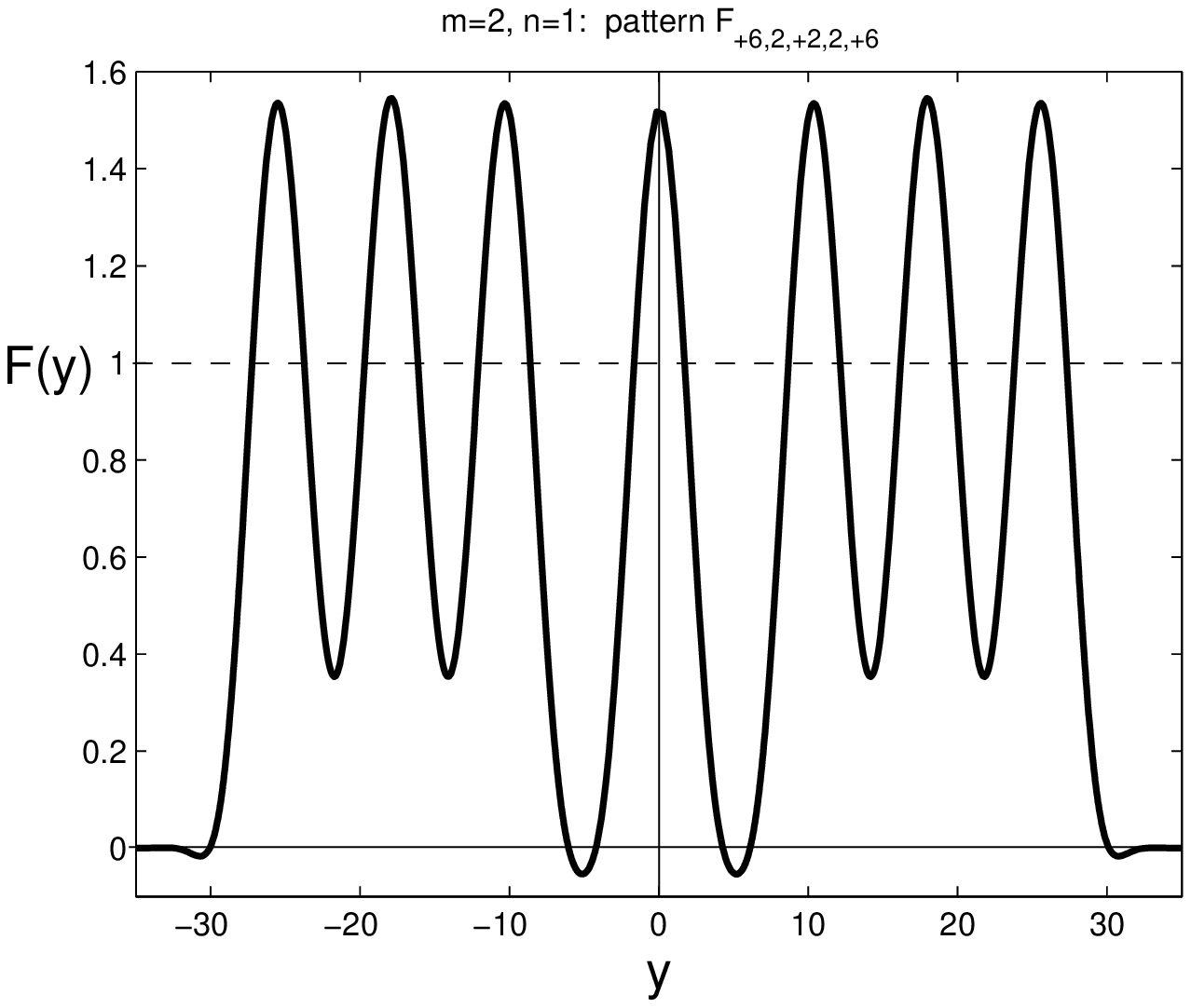}
} \centering \subfigure[$\s=\{+6,2,+4,1,-2,1+2\}$]{
\includegraphics[scale=0.52]{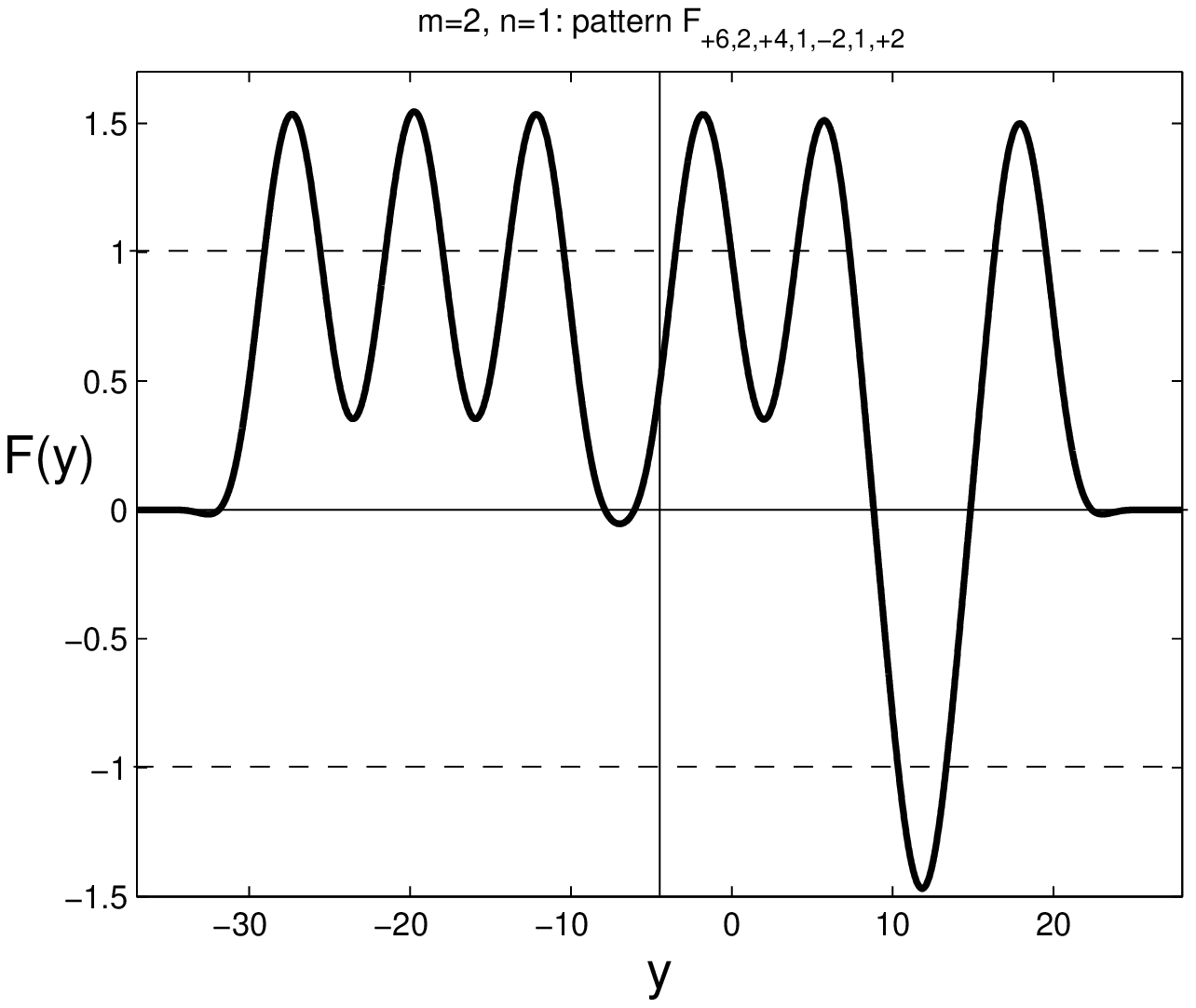}
} \centering \subfigure[$\s=\{+2,2,+4,2,+2,1,-4\}$]{
\includegraphics[scale=0.52]{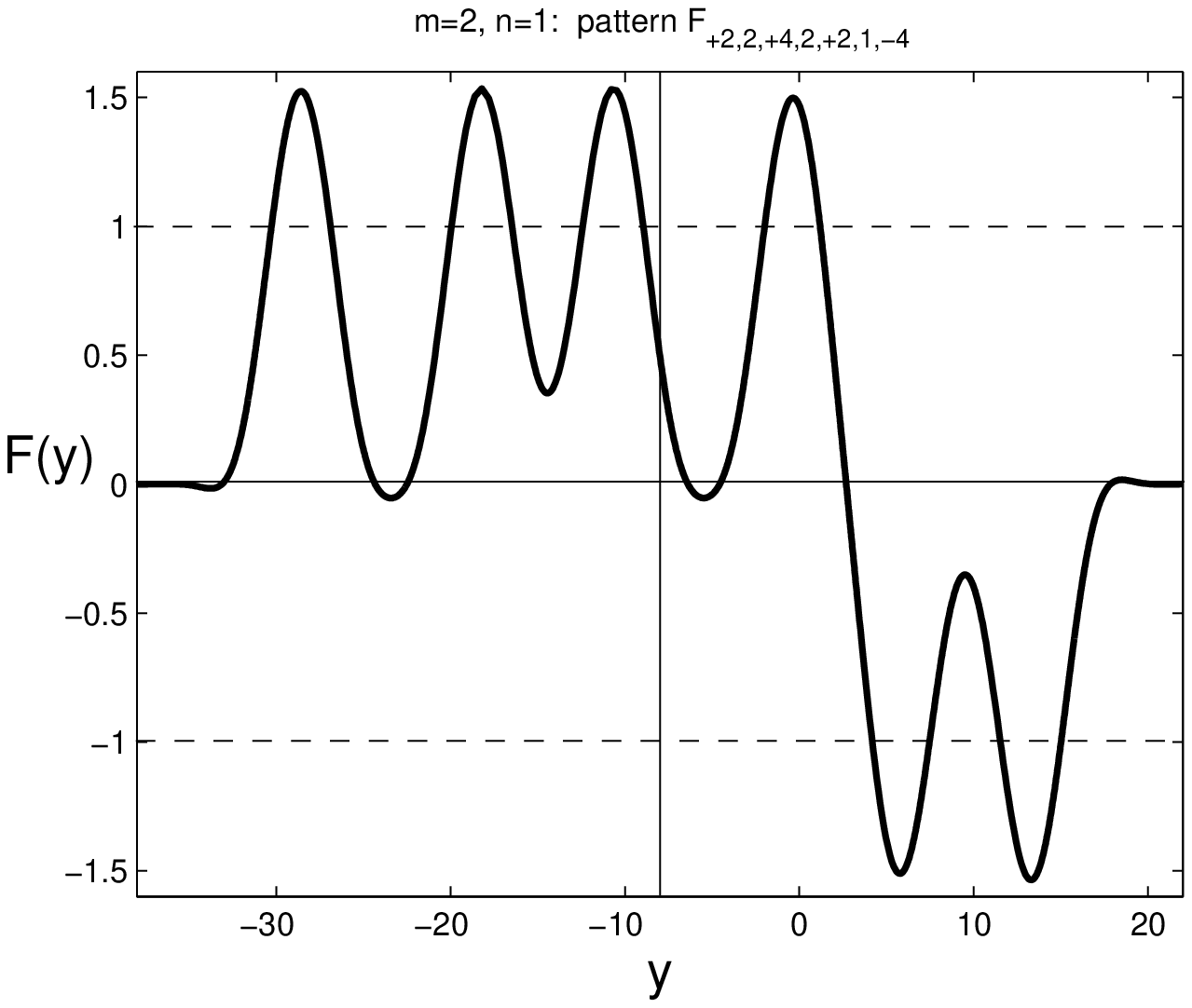}
} \centering \subfigure[$\s=\{+6,3,-4,2,-6\}$]{
\includegraphics[scale=0.52]{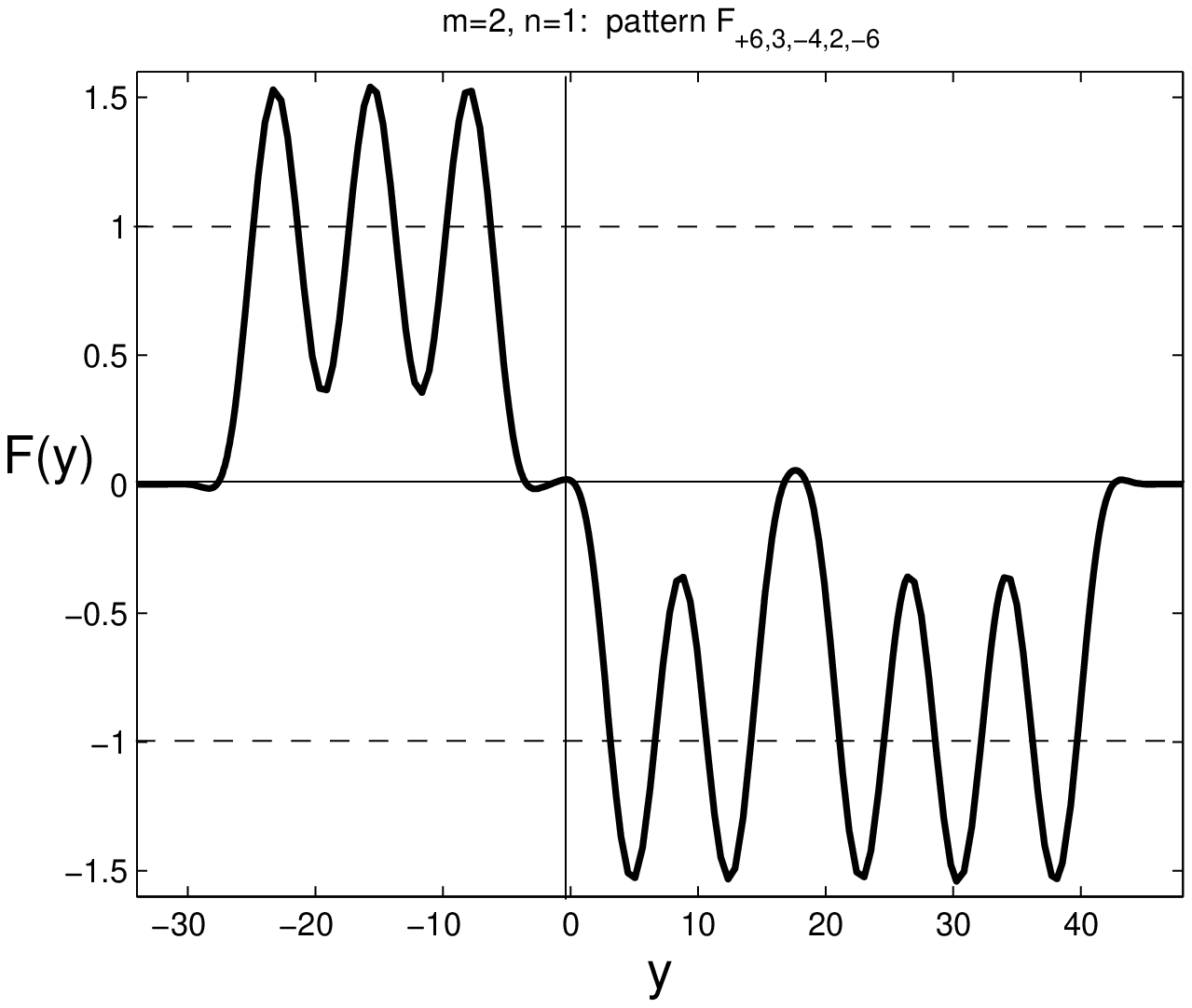}
}
 \vskip -.2cm
\caption{\rm\small  Various patterns for  (\ref{4.3}) for $n=1$.}
 \label{FC1}
\end{figure}

\begin{figure}
 \centering
\includegraphics[scale=0.65]{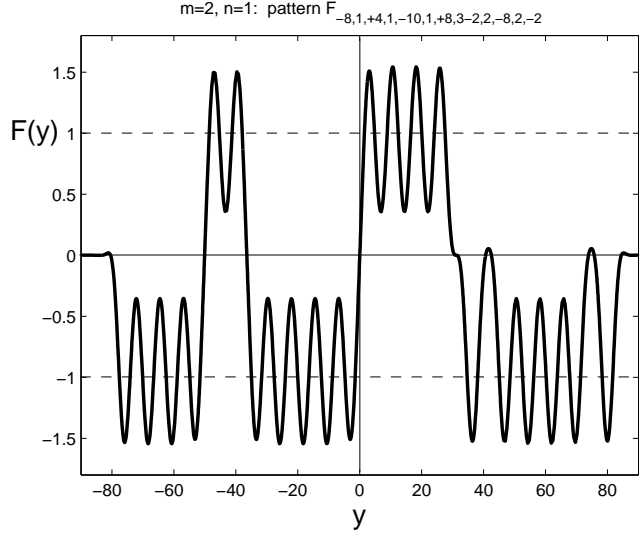}
 \vskip -.4cm
\caption{\rm\small A complicated pattern $F_{\s}(y)$ for
(\ref{4.3}) for $n=1$.}
   \vskip -.3cm
 \label{FC2}
\end{figure}

Using the above rather simple families of patterns, we claim that
a pattern (possibly, a class of patterns) with an arbitrary
multiindex of any length
 \beq
 \label{mm1}
 \s=\{\pm \s_1, \s_2, \pm \s_3, \s_4,..., \pm \s_l\}
  \eeq
  can be constructed.
 Figure \ref{G8}(b)  shows several profiles from the
family with the index
 $
  \s=\{+k,l,-m,l,+k\}.
  $
  In Figure \ref{FC1}, we show further four different patterns, while in
  Figure \ref{FC2}, a single most complicated pattern is
  presented, for which
   \beq
   \label{ch1}
   \s=\{-8,1,+4,1,-10,1,+8,1,3,-2,2,-8,2,2,-2\}.
    \eeq
   All computations are performed for $n=1$ as usual.
Actually, we claim that the multiindex (\ref{mm1}) can be
arbitrary and takes any finite part of any non-periodic fraction.
Actually, this means {\em chaotic features}
 of the whole family of solutions $\{F_\s\}$.
 These chaotic types of behaviour are known for other fourth-order ODEs with
coercive operators, \cite[p.~198]{PelTroy}.

\section{\underline{\bf Problem ``Numerics"}: patterns in 1D in higher-order cases, $m \ge 3$}
\label{Sectm4}

The main features of the pattern  classification by their
structure and computed critical values
 for $m=2$ in the previous section
    can
be extended to arbitrary $m \ge 3$ in the ODEs (\ref{S2}) for
$N=1$, so we perform this in less detail.

\begin{figure}
 \centering
\includegraphics[scale=0.80]{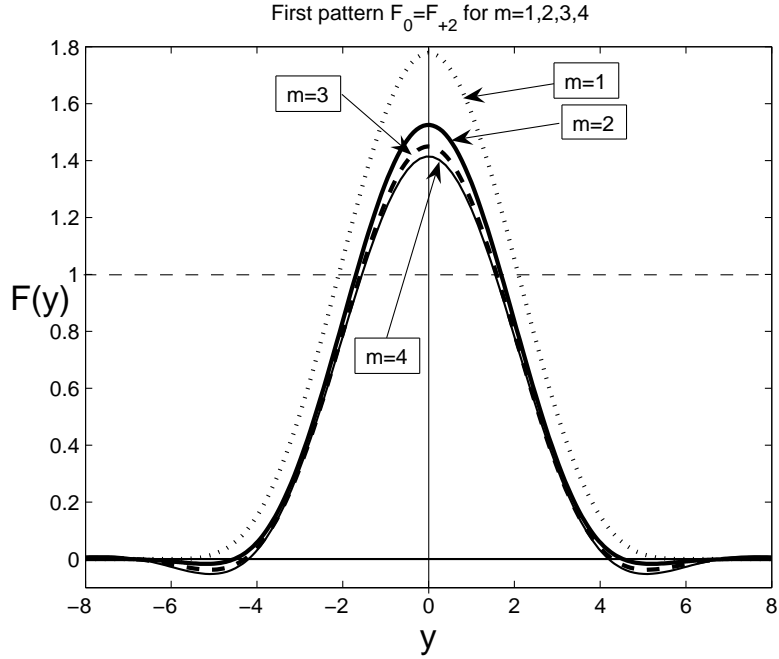}  
 \vskip -.4cm
\caption{\rm\small The first solution $F_0(y)$  of  (\ref{S2}),
$N=1$,
 $n=1$,
for  $m=1,2,3,4$.}
   \vskip -.3cm
 \label{Gmm1}
\end{figure}

 In Figure
\ref{Gmm1}, for the purpose of comparison, we show the first basic
pattern $F_0(y)$ for $n =1$ in four main cases $m=1$ (the only
non-negative profile by the Maximum Principle known from the 1970s
\cite{SZKM2}, \cite[Ch.~4]{SGKM}),
and  2, 3, 4. Next
 Figure \ref{Gmm2}  explains
the oscillatory properties close to the interface. It turns out
that, for $m=4$, the solutions are most oscillatory, so we it is
convenient to use this case  for illustrations.

\begin{figure}
\centering \subfigure[scale $10^{-2}$]{
\includegraphics[scale=0.52]{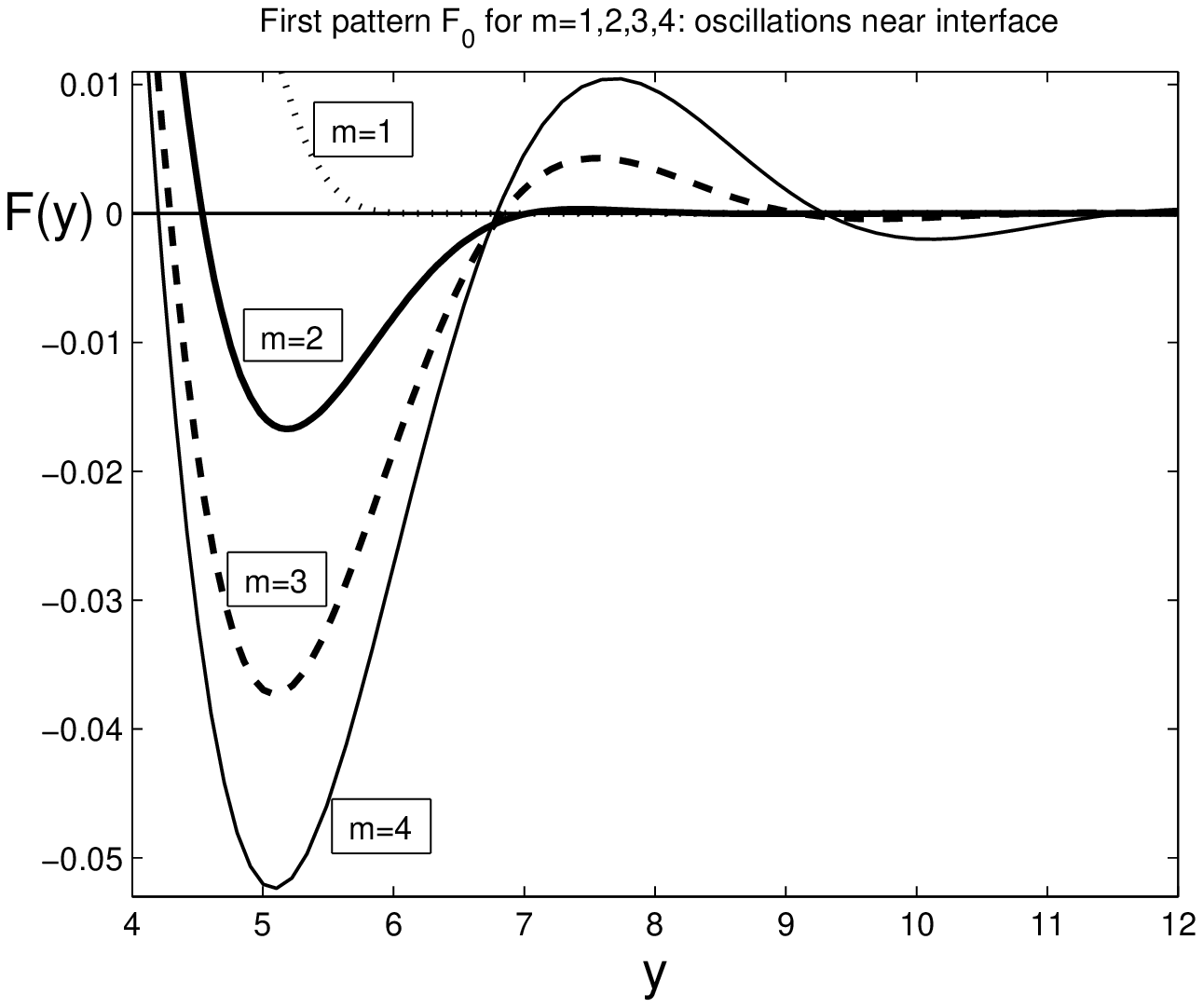}
} \subfigure[scale $10^{-3}$]{
\includegraphics[scale=0.52]{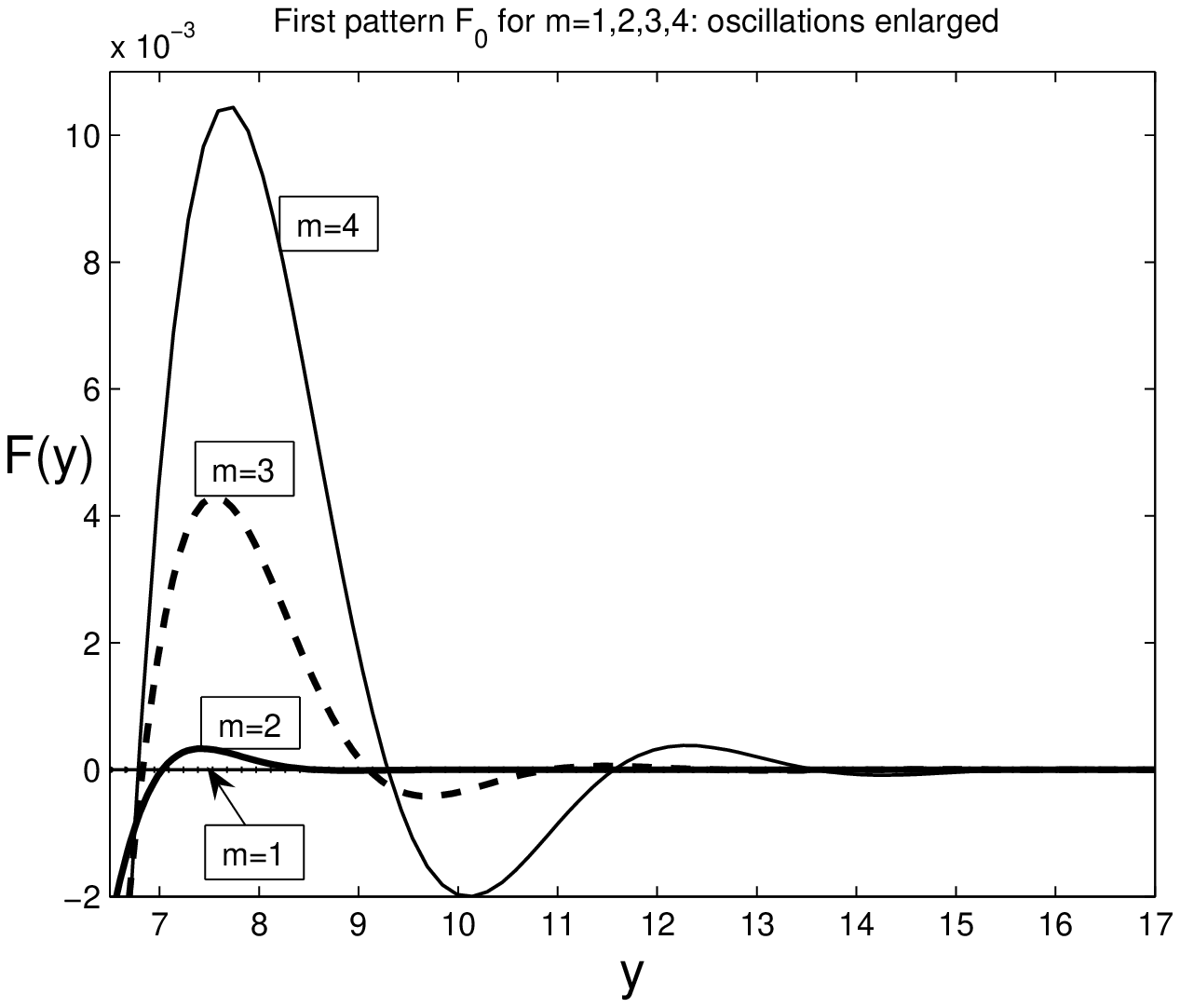}               
}
 \vskip -.2cm
\caption{\rm\small  Enlarged zero structure of the profile
$F_0(y)$ for $n=1$ from Figure \ref{Gmm1}; the linear scale.}
 \label{Gmm2}
\end{figure}



\begin{figure}
 \centering
\includegraphics[scale=0.75]{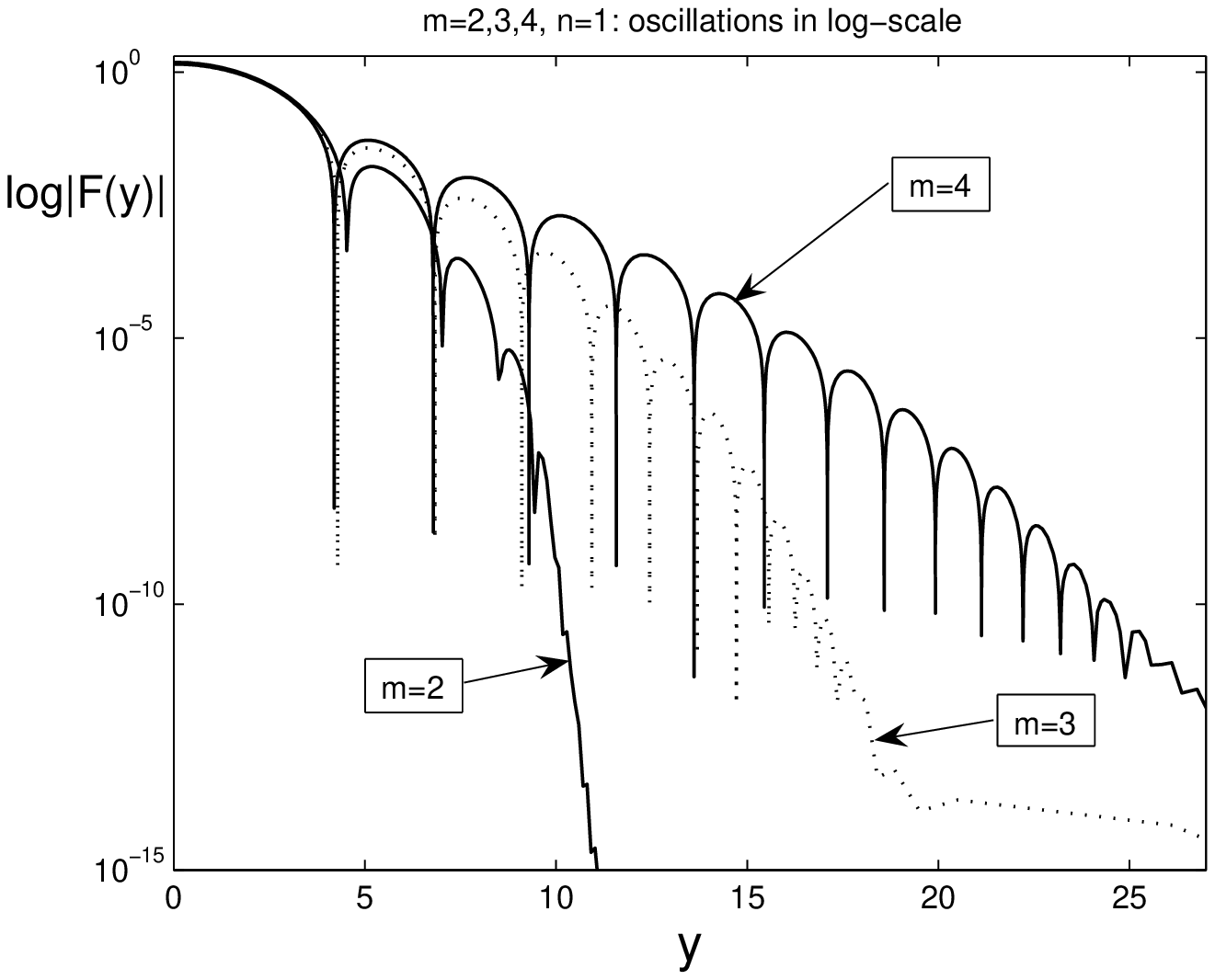}  
 \vskip -.4cm
\caption{\rm\small Behaviour of $F_0(y)$  for
 $n=1$,
for  $m=2,3,4$; the log-scale.}
   \vskip -.3cm
 \label{ZZ1}
\end{figure}

In the log-scale, the zero structure is shown in Figure \ref{ZZ1}
for $m=2,3$, and 4 ($n=1$). For $m=4$ and $m=3$, this makes it
possible to observe a dozen of oscillations that well correspond
to the oscillatory component analytic  formulae (\ref{2.2})
 close to
interfaces. For the less oscillatory case $m=2$, we observe 4
reliable oscillations up to $10^{-10}$, which is our best accuracy
achieved.

The basic countable family satisfying approximate Sturm's property
has the same topology as for $m=2$ in Section \ref{Sect4},
 and we do
not present such numerical illustrations.

\begin{figure}
 \centering
\includegraphics[scale=0.75]{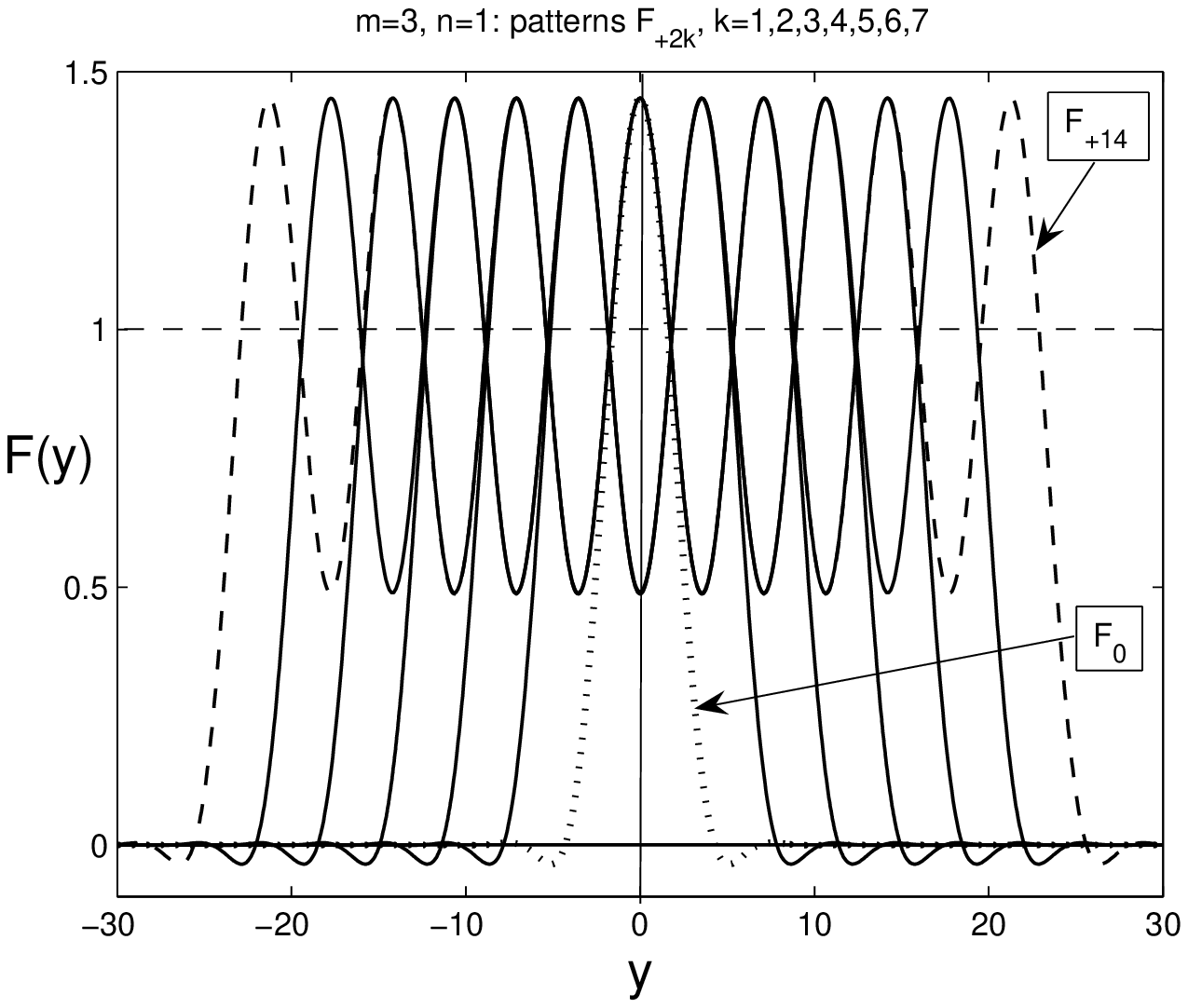}
 \vskip -.4cm
\caption{\rm\small The first seven patterns from the family
$\{F_{+2k}\}$;
 $m=3$ and
$n=1$.}
   \vskip -.3cm
 \label{Gmm6}
\end{figure}

\begin{figure}
 \centering
\includegraphics[scale=0.75]{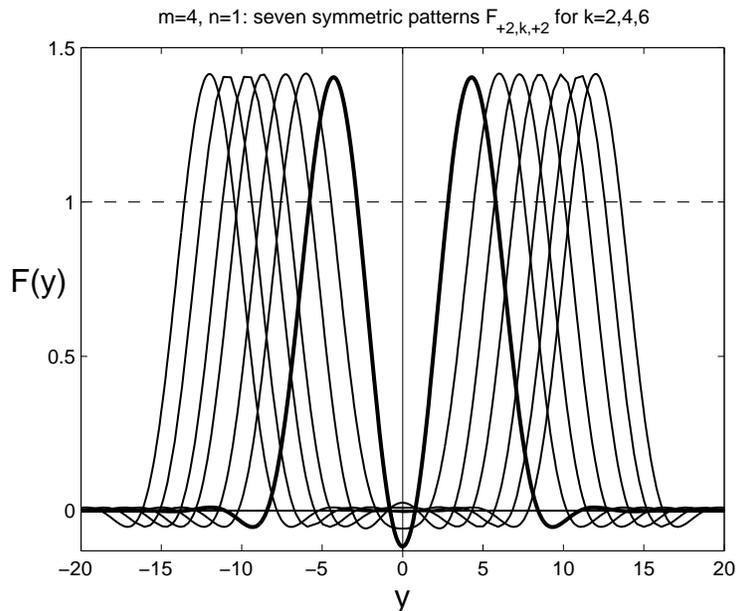}
 \vskip -.4cm
\caption{\rm\small The first  patterns from the family
$\{F_{+2,k,+2}\}$
 of the  $\{F_0,F_0\}$-interactions;
  $m=4$ and $n=1$.}
   \vskip -.3cm
 \label{Gmm68}
\end{figure}

\begin{figure}
\centering \subfigure[profiles]{
\includegraphics[scale=0.52]{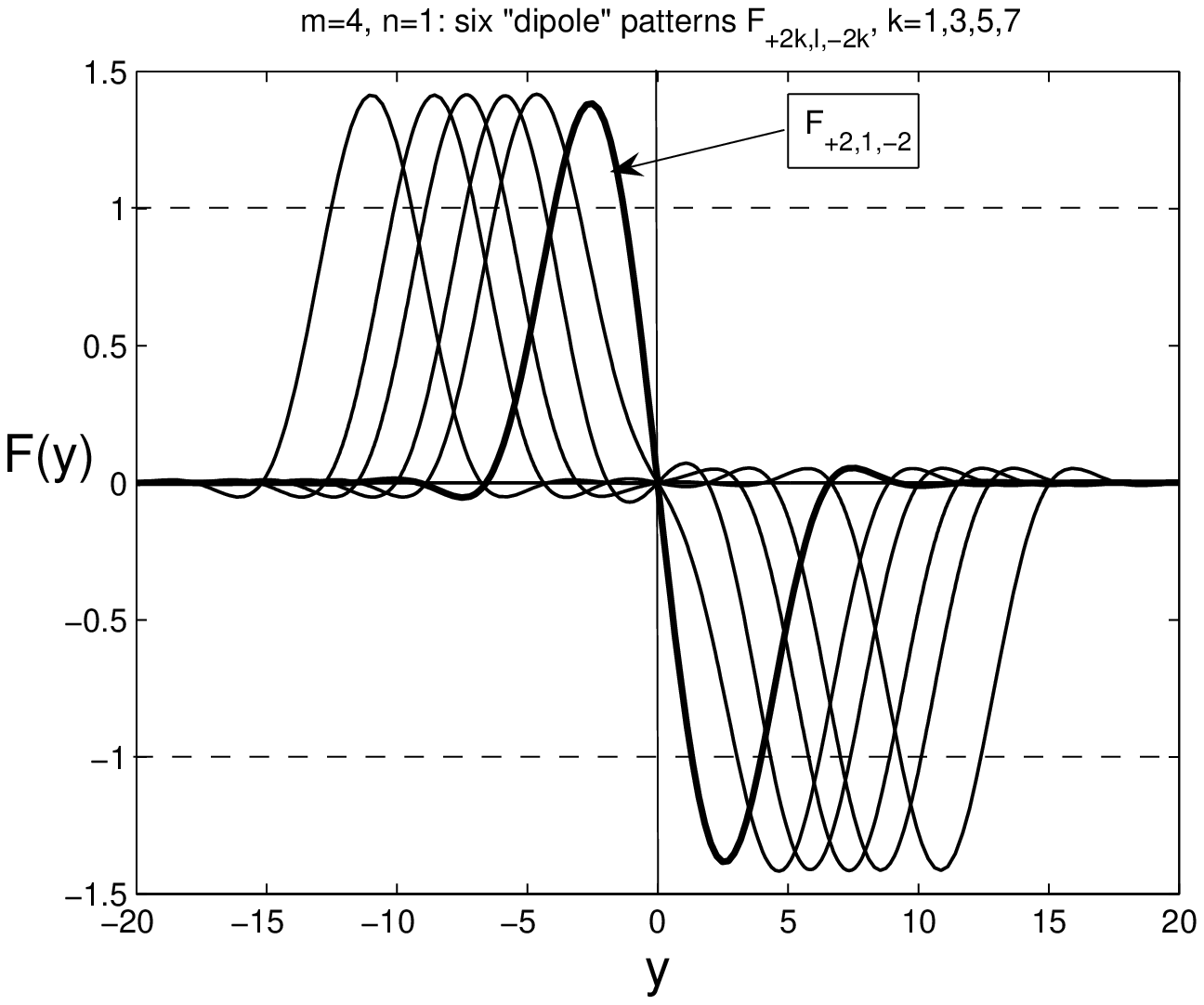}
} \subfigure[zero structure]{
\includegraphics[scale=0.52]{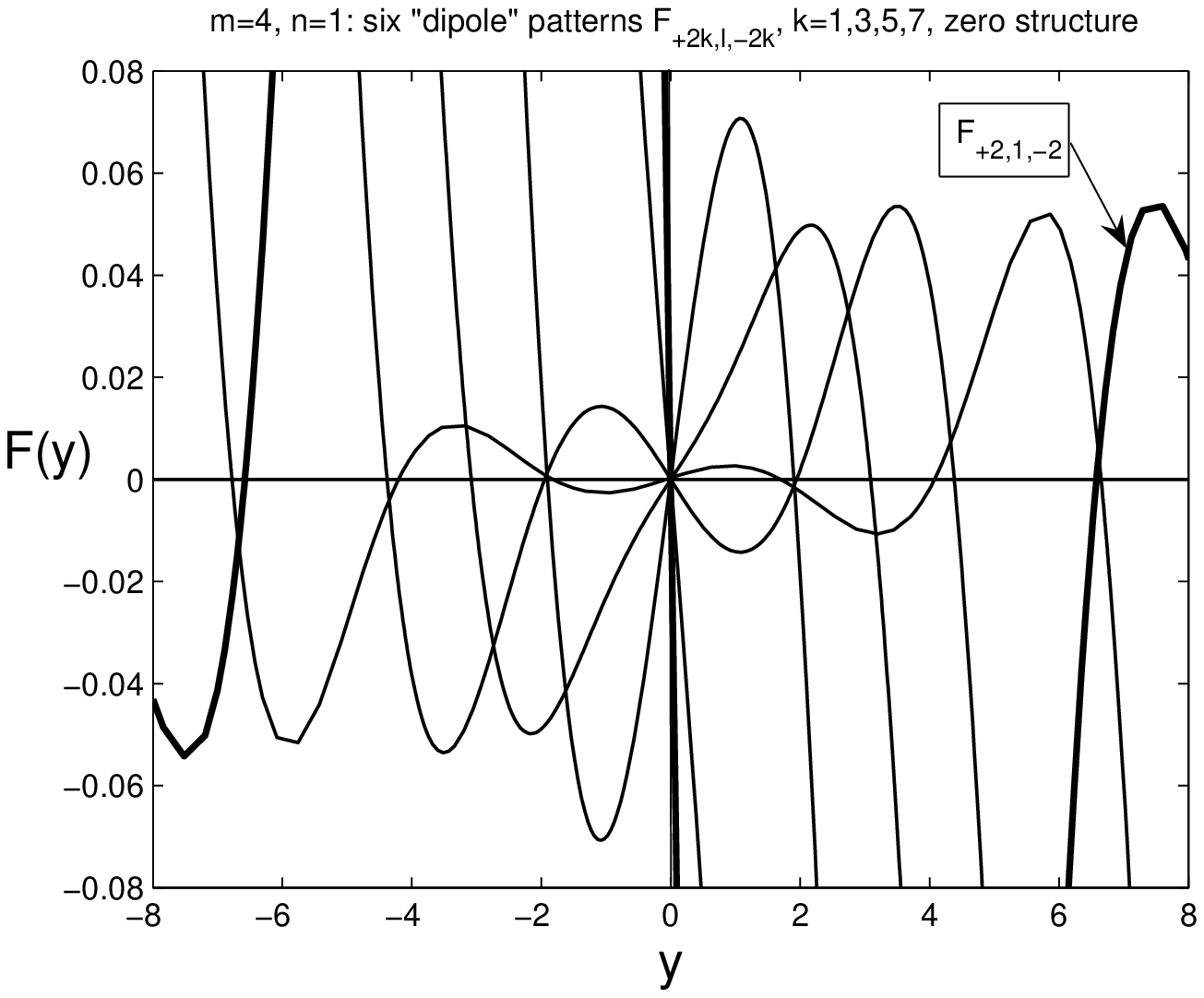}
}
 \vskip -.2cm
\caption{\rm\small The first patterns from the family
$\{F_{+2,k,-2}\}$ of the  $\{-F_0,F_0\}$-interactions, for $m=4$
and $n=1$: profiles (a), and zero structure (b).}
 \label{Gmm7}
\end{figure}

In Figure \ref{Gmm6} for $m=3$ and  $n=1$, we show the first
profiles from the family  $\{F_{+2k}\}$, while Figure \ref{Gmm68}
explains typical structures of $F_{+2,k,+2}\}$ for $m=4$, $n=1$.
In Figure \ref{Gmm7} for $m=4$ and $n=1$, we show the first
profiles from the family $\{F_{+2,k,-2}\}$.

\begin{figure}
 \centering
\includegraphics[scale=0.65]{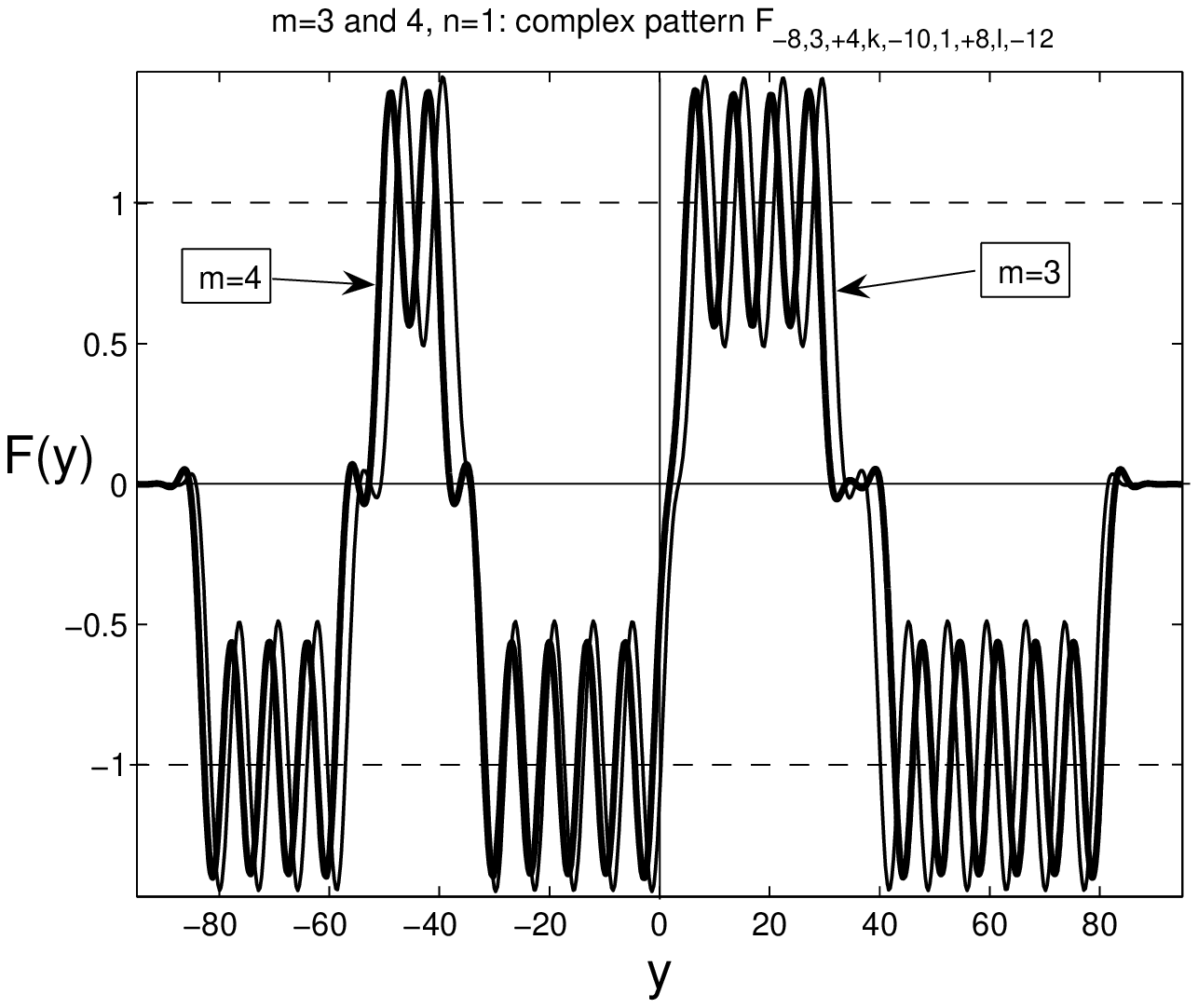}
 \vskip -.4cm
\caption{\rm\small A complicated pattern $F_{\s}(y)$  for $m=3,4$
and $n=1$.}
 \label{FCmm2}
\end{figure}

 Finally, in
  Figure \ref{FCmm2}, for comparison, we present  a complicated pattern
   for $m=3$ and $4$ (the bold line), $n=1$, with the index
   \beq
   \label{ll1}
   \s=\{-8,3,+4,k,-10,1,+8,l,-12\}.
    \eeq
Both numerical experiments were performed starting with the same
initial data. As a result, we obtain quite similar patterns, with
the only difference that, in (\ref{ll1}), $k=1$, $l=3$ for $m=3$,
and for more oscillatory case $m=4$, the number of zeros increase,
so now $k=3$ and  $l=5$.

Our study of other key aspects of these challenging elliptic
problems will be soon continued in \cite{GMPSobII}.



\end{document}